\documentclass[12pt]{article}
\usepackage{amsmath}
\usepackage{graphicx}
\usepackage{enumerate}
\usepackage{natbib}
\usepackage{url} % not crucial - just used below for the URL 
\usepackage{amsthm}
\usepackage{amssymb}

\usepackage[ruled,vlined]{algorithm2e}

\newcommand{\blind}{1}

\def\I{\mathbb I}
\def\d{\frac{1}{2}\,}
\def \tr{\mbox{tr}}

\def\R{\mathbb{R}}
\def\E{\mathbb{E}}
\def\eb{\overline{E}}
\newtheorem{proposition}{Proposition}
\newtheorem{remark}{Remark}
\newtheorem{lemma}{Lemma}
\newtheorem{example}{Example}
\newtheorem{theorem}{Theorem}
\newtheorem{corollary}{Corollary}
\begin{document}
% - - - - - - - - - - - - - - - - - - - - - - - - - - - - - - - - - - - - - - - - - - - - - - - - - - - - - - - - - - - -|

\def\spacingset#1{\renewcommand{\baselinestretch}%
{#1}\small\normalsize} \spacingset{1}

% - - - - - - - - - - - - - - - - - - - - - - - - - - - - - - - - - - - - - - - - - - - - - - - - - - - - - - - - - - - -|
% This work is dedicated to the memory of H\'el\`ene Massam

\if1\blind
{
  \title{\bf Accelerating Bayesian Structure Learning in Sparse Gaussian Graphical Models}
  \author{Reza Mohammadi 
  %\thanks{The authors gratefully acknowledge \textit{please remember to list all relevant funding sources in the unblinded version}}
  \hspace{.2cm} \\
    Department of Operation Management, University of Amsterdam \\
    and \\
    H\'el\`ene Massam \\ 
    Department of Mathematics and Statistics, York University \\
    and \\
    G\'erard Letac \\
    Laboratoire de Statistique et Probabilit\'es, Universit\'e Paul Sabatier \\
    }
  \maketitle
} \fi

\if0\blind
{
  \bigskip
  \bigskip
  \bigskip
  \begin{center}
    {\LARGE\bf Accelerating Bayesian Structure Learning in Sparse Gaussian Graphical Models}	
\end{center}
  \medskip
} \fi

\bigskip
% - - - - - - - - - - - - - - - - - - - - - - - - - - - - - - - - - - - - - - - - - - - - - - - - - - - - - - - - - - - -|
\begin{abstract} % around less than 200 words

Gaussian graphical models are relevant tools to learn conditional independence structure between variables. 
In this class of models, Bayesian structure learning is often done by search algorithms over the graph space. 
The conjugate prior for the precision matrix satisfying graphical constraints is the well-known $G$-Wishart. 
With this prior, the transition probabilities in the search algorithms necessitate evaluating the ratios of the prior normalizing constants of $G$-Wishart. 
In moderate to high-dimensions, this ratio is often approximated using sampling-based methods as computationally expensive updates in the search algorithm. 
Calculating this ratio so far has been a major computational bottleneck. %, barricades the capability of these models to the growing demand for higher-dimensional problems.
We overcome this issue by representing a search algorithm in which the ratio of normalizing constant is carried out by an explicit closed-form approximation. 
Using this approximation within our search algorithm yields significant improvement in the scalability of structure learning without sacrificing structure learning accuracy.
We study the conditions under which the approximation is valid. 
We also evaluate the efficacy of our method with simulation studies. 
%We represent a search algorithm in which the ratio of normalizing constant is carried out by our approximation, which yields significant improvement in the scalability of structure learning. 
We show that the new search algorithm with our approximation outperforms state-of-the-art methods in both computational efficiency and accuracy. 
The implementation of our work is available in the R package \textsf{BDgraph}. 

\end{abstract}
% - - - - - - - - - - - - - - - - - - - - - - - - - - - - - - - - - - - - - - - - - - - - - - - - - - - - - - - - - - - -|

\noindent% 3 to 6 keywords, that do not appear in the title
\vfill

\newpage
\spacingset{1.9} % DON'T change the spacing!
% - - - - - - - - - - - - - - - - - - - - - - - - - - - - - - - - - - - - - - - - - - - - - - - - - - - - - - - - - - - -|
\section{Introduction}
\label{sec:intro}

%Despite rapid developments in stochastic search algorithms, the practicality of Bayesian variable selection methods has continued to pose challenges.

% - >> GGMs ... Bayesian Structure learning ...  structure learning ... Bayesian inference like 2018 JRSS-B
Gaussian graphical models (GGM) have been widely used in many application areas for learning conditional independence structure among a (possibly large) collection of variables.
%Structure learning refers to the problem of estimating unknown graphs from the data and is often carried it out by estimating the sparsity pattern in the precision matrix.
Bayesian structure learning, for these models, while providing a natural and principled way for uncertainty quantification, often lag behind frequentist approaches \citep{friedman2008sparse} in terms of computational efficiency and scalability. 
%Hence a scalable approach to graphical models, with theoretical and computational safegurands, is critical to leveraging the advantages of posterior inference.
%Although substantial progress in computation for these graphical models has been made in recent years, scalability with dimensions remains a significant issue, 
%hindering the ability to adapt these models to the growing demand for higher-dimensional problems.
Despite significant developments of Bayesian structure learning methods in recent years, the scalability of these methods has continued to pose challenges regarding the growing demand for higher dimensions. 
%Despite significant developments in these models over the last few years, the scalability of Bayesian structure learning methods has continued to pose challenges regarding to the growing demand for higher dimensions. 

% - > > G-Wishart .... histroy like Hinne ... till ratio of IG
%To outline the issue with current Bayesian methods more clearly, consider 
An essential element of Bayesian structure learning in GGMs is the prior distribution on the precision matrix $K$ given the graph $G$ constraints. 
Most Bayesian methods use the so-called $G$-Wishart distribution, which is the conjugate prior \citep{roverato2002hyper}. 
For structure learning, more recent Bayesian methods, use versions of search algorithms over the graph space with the capability of jointly estimate graph structure and precision matrix, see \citet{hinne2014efficient, cheng2012hierarchical, lenkoski2013direct, dobra2011copula, dobra2011bayesian, wang2012efficient, mohammadi2015bayesianStructure}. % for arbitrary undirected graphs.
%Each of these search algorithms require substantial computational time due to difficulty with 
A computationally challenging step in these search algorithms is to estimate the ratio of prior normalizing constants for the $G$-Wishart distribution. % with respect to the current graph and with respect to a candidate graph obtained by adding/removing a single edge. 
This ratio, in general, is not available in closed form, except for specific cases, and typically needs to be evaluated using Monte Carlo based approaches.
Until recently, \citet{uhler2014exact} give the exact analytic expression of the normalizing constants of $G$-Wishart, which gave hope of direct evaluation of this ratio. 
The capability of applying this expression in the search algorithms need yet to investigate, since the expression is mathematically rather complex.
%As they point it out, one possibility would be to find more computationally efficient procedures than \citet[Theorem 3.3]{uhler2014exact}. % for computing the normalizing constant for particular classes of graphs.

% - - > state of art for ratio liks ... see refree 1 -JASA ....
To approximate the ratio of normalizing constant, \cite{wang12} introduces %a search algorithm which is based on %the exchange algorithm \citep{murray2006mcmc} and 
the double Metropolis-Hastings algorithm \citep{liang2010double}, by using on the block Gibbs sampler from $G$-Wishart. 
By using direct sampling form $G$-Wishart distribution \citep{lenkoski2013direct}, recently, \cite{hinne2014efficient, lenkoski2013direct} propose more efficient versions of the search algorithms that combine the concept behind the exchange algorithm \citep{murray2006mcmc} with trans-dimensional MCMC algorithm \citep{green2003trans}.
Likewise, \cite{mohammadi2015bayesianStructure} proposed a search algorithm over the graph space based on continuous-time birth-death processes, and following \cite{lenkoski2013direct} combined it with the exchange algorithm. % \citep{murray2006mcmc}. 
%We consider it here as state-of-the-art and we call it BDMCMC-DMH. 
%These algorithms provides significant computational improvements, since they eliminate the unstable MC approximation within the search algorithms.
These algorithms avoid to compute the ratio of normalizing constants by using the exchange algorithm;
Essentially, the ratio of normalizing constants is canceling out in the probabilities of jumping to the proposal graphs, by using exact samples from the $G$-Wishart distribution.
%These algorithm has proven to yield computational improvement  
%While these algorithms provide considerable computational gains compared to earlier approaches, they require exact sample from the $G$-Wishart distribution which is a computationally expensive updates within the search algorithm. 
While these algorithms have clear computational benefits compared to earlier approaches, they require exact samples from the $G$-Wishart distribution, which are computationally expensive updates within the search algorithm. 
We are going to illustrate it in more detail in Section \ref{subsec:review-IG}.

% Each of these approximations has proven to yield computational improvement in certain situations. % check Lenkoski 2013 123
%MCMC is routinely used for posterior computation, in which, the step of normalizing constant approximation often takes a substantial part of the run-time. 

\begin{figure} %[!ht]
\begin{center}
\includegraphics[width=6.0in]{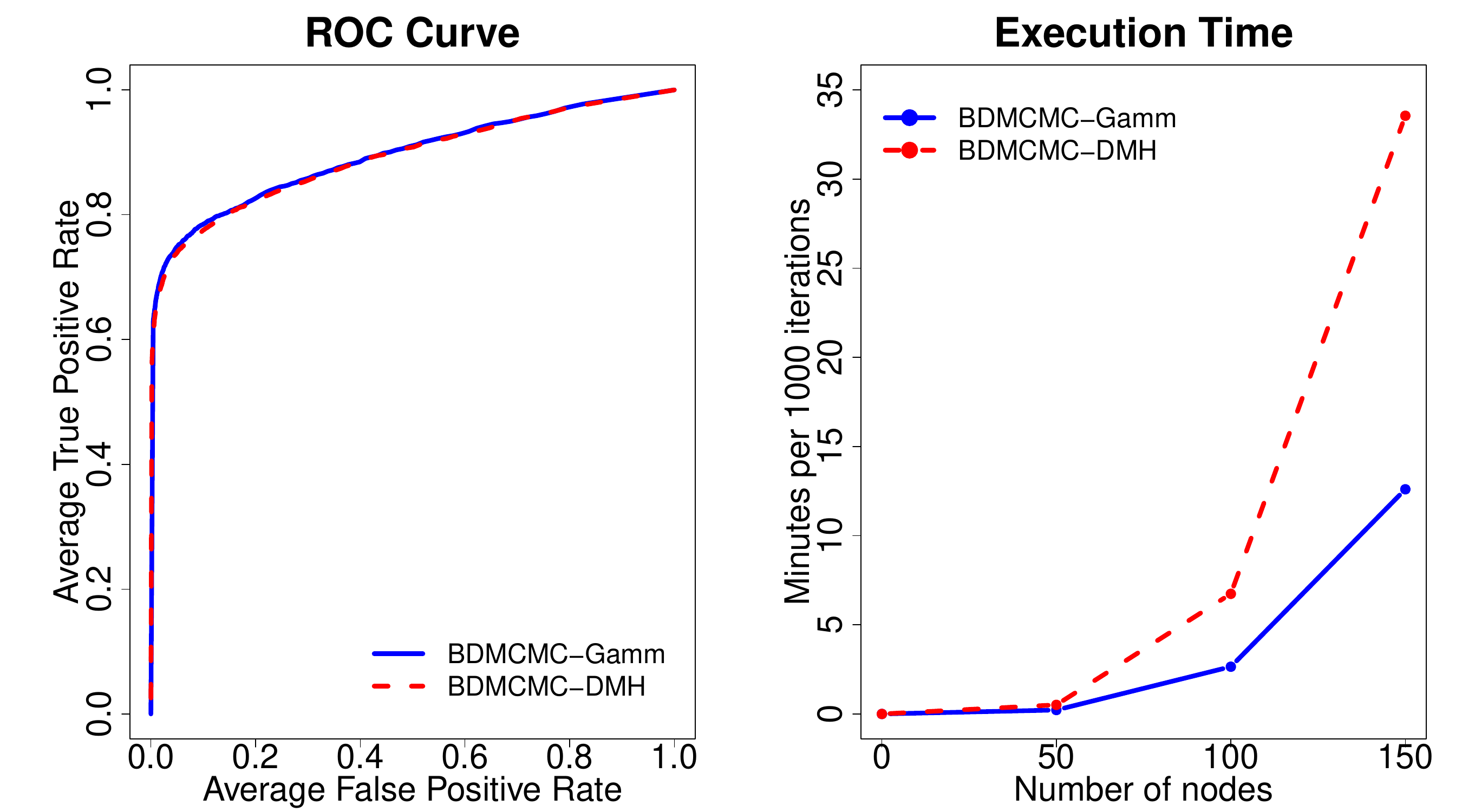}
\end{center}
\caption{ Plots from the simulation study in Section~\ref{sec:simulation} over $50$ replications where $G$ is a random graph with $150$ nodes and sample size $150$. (Left) ROC curve for the BDMCMC search algorithm with our approximation (BDMCMC-Gamm) and the BDMCMC search algorithm with exchange algorithm (BDMCMC-DMH), as state-of-the-art. (Right) execution time for both algorithms where time is per minutes for $1000$ iterations for different number of nodes ($p=50, 100, 150$). \label{fig:Rocplot_time} }
\end{figure}

% - - > > Our solution
%We aim to provide an explicit closed-form approximation to the ratio of the prior normalizing constant of $G$-Wishart, the use of which can lead to significant improvement in the scalability of the search algorithms. 
%For Bayesian structure learning, we hinge our simple approximation into the birth-death Markov Chain Monte Carlo (BDMCMC) search algorithm proposed by \cite{mohammadi2015bayesianStructure}.
We aim to introduce a search algorithm in which the ratio of normalizing constant is evaluated by an explicit closed-form approximation. 
For Bayesian structure learning, we first represent the birth-death Markov Chain Monte Carlo (BDMCMC) search algorithm proposed by \cite{mohammadi2015bayesianStructure}.
Then we provide an explicit closed-form approximation to the ratio of the prior normalizing constant of $G$-Wishart, the use of which leads to significant improvement in the scalability of the search algorithms. 
%Basically, our approximation can apply in any
%We remove the need for evaluating prior normalizing constants in a carefully designed MCMC search algorithm. 
%The aim of this paper is to prove that, for the purpose of Bayesian structure learning, the ratio of the normalizing constant can be approximated as ..... 
%This approximation is exact in some cases, see Remark \ref{remark:equality}.
To immediately illustrate the accuracy, in terms of structure learning, and the computational efficiency of our proposed approximation within the search algorithm, we represent here Figure \ref{fig:Rocplot_time} where $G$ has a random graph structure with $150$ nodes ($p=150$) and a sample size of $150$.
The left-hand side represents the receiver operating characteristic (ROC) plot for comparing the structure learning accuracy of the BDMCMC search algorithm done with our approximation and with the exchange algorithm. %when $G$ has a random graph structure, with $p=150$ and $n=150$. 
We see that our method (BDMCMC-Gamm, in blue) performs at least as well as the state-of-the-art (BDMCMC-DMH, in red). 
The right-hand side represents the execution time of both search algorithms. 
We see that for $p=150$, the execution time when using BDMCMC with our approximation is three times faster than when BDMCMC is done with the exchange algorithm.
More details are given in Section \ref{sec:simulation}. % and in the Supplementary file.

%At this point we recall that, in a remarkable breakthrough, a recent paper, \citet{uhler2014exact}, gives the exact analytic expression of the normalizing constants. 
%The capability of applying this expression for Bayesian structure learning in GGMs has yet to be investigated, since the expression is mathematically complex. 
%As they point it out, one possibility would be to find more computationally efficient procedures than \citet[Theorem 3.3]{uhler2014exact} for computing the normalizing constant for particular classes of graphs.
 
The outline of our paper is as follows. 
In Section \ref{sec:BGGMs}, we introduce background materials for Bayesian structure learning in GGMs. 
After presenting the birth-death MCMC search algorithm in Subsection \ref{subsec:BD}, we review the existing methods for approximating the ratio of normalizing constants in Subsection \ref{subsec:review-IG}, and then we introduce our approximation. 
In Section \ref{sec:ratio}, we provide the technical detail for proving the accuracy of the proposed approximation of the ratio of the normalizing constant. 
In Sections \ref{sec:disjoint} and \ref{sec:joint}, we represent our two main results, Theorems \ref{theorem:disjoint} and \ref{theorem:joint}. 

In Theorem \ref{theorem:disjoint}, we establish the approximation with explicit bounds in the particular case when all paths between the two nodes corresponding to the removed edge are disjoint (Figure \ref{fig:paths} left-side). 
%In that case, we derive the analytic expression of a lower bound for our approximation. 
% The key to proving this result is the fact that crucial random variables can be expressed as a linear product of independent normal variables in this case.
In Subsection \ref{subsec:sim-disjoint}, we verify the accuracy of the approximation by various collections of disjoint paths. 
We compute the theoretical boundary of our approximation as well as the value of the relative error following the Monte Carlo approach of \cite{atay2005monte}. 
We find that, while the theoretical boundary can be as high as $0.30$, the actual value of the relative error hardly goes above $0.10$ (see Figure \ref{fig:boxplot disjoint}). 

In Theorem \ref{theorem:joint}, we consider the general case where paths between the two nodes corresponding to the removed edge are not necessarily disjoint (Figure \ref{fig:paths} right-side). 
In that case, we prove that under a technical assumption, our approximation is accurate. 
%%that a quantity
% we cannot  give an explicit upper bound for the relative error. Instead, first we give an alternative expression \eqref{eq:newi1i2} for $I_1/I_2$. Then, with the help of this expression,  we show that, under the assumption that the variance of
 %% $b_1$, as defined later in \eqref{eq:b}, is small,  $I_1/I_2$ is accurately approximated by 1.
The question is then to know whether this assumption is realistic. 
In Subsection \ref{subsec:sim-joint}, for different types of graphs, we verify  numerically how well this assumption holds. 
%% with up to 30 nodes. 
%We also study numerically the accuracy of our approximation. 
We also evaluate the accuracy of our approximation by simulation. 
To do so, we compute the ratio of the normalizing constant in two ways: first following the Monte Carlo approximation of  \cite{atay2005monte} and, second, using our approximation. 
We see that in all cases, both approximations take the same range of values. 
They  are both reasonably accurate.
When the number of nodes is greater than $30$, due to the limitations of the Monte Carlo approximation in \cite{atay2005monte}, 
%the accuracy of the \eqref{eq:approx} cannot numerically verified directly.
one cannot numerically verify the accuracy of the approximation directly. 
So, in Section \ref{sec:simulation}, we verify it indirectly: we use both our approximation and the  exchange algorithm to compute the ratio in the BDMCMC search algorithm of \cite{mohammadi2015bayesianStructure} 
%% and judge the quality of model selection with the criteria of specificity, sensitivity and Matthews correlation coefficient (MCC) 
for graphs containing $50$, $100$, or $150$ nodes. 
We see that in all cases, our approximation yields results as good or slightly better than the exchange algorithm as a state-of-the-art. 

%We note here that, in both Theorems \ref{theorem:disjoint} and \ref{theorem:joint}, we show that $I_1/I_2$ is always less than 1. This means that
%\begin{equation*}
%\frac{\E(f_{E^{-e}}(\psi_{E^{-e}}))}{\E(f_E(\psi_E))}=
%\frac{I_1}{\E(f_E(\psi_E))}< \frac{I_2}{\E(f_E(\psi_E))},
%\end{equation*}
%so that, by using approximation \eqref{eq:approx}, we replace the Bayes factor
%$$\frac{P(G^{-e}\mid S)}{P(G\mid S)}=\frac{I_{G^{-e}}(\delta+n,I_p+S)}{I_G(\delta+n, I_p+S)}\frac{I_G(\delta, I_p)}{I_{G^{-e}}(\delta,I_p)}$$
%by a bigger quantity, which  means that by using \eqref{eq:approx}, we will slightly favor the sparser graph $G^{-e}$.

%Early works on G-Wishart distribution was largely computational in nature ...... and was predicated on three assumptions:
%Structure learning faces several challenges: the number of possible structures is super-exponential in the number of variables while the required sample size might be even exponential. 
%Therefore, finding good regularization techniques is very important in order to avoid overfitting and to achieve a better generalization performance.

% - - - - - - - - - - - - - - - - - - - - - - - - - - - - - - - - - - - - - - - - - - - - - - - - - - - - - - - - - - - -|

\section{Bayesian structure learning in GGMs}
\label{sec:BGGMs}

Graphical models \citep{lauritzen1996graphical} are powerful tools to express the conditional dependence structure among random variables by a graph in which each node corresponds to a random variable.
For the case of undirected graphs, also known as Markov random field \citep{rue2005gaussian}, an edge between two nodes determines the conditional dependence of the regarding variables.
Let $G=(V,E)$ be an undirected graph where $V$ contains $p$ nodes corresponding to the $p$ coordinates and the edges $E$ describe the conditional independence relationships among variables; We use the convention that if $(i,j) \in E$ then $i < j$.
Let $\eb$ be the complement of $E$ that indexes the missing edges of $G$. 

A Gaussian graphical model for the Gaussian random vector $\mathbf{X} = ( X_1, ..., X_p) \sim \mathcal{N}_p (\mu,K^{-1})$ is represented by an undirected graph $G=(V,E)$.
Variables $X_i$ and $X_j$ are independent given all the other variables if and only if there is no edge $(i,j)$ in $E$.
It is well-known \citep{lauritzen1996graphical} that in that case, the precision matrix $K=\Sigma^{-1}$ belongs to the cone $P_{G}$ of positive definite matrices with $K_{ij}=0$ whenever $(i,j) \in \eb$. 
In other words, the zero entries in the off-diagonal of the precision matrix correspond to conditional independencies in the graph;
It is an essential property of the precision matrix for model selection \citep{dempster1972covariance}.
One can then define the GGM for a given graph $G$ as the family of distributions
\[
{\cal N}_G = \{N(0, \Sigma) : K=\Sigma^{-1}\in P_G\}.
\]
The likelihood based on a random sample $\mathbf{X} = ( \mathbf{X}^{(1)}, ..., \mathbf{X}^{(n)} )^\top $ from ${\cal N}_G$ is
\begin{eqnarray*}
%\label{eq:likelihood}
P( \mathbf{X} | K,G ) \propto |K|^{n/2} \exp \left\{ \frac{-1}{2} \mbox{tr}( KS ) \right\},
\end{eqnarray*}
where $ S = \mathbf{X}^\top \mathbf{X}$.

In GGMs, for Bayesian structure learning, the standard conjugate prior for the precision matrix $K$ of the Gaussian distribution is the $G$-Wishart distribution \citep{roverato2002hyper, letac2007wishart}. % as the conjugate prior for multivariate data \eqref{likelihood}. 
The G-Wishart is the Wishart distribution restricted to the space of precision matrices with zero entries specified by a graph $G$. 
The G-Wishart density $W_G(b, \Omega)$ is %can be written as
\begin{equation*}
%\label{eq:gwishart}
P(K\mid G)=\frac{1}{I_G(\delta,\Omega)} |K|^{\frac{\delta-2}{2}}\exp \{\frac{-1}{2} \tr( K \Omega ) \}{\bf 1}_{P_G}(K),
\end{equation*}
where $|K|$ denotes the determinant of $K$ and the symmetric positive definite matrix $\Omega$ and the scalar $\delta>2$ are called, respectively, the scale and shape parameters. 
%As is the case most of the time, in the absence of prior information, the parameter $\Omega$ is taken to be the $p$-dimensional identity matrix $\I_p$; 
%Throughout, we set $\Omega = \I_p$.
The normalizing constant
\begin{equation}
\label{eq:I_G}
I_G(\delta,\Omega) = \int_{K \in P_G} |K|^{\frac{\delta-2}{2}} \exp \{\frac{-1}{2} \tr( K \Omega ) \} dK
\end{equation}
is of central interest to us. % here.
%The quantity $I_G(\delta,\Omega)$ is the normalizing constant and of central interest to us. % here.
For arbitrary graphs, the explicit formula for this normalizing constant is given in Proposition \ref{roposition:atay};
We return to the computations of this fact in Section \ref{sec:ratio}.

The joint posterior distribution of the graph $G$ and the precision matrix $K$ is given as
\begin{eqnarray}
\label{eq:posterior}
P(K, G \mid \mathbf{X}) & \propto & P(\mathbf{X} \mid K, G) ~ P( K \mid G) ~ P(G) \nonumber \\
    &\propto& P(G) ~ \frac{1}{ I_G(\delta,\Omega) } |K|^{\frac{\delta+n-2}{2}} \exp \{ \frac{-1}{2} \tr( K (\Omega+S) ) \},
\end{eqnarray}
where %$\Omega^* = S + \Omega$ and 
$P(G)$ is the prior distribution of the graph $G$, which here we consider a uniform distribution over all graphs with fixed $p$ nodes, as a non-informative prior; 
For other options, see \cite{dobra2011bayesian, hinoveanu2018loss, mohammadi2015bayesianStructure}.
    
% - - - - - - - - - - - - - - - - - - - - - - - - - - - - - - - - - - - - - - - - - - - - - - - - - - - - - - - - - - - -|

\subsection{Structure learning via birth-death MCMC algorithm}
\label{subsec:BD}

\begin{figure} %[t]
\begin{center}
 \includegraphics[width=6.3in]{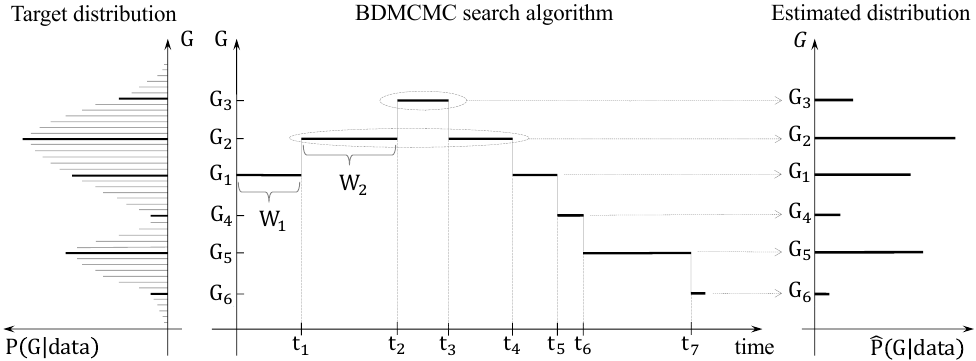} 
\end{center}
    \caption{ The left panel represents the target posterior distribution of the graphs. The middle panel shows the sampling scheme of BDMCMC  algorithms in which $\left\{ W_1, W_2,... \right\}$ stand for waiting times and $\left\{ t_1, t_2,... \right\}$ stand for jumping times of the BDMCMC algorithm. The right panel shows the estimated posterior distribution of the graphs based on the BDMCMC sampler which are the proportional to the total waiting times of the visited graphs. \label{fig:bdmcmc} }
\end{figure}

Bayesian structure learning in GGMs which revolves around the joint posterior distribution of the precision matrix and graph \eqref{eq:posterior} requires carefully designed MCMC search algorithms over the graph space. 
A common way to explore the graph space is by using a search algorithm known as reversible jump MCMC (RJMCMC) \citep{green1995reversible} which is based on a \emph{discrete}-time Markov chain. 
%These types of algorithms often suffer from low acceptance rates particularly when the model space is large, as the case for the graph models. 
%The issue is surprisingly common since the graph space is enormous and proposals with low probabilities are frequent. 
These kinds of algorithms often suffer from low acceptance rates since the graph space is enormous and proposals with low probabilities are frequent. 
\cite{mohammadi2015bayesianStructure} addressed this issue by developing a \emph{continuous}-time Markov chain process---or a BDMCMC search algorithm---as an alternative to RJMCMC. 
The BDMCMC search algorithm explores the graph space by either jumping to a larger dimension (birth) or lower dimension (death). 
The birth/death events are modeled as independent Poisson processes, thus the time between two successes events is exponentially distributed. 
%The events occur in continuous time and their rates determine the stationary distribution of the process; 
The stationary distribution of the process is determined by the rates of the birth and death events that occur in continuous time;
See Figure~\ref{fig:bdmcmc} for a graphical representation of birth and death events from a given graph. 
%Unlike the RJMCMC, in the BDMCMC search algorithm the moves between models are \emph{always} accepted making the algorithm more efficient. 

In the birth and death process, given the current state $(G, K)$, each edge is added/deleted independently of the rest as a Poisson process with birth/death rate $R_e( G, K )$ for each $e \in \{E \cup \eb \}$. 
Since birth and death events are independent Poisson processes, the time between two consecutive events has an exponential distribution with mean
\begin{equation}
\label{eq:waiting time}
W( G, K) = \frac{1}{\sum R_e( G, K ) }
\end{equation} 
which is the \textit{waiting time}. 
The waiting times capture all the possible moves of each step of the BDMCMC search algorithm. %, since they are calculated based on all the possible birth and death rates of the current state $(G, K)$. 
%Therefore, if the waiting time from $(G, K)$ is small then the process tends to stay shorter in the current state while if the waiting time is large then the process tends to stay longer.
Essentially, the birth-death process tends to stay shorter in the current state for a small waiting time, while the process tends to stay longer for a large waiting time. 
%Note that the waiting times are calculated based on all the possible birth and death moves from the current state $(G, K)$ to a new state which would be a graph with one more or less edge regarding to the birth/death rates. 
%Therefore, the waiting times essentially capture all the possible moves of each step of the BDMCMC search algorithm. 
%If the waiting time from $(G, K)$ is large then the process tends to stay longer in the current state while if the waiting time is small, the process tends to transition away from the current state. 
The birth and death probabilities involved are 
\begin{equation}
\label{eq:prop rate}
P(\mbox{birth/death of edge  } e \in \{ E \cup \eb \} ) \propto R_e( G, K ).
\end{equation}
%\begin{equation}
%\label{prop death}
%P(\mbox{death of edge  } e \in E ) \propto D_e( G, K ).
%\end{equation}
%For the birth of edge $e \in \eb$ we take $G^*=(V, E \cup {e} )$ and for the death of edge $e \in E$ we take $G^* = ( V, E \setminus {e} )$ and with reagrding .  

The BDMCMC seaerch algorithm converges to the joint posterior distribution \eqref{eq:posterior} given the birth and death rates as a ratio of the joint posterior distributions as follows
%\begin{equation}
%\label{birth rate}
%B_e( G, K ) = \min \left\{ \frac{P(G^{+e},K^{+e} |\mathbf{x})}{P(G,K |\mathbf{x})}, 1 \right\},  \text{ for each  } e \in \eb,
%\end{equation}
%\begin{equation}
%\label{death rate}
%D_e( G, K ) = \min \left\{ \frac{P(G^{-e},K^{-e} |\mathbf{x})}{P(G,K |\mathbf{x} )}, 1 \right\} ,  \text{ for each  } e \in E.
%\end{equation}
%We show the birth and death rates as follows
\begin{equation}
\label{eq:rate}
R_e( G, K ) = \min \left\{ \frac{P(G^{*},K^{*} |\mathbf{x})}{P(G,K |\mathbf{x})}, 1 \right\},  \text{ for each  } e \in \{ E \cup \eb \}.
\end{equation}
For the birth of edge $e \in \eb$ we take $G^*=(V, E \cup {e} )$ and for the death of edge $e \in E$ we take $G^* = ( V, E \setminus {e} )$ and with the regarding preposition matrix is $K^{*}$.  
Algorithm~\ref{alg:BDMCMC} represents the pseudo-code for the BDMCMC search algorithm. 
%Basically, in the CT-MCMC search algorithm, we only simulate the jump chain and store each tree which visits and the corresponding waiting time. 

%The ratio \eqref{eq:approx} can be used to speed up Bayesian structure learning in high-dimensional GGMs. 
%Basically, the ratio can potentially be used in most search algorithms that are constrained to move from one graph to another that differs by one edge, 
%see for example \cite{dobra2011bayesian, lenkoski2011computational, lenkoski2013direct, dobra2011computational, dobra2011copula, wang2012efficient, cheng2012hierarchical}. 

%Here we use the BDMCMC algorithm search algorithm developed by \cite{mohammadi2015bayesianStructure}.
%The BDMCMC algorithm explores the graph space by adding/removing edges via birth/death events which are modelled as independent Poisson processes.

%Given the results provided above, the BDMCMC search algorithm for the posterior sampling from a graph model is represented in Algorithm \ref{alg:BDMCMC}.

%  - - - - - - - - - - - - - - -  - - - - - - - - - - - - - - - - - - - - - - - - - - - - - -   |
\begin{algorithm} %[H]
\label{alg:BDMCMC}
\caption{ BDMCMC search algorithm }
 \KwIn{ A graph $G = (V,E)$ with a precision matrix $K$ and data $X$. }
 %initialization\;
 \For{$N$ iteration}{
  \For{all the possible moves in parallel}{
   Calculate the birth and death rates by Equation \ref{eq:rate}\;
   }
   Calculate the waiting time by Equation \ref{eq:waiting time}\;
   Update the graph based on the birth/death probabilities in Equation \ref{eq:prop rate}\;
   Update the precision matrix\;
 }
 \KwOut{ Samples from the joint posterior distribution \eqref{eq:posterior}. } 
\end{algorithm}
%  - - - - - - - - - - - - - - -  - - - - - - - - - - - - - - - - - - - - - - - - - - - - - -   |

The essential element of the BDMCMC search algorithm is that a continuous-time jump process is associated with the birth and death rates. 
Whenever a jump occurs, the corresponding move is always accepted, which can consider as more intelligent navigation of the graph space. 
The acceptance probabilities of commonly used RJMCMC algorithms are replaced by the waiting times in the BDMCMC algorithm. 
Correspondingly, graphs with high posterior probabilities have larger waiting times while graphs with low posterior probabilities have small waiting times and as a result, die quickly.
Another computational advantage of the BDMCMC algorithm is that the nested for loop, as a computationally expensive part of the algorithm, for computing the birth/death rates can be implemented in parallel since the rates associated with each edge can be calculated independently of each other.
We have implemented this part in parallel in the current version of the \textit{R} package  \textbf{BDgraph} \citep{BDgraph}.
These properties make the BDMCMC algorithm an efficient search approach to explore the graph space to identify the high posterior probability regimes, particularly for high-dimensional graphical models. 

%Another computationally advantage of the BDMCMC search algorithm is that the nested for loop, as a computationally expensive part of the algorithm, for computing the birth/death rates can be implemented in parallel as it is the case for our implementation in the BDgraph package. 

The main computational bottleneck of Algorithm \ref{alg:BDMCMC} is to evaluate the birth/death rates, which are based on the ratio of the posterior probabilities. 
These birth/death rates can be considered as the conditional Bayes factor of the comparison between graph $G$ and $G^{+e}$/$G^{-e}$, similar to \cite{hinne2014efficient}.
These ratios in Equation \ref{eq:rate} can be derived as
\begin{equation*}
%\label{eq:postrior ratio}
\frac{ P( G^{*}, K^{*} |\mathbf{x} ) }{ P( G,K |\mathbf{x} )} = \frac{I_G(\delta,\Omega) }{ I_{G^{*}}(\delta, \Omega) } H(K, \Omega+S, \delta+n,e),
\end{equation*}
where 
\begin{eqnarray}
\label{eq:Hij}
H(K,\Omega+S, \delta+n,e) = \left( \frac{|K^{*}|}{|K|} \right)^{\frac{\delta+n}{2} - 1}  \frac{ \exp \{\frac{-1}{2} \tr( K^{*} (\Omega+S) ) \} }{ \exp \{\frac{-1}{2} \tr( K ( \Omega+S ) ) \} }.
\end{eqnarray}
%and $\Omega^*=\Omega+S$. 
For details regarding how to compute the above function, see \cite{cheng2012hierarchical, mohammadi2015bayesianStructure, hinne2014efficient}.
We see that computing the ratio of posteriors requires evaluating the ratio of prior normalizing constants.
%\begin{equation}
%\label{eq:ratio}
%\frac{I_G(\delta,\Omega)}{I_{G^{-e}}(\delta, \Omega)}.
%\end{equation}
That is the main computational bottleneck of these types of search algorithms. 
%The computationally efficient way hitherto to evaluate this ratio is by exchange algorithm, see for example \cite{mohammadi2015bayesianStructure}. 

%  - - - - - - - - - - - - - - -  - - - - - - - - - - - - - - - - - - - - - - - - - - - - - -   |
\subsection{Existing methods to compute the normalizing constant}
\label{subsec:review-IG}

\textbf{Exact formula:} 
Recently, \citet{uhler2014exact} certify that it is possible to drive an explicit expression for the intractable normalizing constant for general graphs. 
Since the expression is (by its nature) mathematically complex, the capability of applying this intricate expression for Bayesian structure learning has yet to be investigated. 
%Recently, \citet[Theorem 3.3]{uhler2014exact} show that an explicit representation of the normalizing constant for general graphs is indeed possible and they drive an explicit expression for the intractable normalizing constant. 
%The capability of applying this expression for Bayesian structure learning in GGMs has yet to be investigated, since the expression is mathematically complex.
One possibility, as they point it out,  would be to find more computationally efficient procedures than \citet[Theorem 3.3]{uhler2014exact} for computing the normalizing constant for particular classes of graphs. 
%For example, one can develop specialized methods for computing the normalizing constants of different classes of graphs that are important for applications.
%Their work is a first step in the direction of developing new methodology for structure learning in GGMs based on the exact normalizing constant of the G-Wishart distribution. 

\textbf{Monte Carlo approximation:} 
\cite{atay2005monte} developed a Monte Carlo (MC) approach to approximate the normalizing constant based on the decomposition described in Section \ref{sec:ratio}. 
Although the MC approximation is accurate, it can be computationally expensive. 
In our simulation of Sections \ref{subsec:sim-disjoint} and \ref{subsec:sim-joint}, we faced numerical and computational issues of MC approximation for $p$ higher than 30. 

\textbf{Laplace approximation:} 
\cite{lenkoski2011computational} developed a Laplace approximation to compute $I_G(\delta, \Omega)$. 
Their approximation is based on using the iterative proportional scaling algorithm for computing the mode of the integral in Equation \ref{eq:I_G}.
This approximation is computationally faster than the MC approach, though it tends to be accurate only for the case of computing the posterior normalizing constant.%Thus, they recommend to apply the computationally expensive but more accurate MC integration for the prior normalizing constant and the computationally fast but less accurate Laplace approximation for the posterior normalizing constant.
Thus, they suggest using the Laplace approximation (as a computationally fast but less accurate approach) for the posterior normalizing constant and the MC integration (as a computationally expensive but more accurate approach) for the prior normalizing constant.

\textbf{Exchange algorithm:} 
\cite{murray2006mcmc} proposed the exchange algorithm for simulating from distributions, where prior distributions--like $G$-Wishart-- have intractable normalizing constants that varies according to the model.
These types of algorithms are also known as auxiliary variable approaches since they require exact sampling from the auxiliary variable to canceling out the ratio of normalizing constant in the Metropolis-Hastings acceptance probabilities \citep{park2018bayesian}. 
\cite{hinne2014efficient, lenkoski2013direct, mohammadi2015bayesianStructure} have implemented this algorithm in GGMs to avoid normalizing constant calculation by using the  exact sampler algorithm from $G$-Wishart distribution, proposed by \cite{lenkoski2013direct}. 
As state-of-the-art, this development has proven to yield significant computational improvement as it avoids the need for expensive approximations within the search algorithm. 
We briefly review the implementation of the exchange algorithm within the search algorithm;
For more details, see \cite[Section 5.2]{wang12}.

Suppose we want to compute the birth/death rate \eqref{eq:rate} for graph $G=(V,E)$ with the precision matrix $K$ as a current state of the search algorithm. 
By using the exchange algorithm, we can replace the intractable normalizing constant ratio with an estimate from a single sample at each parameter setting as
\begin{equation*}
%\label{ratio1}
\frac{I_G(\delta,\Omega)}{I_{G^{*}}(\delta, \Omega)}  \approx \frac{ |\widetilde{K}|^{\frac{\delta}{2} -1 } \exp \{\frac{-1}{2} \tr( \widetilde{K} \Omega ) \}  }{ |\widetilde{K}^{*}|^{\frac{ \delta }{ 2 } -1 } \exp \{\frac{-1}{2} \tr( \widetilde{K}^{*} \Omega ) \} }
\end{equation*}
where $\widetilde{K}$ has to be an exact sampler from the prior distribution, $W_G(\delta, \Omega)$. 
The exchange algorithm replaces the ratio of the intractable normalizing constants with an estimate from a single sample at each parameter setting. 
By using the above approximation, the birth/death rates will be
\begin{equation}
\label{eq:rate-DMH}
R_e( G, K ) \approx \min \left\{ \frac{ H(K, \Omega+S, \delta+n,e) }{ H( \widetilde{K},\Omega, \delta,e) }, 1 \right\},  \text{ for each  } e \in \{ E \cup \eb \},
%\frac{ P( G^{*}, K^{*} |\mathbf{x} ) }{ P( G,K |\mathbf{x} )} \approx \frac{ H(K, \Omega+S, \delta+n,e) }{ H( \widetilde{K},\Omega, \delta,e) },
\end{equation}
where function $H$ is given in Equation \ref{eq:Hij}. 
Essentially, the intractable prior normalizing constants have been replaced by an evaluation of function $H$ at $\widetilde{K}$ as an exact sample from the prior distribution $W_G(\delta, \Omega)$.

Algorithm \ref{alg:BDMCMC-DMHJ} represents the pseudo-code for the BDMCMC search algorithm combined with the exchange algorithm to compute the ratio of normalizing constant.
We call it a \textit{double} BDMCMC algorithm and consider it here as state-of-the-art.
For more details, see \cite{mohammadi2015bayesianStructure, hinne2014efficient}.
 
\begin{algorithm} %[H]
\label{alg:BDMCMC-DMHJ}
\caption{ \textit{Double} BDMCMC algorithm }
 \KwIn{ A graph $G = (V,E)$ with a precision matrix $K$ and data $X$. }
 %initialization\;
 \For{$N$ iteration}{
 Draw $\widetilde{K}  \thicksim W_G(\delta, \Omega) $\;
  \For{all the possible moves in parallel}{
   Calculate the birth and death rates by Equation \ref{eq:rate-DMH} \;
   }
   Calculate the waiting time by Equation \ref{eq:waiting time} \;
   Update the graph based on birth/death probabilities in Equation \ref{eq:prop rate} \;
   Update the precision matrix\;
 }
 \KwOut{ Samples from the posterior distribution \eqref{eq:posterior}. } 
\end{algorithm}

\begin{remark}
Algorithm \ref{alg:BDMCMC-DMHJ} requires exact sampling from the prior distribution of $G$-Wishart as a computationally expensive update within the BDMCMC search algorithm. 
Exact sampling from $G$-Wishart distribution, following \cite{lenkoski2013direct}, can be done by first sampling a standard Wishart variable from a full model and then using the iterative proportional scaling algorithm to place the variable in the correct space. 
It requires the solution of systems involving large matrices, in particular the inverse calculation of matrix $K$. 
\end{remark}
%In the following section we going to represent a simple an explicit formula to approximate the normalizing constant which does not require expansive exact sample form G-Wishart.  

%  - - - - - - - - - - - - - - -  - - - - - - - - - - - - - - - - - - - - - - - - - - - - - -   |
\subsection{Proposed method to compute the normalizing constant}
\label{subsec:formula-IG}

To bypass the computational bottleneck from the intractable normalizing constant in Algorithm \ref{alg:BDMCMC}, we represent a simple explicit analytic formula to approximate the normalizing constant as
\begin{equation}
\label{eq:approx}
\frac{I_{G^{-e}}(\delta,\I_p)}{I_{G}(\delta,\I_p)} \approx \frac{1}{2 \sqrt{\pi }}
%\frac{\Gamma(\delta/2)\Gamma((\delta+d)/2)\Gamma((\delta+1)/2)}{\Gamma((\delta+1)/2)\Gamma(\delta/2)\Gamma((\delta+d+1)/2)}  
\frac{ \Gamma( \frac{\delta+d}{2} ) }{\Gamma( \frac{\delta+d + 1}{2} )}
%\frac{\Gamma((\delta+d)/2)\Gamma((\delta+1)/2)}{\Gamma(\delta/2)\Gamma((\delta+d+1)/2)}=\frac{\Gamma((\delta+d)/2)}{\Gamma((\delta+d+1)/2)},
\end{equation}
where $d$ is the number of paths of length two linking the endpoints of $e$.
As is the case most of the time, in the absence of prior information, the parameter $\Omega$ is taken to be the $p$-dimensional identity matrix $\I_p$;
Throughout, we set $\Omega = \I_p$.
This approximation is exact in some cases, as we mentioned in Remark \ref{remark:equality}.
%The remainder of this paper is therefore  devoted to proving  this approximation and  analyzing its accuracy. 
The following sections are therefore devoted to proving  this approximation and  analyzing its accuracy.

%Here, we will perform Bayesian structure learning on high-dimensional simulated data using the BDMCMC approach developed by \cite{mohammadi2015bayesianStructure}. 
%Ratio \eqref{eq:ratio} is central to all the computations in this method and an approximation that bypasses the lengthy Monte Carlo computations is essential. 
%We bypass the computation of this ratio by replacing ratio \eqref{eq:ratio}  by approximation \eqref{eq:approx} or by using the double Metropolis-Hastings.
%Using the ratio \eqref{eq:ratio}  within the BDMCMC aglorithm is implemented in the \textbf{BDgraph} \textit{R} package \citep{BDgraph} in the function \textit{bdgraph()}. 

% - - - - - - - - - - - - - - - - - - - - - - - - - - - - - - - - - - - - - - - - - - - - - - - - - - - - - - - - - - - -|
\section{The ratio of normalizing constants}
\label{sec:ratio}

We first  recall a result by \cite{atay2005monte} which expresses $I_G(\delta,\Omega)$ as the product of a constant and an expectation. 
Let $K$ be the precision matrix and $K=\Psi^t \Psi$ its Cholesky decomposition where $\Psi$ is upper triangular  with positive diagonal elements. 
%% Let $\overline{E}$ be the complement of $E$ in the set of all possible edges in a graph with node set $V$, i.e. $\eb$ indexes the missing edges of $G$. 
%Let $\eb$ be the complement of $E$ that indexes the missing edges of $G$. 
%Similarly $\overline{E^{-e}}$ denotes the missing edges of $G^{-e}$. It was shown that 
Given the fact $K_{ij}=0$ for $(i,j)\in \overline{E}$, through simple matrix multiplication, we can verify
\[
\psi_E = \{ \psi_{ij} : (i,j)\in E \; \& \; \psi_{ii} : i \in V \}
\]
is in 1-1 correspondence with $K_E = \{ K_{ij} : (i,j) \in E  \; \& \; K_{ii} : i\in V \}$.
Also, the entries of $\psi_{\eb} = \{ \psi_{ij} : (i,j)\in \eb \}$ can be expressed in terms of $\psi_E$, a fact used in Proposition \ref{roposition:atay} below. 
Thus the entries of $\psi_E$ are called free variables while the entries of $\psi_{\eb}$ are non-free variables. 
Using the change of variables from $K_E$ to $\psi_E$, \cite{atay2005monte} prove the normalizing constant $I_G(\delta,\Omega)$ can be expressed as a known constant multiplied by the expected value of a function of $\psi_E$. 
In the particular case where $\Omega=\I_p$, which is of concern to us, the result is as follows.

\begin{proposition}
\label{roposition:atay}
For each node $i = \{ 1, \ldots, p \}$ of the undirected graph $G$, let $\nu_i$ be the number  of neighbours of $i$ which have a numbering larger than or equal to $i+1$. 
Then we have
\begin{equation*}
%\label{eq:IG_pro}
I_G(\delta,\I_p) = \left[ \prod_{i=1}^p \pi^{\frac{\nu_i}{2}} 2^{\frac{\delta}{2}+\nu_i} \Gamma \left( \frac{\delta+\nu_i}{2} \right) \right] \E \left( e^{ -\frac{D}{2} } \right)
\end{equation*}
where
\begin{equation*}
%\label{eq:D}
%f_E( \psi_E ) = \exp \left( \frac{-1}{2} \sum_{(i,j)\in \eb} \psi_{ij}^2 \right)
D = \sum_{(i,j)\in \eb} \psi_{ij}^2.
\end{equation*}
%and  the expected value is taken with respect to a product of independent $N(0,1)$ random variables $\psi_{ij}, (i,j)\in E$ and $\chi^2_{(\delta+\nu_i)/2}$ random variables $\psi_{ii}^2, i=1,\ldots,p$.  
The expected value is taken with respect to a product of independent random variables $\psi_{ij} \sim N(0,1)$ where $(i,j)\in E$ and random variables $\psi_{ii}^2 \sim \chi^2_{\delta+\nu_i}$ where $i = \{ 1, \ldots, p \}$.  
\end{proposition}

The value of $I_G(\delta,\I_p)$ is independent of the ordering of the nodes, so without loss of generality, in the remainder of this paper, we assume the nodes defining the edge $e$ are $q=p-1$ and $p$, that is the endpoints of $e$ are numbered last. 
For convenience, we write $\psi_e=\psi_{qp}$,
%\[ \psi_e=\psi_{qp} \]
 which is a non-free variable in the graph $G^{-e}$.  
 %, obtained from $G$ by removing the edge $e=(q,p)$, while it is  a free variable in  $G$. We have the following result.
\begin{corollary}
\label{prop:ratio}
Let $G^{-e}$ be the graph obtained from $G$ by removing the edge $e=(q,p)$. The ratio of the prior normalizing constants for $G^{-e}$ and $G$  is %then  equal to
\begin{equation}
\label{eq:ratio_pro}
\frac{I_{G^{-e}}(\delta,\I_p)}{I_{G}(\delta,\I_p)} = \frac{1}{2 \sqrt{\pi}}\frac{\Gamma(\frac{\delta}{2})}{\Gamma(\frac{\delta+1}{2})} 
\frac{\E \left( e^{ \frac{-1}{2} \left( D +\psi_{e}^2 \right) } \right) }{ \E \left( e^{ \frac{-1}{2} D } \right)  }.
%\frac{\E \left( \exp \left\{ \frac{-1}{2} \left( D +\psi_{e}^2 \right) \right\} \right) }{ \E \left( \exp \left\{ \frac{-1}{2} D \right\} \right)  }.
\end{equation}
%where $D$ is defined in \eqref{eq:D}.
\end{corollary}

Let $nb(i)$ denote the set of neighbours of $i\in V$ for $i = \{ 1, \ldots, p \}$. 
The proof of Corollary \ref{prop:ratio} is immediate if we observe that, since $\nu_i = |nb(i)\cap \{i+1,\ldots,p\}|$, the only $\nu_i$ that changes between $G^{-e}$ and $G$ is the node $\nu_q$ and, clearly, $\nu_q^{G^{-e}}=0$ while $\nu_q^{G}=1$. 
%Then \eqref{eq:ratio_pro} follows immediately from \eqref{eq:IG_pro}.

% - - - - - - - - - - - - - - - - - - - - - - - - - - - - - - - - - - - - - - - - - - - - - - - - - - - - - - - - - - - -|
\subsection{Reformulation of the ratio of normalizing constants}
\label{subsec:reformation}

We can drive the non-free entries of $\psi$ as
%From simple matrix multiplication (see also Proposition 2 in \citealp{atay2005monte}), we know that the non-free entries of $\psi$ are such that
\begin{equation}
\label{eq:aliye}
\psi_{1j}=0 \;\;\;  \text{and} \;\;\; \psi_{ij} = \frac{ -1 }{ \psi_{ii}} \sum_{l=1}^{i-1}\psi_{li}\psi_{lj}, \;i\not=1.
\end{equation}
The variables $\psi_{li}$ or $\psi_{lj}$ in the expression of $\psi_{ij}$ above may be free or non-free variables; see also \citet[Proposition 2]{atay2005monte}.

\begin{remark}
%\label{remark:up}
If $\psi_{ij}$ is non-free, it follows from Equation \ref{eq:aliye} that $\psi_{ij}$ can only be function of free variables $\psi_{lk}, l\not =k$ such that $l\leq i$ and $k<j$ and  $\psi_{ll}, l\leq i$.
\end{remark}

%We will say that there is a path of length 2 between nodes $q$ and $p$ if there exists $l\in \{1,\ldots, p-2\}$ such that both $(l,q)$ and $(l,p)$ belong to $E$. 
%A path of length 2 between nodes $q$ and $p$ if there exists a node in the graph that connect to both $q$ and $p$.
%Let $d$ be the number of paths of length 2 between nodes $q$ and $p$. 
Since the value of $I_G(\delta,\I_p)$ does not depend upon the order of the nodes, from now on in this paper, we assume the nodes which are neighbours to both $q$ and $p$, are numbered $p-1-d,p-1-(d-1),\ldots, p-1-1$ where $d$ is the number of paths of length 2 between nodes $q$ and $p$;
See for example the node orders in Figure \ref{fig:paths}. 
With this convention, we have $\psi_e=A+b$ where
\begin{align} 
\label{eq:A}
A &=  \frac{-1}{\psi_{qq}} A_1 \;\; \text{ where } \;\; A_1 = \sum^{q-1}_{l = q - d}\psi_{lq}\psi_{lp},  \\ 
\label{eq:b}
b &=  \frac{-1}{\psi_{qq}} b_1  \;\;\; \text{ where } \;\;\ b_1 = \sum^{q-(d+1)}_{l=1} \psi_{lq}\psi_{lp}.
\end{align}

% For the sake of brevity, we will also use the notations
%\begin{equation}
%\label{D}
%b_1=\psi_{qq}b \;\;\; \text{and}  \;\;\;  D=\sum_{(i,j)\in \eb}\psi_{ij}^2.
%\end{equation}

\begin{remark}
%\label{remark:independence}
The numbering we have adopted for nodes that are neighbours both to $q$ and $p$ ensures that $A$ is independent of $b$ and $D$.
\end{remark}

With the notations above, Equation \ref{eq:ratio_pro} can be written
\[
\frac{I_{G^{-e}}(\delta,\I_p)}{I_{G}(\delta,\I_p)} = \frac{1}{2 \sqrt{\pi }} \frac{ \Gamma( \frac{\delta}{2}) }{\Gamma( \frac{\delta + 1}{2} )} \frac{\E \left( e^{-\frac{D}{2}} e^{-\frac{(A+b)^2}{2}} \right) }{ \E \left( e^{-\frac{D}{2}} \right) }.
\]
Our aim is to approximate this ratio and, towards this goal, we have the following approximation
\begin{equation}
\label{eq:aim}
\E \left( e^{-\frac{D}{2}} e^{-\frac{(A+b)^2}{2}} \right) \approx \E \left( e^{-\frac{D}{2}} \right) \E \left( e^{-\frac{A^2}{2}} \right).
\end{equation}
If we prove that the above approximation holds, then we will have
\begin{equation*}
%\label{eq:nearly}
\frac{I_{G^{-e}}(\delta,\I_p)}{I_{G}(\delta,\I_p)} \approx \frac{1}{2 \sqrt{\pi }} \frac{ \Gamma( \frac{\delta}{2}) }{\Gamma( \frac{\delta + 1}{2} )} \E \left( e^{-\frac{A^2}{2}}\right).
\end{equation*}
%By \eqref{eq:EA} in Section \ref{sec:proposition EA} of the Supplementary File, we have the analytic expression
Regarding Proposition \ref{proposition:EA} of the Supplementary File, we have the analytic expression
\[
\E \left(e^{-\frac{A^2}{2}}\right)= \frac{ \Gamma( \frac{\delta + 1}{2}) }{\Gamma( \frac{\delta}{2} )}  \frac{ \Gamma( \frac{\delta+d}{2} ) }{\Gamma( \frac{\delta+d + 1}{2} )} 
\]
%and thus \eqref{eq:nearly} becomes
and thus we have
\begin{equation*}
%\label{approx}
\frac{I_{G^{-e}}(\delta,\I_p)}{I_{G}(\delta,\I_p)} \approx \frac{1}{2 \sqrt{\pi }}
%\frac{\Gamma(\delta/2)\Gamma((\delta+d)/2)\Gamma((\delta+1)/2)}{\Gamma((\delta+1)/2)\Gamma(\delta/2)\Gamma((\delta+d+1)/2)}  
\frac{ \Gamma( \frac{\delta+d}{2} ) }{\Gamma( \frac{\delta+d + 1}{2} )}
%\frac{\Gamma((\delta+d)/2)\Gamma((\delta+1)/2)}{\Gamma(\delta/2)\Gamma((\delta+d+1)/2)}=\frac{\Gamma((\delta+d)/2)}{\Gamma((\delta+d+1)/2)},
\end{equation*}
which is the approximation \eqref{eq:approx} that we want to prove. 
%Therefore the remainder of this paper is devoted to proving the approximation in Equation \ref{eq:aim} and analyzing its accuracy. 

\begin{remark}
\label{remark:equality}
It is important to note that in Equation \ref{eq:aim} if $b=0$, then our approximation \eqref{eq:approx} is exact. 
This means that when there are only paths of length 2, or no path, between nodes $q$ and $p$, the approximation is exact. 
It is interesting to note that this happens also in other cases. 
In fact, \citet[Theorem 2.5]{uhler2014exact} show that if $G^{-e}$ is such that $G$ is decomposable, then our approximation \eqref{eq:approx} is exact.
\end{remark}

%% We will show that always $I_1/I_2 \leq 1$ and thus the relative error $(I_2-I_1)/I_2=1-I_1/I_2$ is always positive.
%% In the case where the paths between $q$ and $p$ are disjoint, we are able to give an analytic upper bound for the relative error. 
%This is done in Theorem \ref{theorem:disjoint}. In Section \ref{sec:joint}, we consider the case where the paths are not necessarily disjoint.  
%We then express $I_1/I_2$ as a single integral in \eqref{relationship3}, and give an approximation \eqref{eq:newi1i2} to this integral which, under some conditions, is close to 1.
%%  As mentioned before, this approximation could be far from  1 and yet $I_1/I_2$ could be close to 1: see Figure \ref{fig:boxplot disjoint}.

\begin{lemma}
\label{lemma:I_1}
Using the quantities, $D$, $A$, $b$, and $b_1$ defined above, we have
\begin{equation*}
%\label{eq:E1}
\E \left( e^{-\frac{ D}{2} - \frac{(A+b)^2}{2}} \right) 
    = \E \left( e^{-\frac{ A^2}{2}} \right) \E \left( e^{-\frac{ D}{2} }  \E \left( h \left( b_1, \delta^* \right) \big| \Psi_{\cup}^- \right) \right),
\end{equation*}
where $\delta^* = \frac{\delta+d}{2}$ and
\begin{equation}
\label{eq:int-b1}
 h(b_1, \delta^*) =  \frac{ 2^{-\delta^*} }{ \Gamma( \delta^* )}  \int_0^{+\infty} y^{\delta^* - 1} e^{ \frac{-1}{2} \left( y + \frac{b_1^2}{y} \right) }  dy,
\end{equation}
and 
\begin{equation}
\label{eq:smallD}
% \Psi_{E_{qp}^+} = 
%\psi_{\tilde{D}_1} = 
\Psi_{\cup}^- = \{ \psi_{ij} : (i,j) \in E \setminus \left( E_q \cup E_p \right) \}.
\end{equation}
where $E_q = \{ (i,j): (i,q) \in E \}$ and $E_p = \{ (i,j): (i,p) \in E \}$.
$\Psi_{\cup}^-$ includes all the free elements of the matrix $\Psi$ except those are the neighbors of nodes $p$ and $q$. % in graph $G$.
\end{lemma}
%Note that $D$ and $b_1$ are $\psi_{\tilde{D}}$-measurable, that is, are functions of the elements of $\psi_{\tilde{D}}$ only. 
The proof is given in Section \ref{sec:lemma:I_1} of the Supplementary file.
Regarding to the above lemma, proving
\begin{align*}
\E \left( e^{-\frac{ D}{2} }  \E \left( h \left( b_1, \delta^* \right) \big| \Psi_{\cup}^- \right) \right) \approx \E \left( e^{-\frac{D}{2}} \right) 
\end{align*}
leads to the approximation in Equation \ref{eq:aim}. 
For convenience, we will also adopt the notation
\begin{equation*}
%\label{eq:i1i2}
I_1 = \E \left( e^{-\frac{ D}{2} }  \E \left( h \left( b_1, \delta^* \right) \big| \Psi_{\cup}^- \right) \right) \;\ \text{and} \;\; I_2 = \E \left( e^{-\frac{D}{2}} \right) .
\end{equation*}
and therefore
\begin{equation}
%\label{eq:relationship2}
\label{eq:i1i2}
\frac{I_1}{I_2} = \frac{ \E \left( e^{-\frac{ D}{2} }  \E \left( h \left( b_1, \delta^* \right) \big| \Psi_{\cup}^- \right) \right) }{ \E \left( e^{-\frac{ D}{2}} \right) }.
\end{equation}
Note that the accuracy of our approximation in Equation \ref{eq:approx} is represented by how close is the above ratio $I_1/I_2$ to 1. 
Thus, proving that our approximation is accurate is equivalent to prove that $I_1/I_2$ can accurately be approximated by 1.
For example, for the cases that $I_1/I_2$ is equal to 1, our approximation is exact.
%Therefore, the remainder of this paper is devoted to proving the ratio $I_1/I_2$ is close to 1 and analyzing its accuracy. 

%\begin{remark}
%\label{if}
%If we could show that, whatever the value of $\psi_{\tilde{D}}$, the expectation $\E(e^{-\d b_1^2/Y}\mid \psi_{\tilde{D}})$ can uniformly be approximated by $1$, then from \eqref{eq:relationship2}, it would follow that $I_1/I_2$ can also be approximated by $1$. 
%We are not able to quite achieve this goal but in the next few subsections, we are going to show that, conditional on $\psi_{\tilde{D}_1}\subset \psi_{\tilde{D}}$ defined below in \eqref{eq:smallD}, the distribution of $b_1$ is a scale mixture of normal distributions. 
%We then use this scale mixture of distributions admits a unique $N(0, v_D)$ approximation in the  $ L^2(\R)$ sense. And finally, we show that a sufficient condition for $\E \left( e^{ \frac{-b_1^2 }{2 Y} } \big| \psi_{\tilde{D}} \right)$ to be close to $1$ is that $v_D$ is close to $0$. 
%We verify this result numerically. Besides, our numerical results show our approximation $I_1/I_2$ (given by $I_{1/2,\mbox{Gamm}}$) is accurate (close to 1) even for the cases that $v_D$'s are not close to $0$. 
%In fact, both set of values for $I_{1/2,\mbox{MC}}$ and $I_{1/2,\mbox{Gamm}}$ seem to be affected by the size of $v_D$ but are reasonably close to 1, whatever the value of $v_D$.
%\end{remark}

\begin{remark}
\label{remark:I1I2}
It is important to mention that $I_1/I_2$ is always equal to or less than 1 ($I_1/I_2 \leq 1$).
It follows immediately from Equation \ref{eq:int-b1} since $b_1^2/Y$ is always positive and $e^{-b_1^2 / Y} \leq 1$.
\end{remark}

\begin{remark}
%\label{remark:if}
If we could show, whatever the value of $\Psi_{\cup}^-$,  the expectation $\E \left( h \left( b_1, \delta^* \right) \big| \Psi_{\cup}^- \right)$ can uniformly be approximated by $1$, it would follow that $I_1/I_2$ can also be approximated by $1$. 
We are not able to quite achieve this goal but, first, in the next Section (Theorem \ref{theorem:disjoint}), we establish the approximation with explicit bounds in the special case when all paths between $q$ and $p$ are disjoint. 
%The key to proving this result is the fact that $b_1$ can be expressed as a linear product of independent normal variables in this case. 
The key to proving this result is the fact that $b_1$ can be expressed as a linear product of independent normal variables, for the cases of disjoint paths. 
Then, in Section \ref{sec:joint}, we show, conditional on $\Psi_{\cup}^-$ defined in Equation \ref{eq:smallD}, the distribution of $b_1$ is a scale mixture of normal distributions. 
We then use this scale mixture of distributions to admit a unique $N(0, v_D)$ approximation. % in the  $ L^2(\R)$ sense. 
Finally, we show that a sufficient condition for $\E \left( h(b_1, \delta^*) \big| \Psi_{\cup}^- \right)$ to be close to $1$ is that $v_D$ is close to $0$.
%We verify this result numerically. 
%Besides, our numerical results show our approximation $I_1/I_2$ (given by $I_{1/2,\mbox{Gamm}}$) is accurate (close to 1) even for the cases that $v_D$'s are not close to $0$.
%%If $v_D$ is not close to $0$, we see that numerically still, the approximation $I_1/I_2$, given by the values of $I_{1/2,\mbox{Gamm}}$  is close to 1. 
%In fact, both set of values for $I_{1/2,\mbox{MC}}$ and $I_{1/2,\mbox{Gamm}}$ seem to be affected by the size of $v_D$ but are reasonably close to 1, whatever the value of $v_D$.
\end{remark}

% - - - - - - - - - - - - - - - - - - - - - - - - - - - - - - - - - - - - - - - - - - - - - - - - - - - - - - - - - - - -|

\section{The ratio for the case disjoint paths}
\label{sec:disjoint}

\begin{figure} %[t]
\centering
\begin{center}
\includegraphics[width=4in]{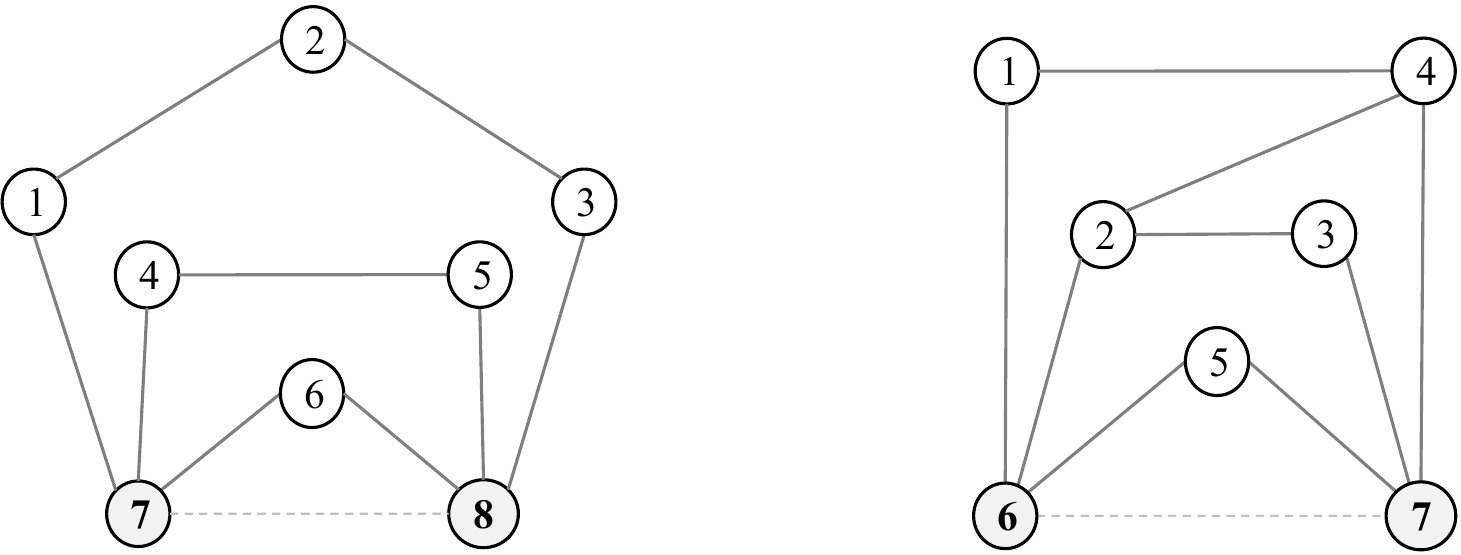} 
\end{center}
  \caption{ (Left) A graph  with disjoint paths between $q$ and $p$. (Right) A graph with several non-disjoint paths between $q$ and $p$. \label{fig:paths} }  
\end{figure}

%An ordered subset of nodes in an undirected graph is a \textit{path} if each set of consecutive nodes in the ordering are joined by an edge and the \textit{path length} is the number of edges. 
A \textit{path} is a sequence of nodes in which each node is connected by an edge to the next and the \textit{path length} is the number of edges between them.
Two paths between $q$ and $p$ are \textit{disjoint} if they have no node other than $p$ and $q$  in common. 
For example, in the left-hand side graph of Figure \ref{fig:paths}, the paths between $q=7$ and $p=8$ are
%Figure \ref{fig:paths} (left) shows that the only paths between $q=7$ and $p=8$ are
\[
 \lambda_1 = \{q, 1, 2, 3, p \}, \;\; \lambda_2 = \{ q, 4, 5, p \}, \;\; \lambda_3 = \{ q, 6, p \},
\]
and they are \textit{disjoint} paths.  %, while the graph in Figure \ref{fig:paths} (right) has two joint paths paths between $q=6$ and $p=7$ as follow
%\[ \{q, 1, 4, p \} \;\; \text{and} \;\; \{ q, 2, 3, p \}. \]
%and they are clearly joint.

%We now set up notations specific to this case. 
%A path $\lambda=\{q,1_{\lambda},\ldots, i_{\lambda}, \ldots, j_{\lambda},\ldots, p\}$ is said to be minimal if it has no chord, i.e., there is no edge $(i_{\lambda}, j_{\lambda})$ in $E$. If the paths between $q$ and $p$ are disjoint, they are necessarily minimal. 
 
%Next, we note that the prime component $P_{qp}$ containing $q$ and $p$ is necessarily the graph induced from $G$ by $\Lambda$. 
%Indeed, if $v$ is a node in $P_{qp}$ but does not belong to $\lambda$, then it is linked to at most one node $i_1$ in $\Lambda$ because if it was linked to more than one node, say it was linked to both $i_1$ and $i_2$, 
%then there would be an additional path using the sub-path $\{i_1,v,i_2\}$: this is impossible since $v$ does not belong to $\Lambda$.

A path $\lambda \in \Lambda$ of length $\ell_{\lambda}+1$ is a sequence of distinct nodes as
$\lambda = \{q,1_{\lambda}, 2_{\lambda},\ldots, \ell_{\lambda},p\}$
where $(q,1_{\lambda}), \ldots, ( i_{\lambda}, (i+1)_{\lambda} ), \ldots,  (\ell_{\lambda},p)$ are edges of $G$;
The set of all such paths $\lambda$ between $q$ and $p$ is denoted $\Lambda$.  
 We let $E_{\lambda}$ and $V_{\lambda}$ be, respectively, the set of edges, the set of interior nodes of $\lambda$ and the set of interior points deprived of $1_{\lambda}$, i.e.
\[
E_{\lambda}=\left\{ (q, 1_{\lambda}), (1_{\lambda}, 2_{\lambda}),\ldots, (\ell_{\lambda},p) \right\}, \;\; V_{\lambda}=\{1_{\lambda},2_{\lambda},\ldots, \ell_{\lambda}\}.
%,\;\;V^{(-1)}_{\lambda} = V_{\lambda}\setminus \{1_{\lambda}\}.
\]
If $L=|\Lambda|$ is the total number of  paths, we set an arbitrary order $\lambda_1, \ldots, \lambda_L$ of the paths where, for convenience, we list the paths of length 2, i.e. $\ell_{\lambda}=1$ last. The nodes $q$ and $p$ are ranked last so that the order of the nodes in $V$ is
\begin{equation*} %\label{order}
1_{\lambda_1}, \ldots,  \ell_{\lambda_1}, 1_{\lambda_2}, \ldots, \ell_{\lambda_2}, \ldots,\ldots, 1_{\lambda_L}, \ldots,  \ell_{\lambda_L}, q,p.
\end{equation*}
Using these notations, the following lemma gives the expression for $\psi_e$ in terms of the free variables $\psi_E$.

\begin{lemma}
\label{lemma:psie}
In the model with underlying graph $G^{-e}$, the variables $\psi_{qp}=\psi_e$ of the Cholesky decomposition of the precision matrix $K$ is expressed in terms of $\psi_E$ as
\begin{equation}
\label{eq:psie}
\psi_e = \frac{1}{\psi_{qq}} \sum_{\lambda\in \Lambda}(-1)^{\ell_{\lambda}}\frac{\prod_{a\in E_{\lambda}}\psi_{a}}{\prod_{v \in V_{\lambda}\setminus \{1_{\lambda}\}} \psi_{vv}}.
\end{equation}
\end{lemma}

The proof relies on a repeated application of Equation \ref{eq:aliye}. 
The proof is given in Section \ref{sec:lemma:psie} of the Supplementary file. 
We illustrate these calculations with the following example.

\begin{example}
Consider the graph of Figure \ref{fig:paths} (left). 
The upper triangular matrix $\Psi$ is 
\begin{equation*}
\Psi = \begin{bmatrix}
\psi_{11} & \psi_{12} & 0            & 0             & 0             & 0             & \psi_{17} & 0             \\
               & \psi_{22} & \psi_{23}& 0             & 0             & 0             & *             &  0             \\
               &                & \psi_{33}& 0             & 0             & 0             & *             & \psi_{38}  \\
               &                &               & \psi_{44} & \psi_{45} & 0             & \psi_{47} & 0              \\
               &                &               &                & \psi_{55} & 0             & *             & \psi_{58}  \\
               &                &               &                &                & \psi_{66} & \psi_{67} & \psi_{68}  \\
               &                &               &                &                &                & \psi_{77} & *               \\
               &                &               &                &                &                &                & \psi_{88} 
\end{bmatrix}
\end{equation*}
where the entries marked with a ``$*$" are the non-free entries and are given as
\begin{eqnarray*}
\psi_{27} = -\frac{\psi_{12}\psi_{17}}{\psi_{22}}, \;\; \psi_{37} = \frac{\psi_{17}\psi_{12}\psi_{23}}{\psi_{22}\psi_{33}}, \;\; 
\psi_{57} = - \frac{\psi_{45}\psi_{47}}{\psi_{55}},
\end{eqnarray*}
and 
\begin{align*}
\psi_{78} &=-\frac{1}{\psi_{77}} \left(   \psi_{67}\psi_{68} + \psi_{57}\psi_{58} + \psi_{37}\psi_{38} \right)\\
               &= \frac{1}{\psi_{77}} \left( - \psi_{67}\psi_{68} + \frac{ \psi_{47} \psi_{45} \psi_{58} }{ \psi_{55} }
                    - \frac{ \psi_{17} \psi_{12} \psi_{23} \psi_{38} }{ \psi_{22} \psi_{33} } \right) .
\end{align*}
Equation \ref{eq:psie} is verified. 
We see that the different terms in $\psi_{qp}=\psi_{78}$ above concern, successively, the paths of length 2, 3, and 4.
\end{example}

We are now in a position to state and prove the first of our two main results regarding the error made of our approximation in Equation \ref{eq:approx} or equivalently the approximation in Equation \ref{eq:aim}.

\begin{theorem}
\label{theorem:disjoint}
%For the case where  in the graph $G$ the paths between the end points of the edge $e=(q,p)$ are disjoint, and for $I_1$ and $I_2$ as defined in \eqref{eq:i1i2}, the ratio $I_1/I_2$ 
For the case where in the graph $G$ the paths between the endpoints of the edge $e=(q,p)$ are disjoint, the ratio $I_1/I_2$ \eqref{eq:i1i2} 
%which is equal to the ratio of $\frac{I_{G^{-e}}(\delta,\Omega)}{I_{G}(\delta,\Omega)} $ and  its approximation \eqref{eq:approx},  
is such that
% we have
\begin{equation}
\label{ineq:approx}
B(\delta,d, \ell_{\lambda}) \leq \frac{I_1}{I_2}\leq 1,
\end{equation}
where
\begin{equation}
\label{eq:B}
B(\delta,d, \ell_{\lambda}) = 1 - \frac{\delta^2}{\pi (\delta + 2)} \left( \frac{ \Gamma( \frac{\delta}{2} ) }{ \Gamma( \frac{\delta+1}{2} )} \right)^2 r(\delta+d-1) \sum_{\lambda\in \Lambda}{r(\delta)}^{\ell_{\lambda}},
\end{equation}
with $\Lambda$ being the set of paths between $q$ and $p$, $d$ the number of paths of length 2, and
\[
r(\delta)=\frac{\Gamma(\frac{\delta}{2})}{\sqrt{\pi}\Gamma(\frac{\delta+1}{2})}.
%% \;\; \text{and} \;\; r(\delta)= \frac{2}{\delta} \left[ \frac{\Gamma( \frac{\delta+1}{2} )}{\Gamma( \frac{\delta}{2} )} \right]^2.
\]
%\newline Moreover, $A$ is independent of $\sum_{(i,j)\in \eb}\psi_{ij}^2$ and we have
%\begin{equation}
%\E(e^{-\frac{A^2}{2}})=\frac{\Gamma((\delta+d)/2)\Gamma((\delta+1)/2)}{\Gamma(\delta/2)\Gamma((\delta+d+1)/2)}.
%\end{equation}
With an accuracy given by Equation \ref{ineq:approx}, we have the approximation
\begin{equation*} %\label{eq:approx}
\frac{I_{G^{-e}}(\delta,\I_p)}{I_G(\delta,\I_p)}\approx \frac{1}{2\sqrt{\pi}}\frac{\Gamma(\frac{\delta+d}{2})}{\Gamma(\frac{\delta+d+1}{2})}.
\end{equation*}
\end{theorem}

\textbf{Proof.} 
The proof is given in Section \ref{sup:proof-theorem1} of the Supplementary file.
The proof is based on the fact that the expression of $b_1$ \eqref{eq:b} can be expressed as a linear product of independent normal variables in the case the paths between $q$ and $p$ are disjoint.

%The reader not interested in the details of the proof of Theorem \ref{theorem:disjoint} can now move to Section \ref{sec:joint} but should remember the definitions of the quantities $A, A_1, b, b_1, b_{\lambda}$ and $Q_{\delta}$ defined  below.

%Before leaving the case of disjoint paths, let us note that  this case  holds  also if the following two conditions are satisfied:
%\begin{enumerate}
%\item The graph $G_{qp}$ induced by a set of disjoint paths between $q$ and $p$ can be separated by a complete separator from the remaining nodes.
%\item The nodes that do not belong to this induced graph are labelled with numbers strictly less than those in $G_{qp}$.
%\end{enumerate}
%An example of such graph is given in Figure \ref{paths} (c). Indeed, it is immediate to verify that the computation of $\psi_e$ using \eqref{eq:aliye} cannot use any path involving nodes with a numbering less than those in $G_{qp}$. 
%The fact that $\psi_e$ only depends on the paths in $G_{qp}$ also follows from the fact that the nodes with a numbering less than those in $G_{qp}$ belong to a prime component different from that containing $G_{qp}$. 

% - - - - - - - - - - - - - - - - - - - - - - - - - - - - - - - - - - - - - - - - - - - - - - - - - - - - - - - - - - - -|
\subsection{Simulated experiments for the case of disjoint paths} 
\label{subsec:sim-disjoint}

\begin{figure} %[!ht]
\begin{center}
	  \includegraphics[width=6in]{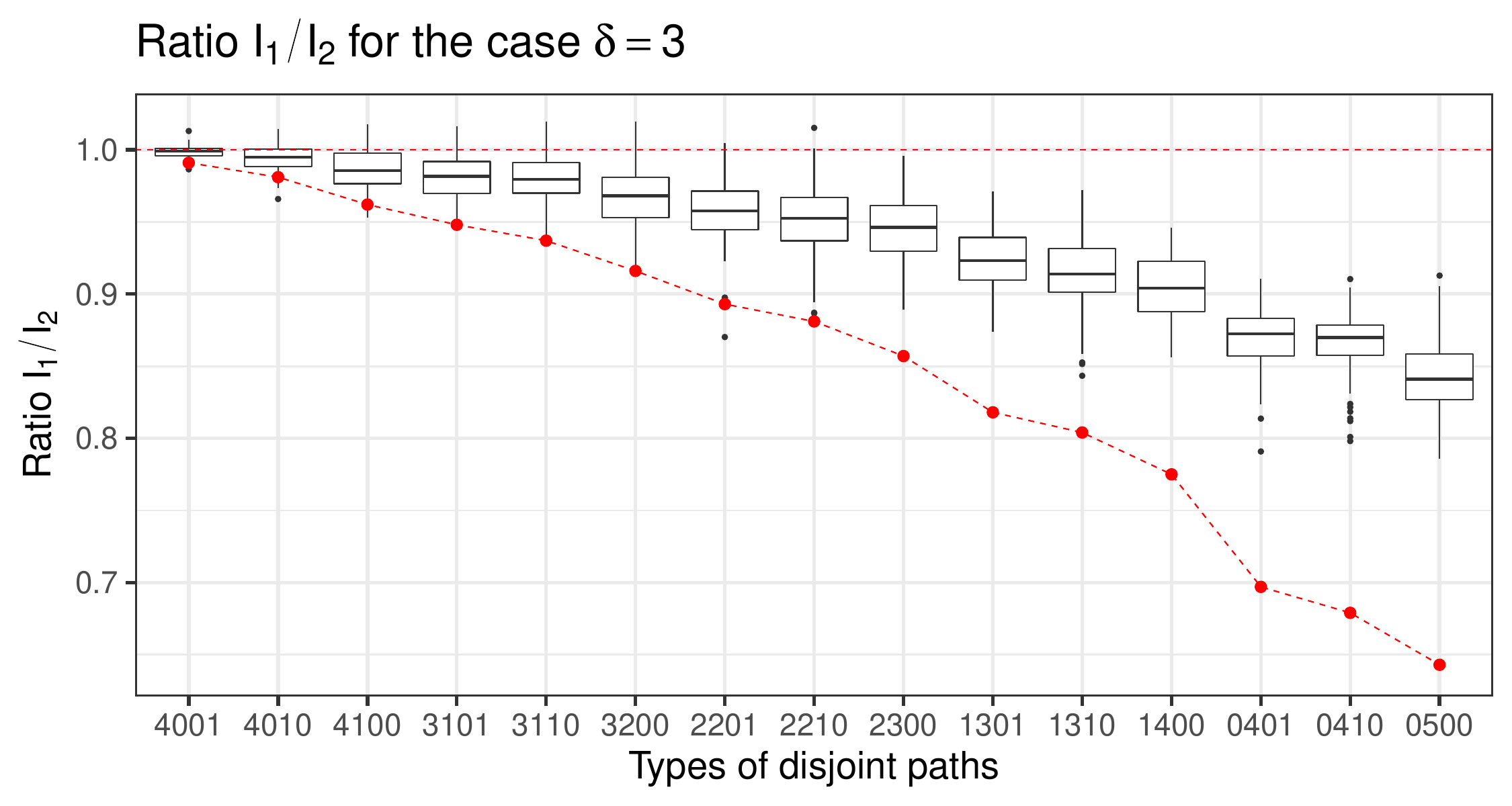}
	  \includegraphics[width=6in]{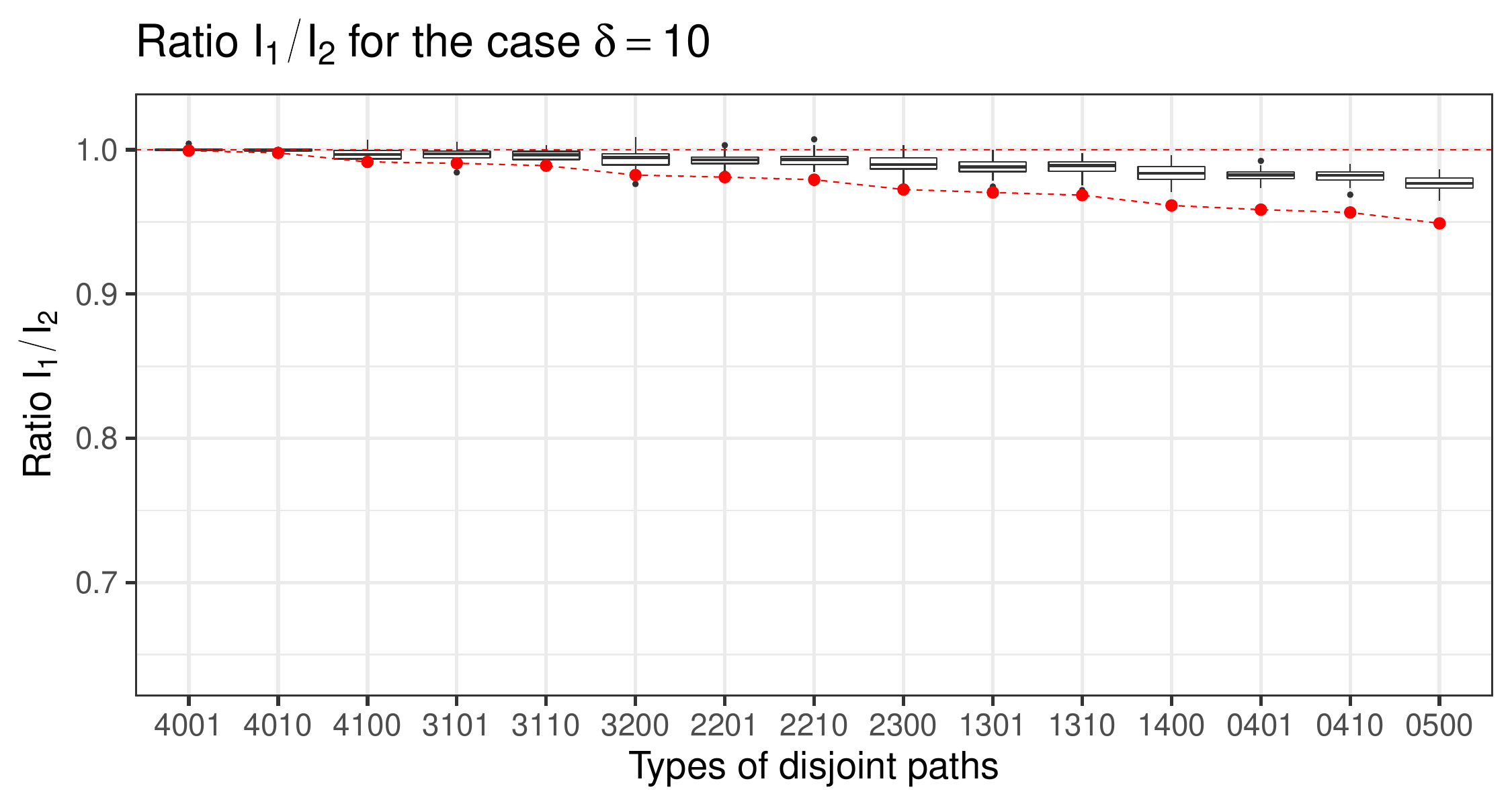}
\end{center}
\caption{ The ratio $I_1/I_2$ and its bound $B$ in Equation \ref{eq:B} for $\delta=3$ (top) and $\delta=10$ (bottom). The red dotted line is the lower bound $B$ and the boxplots are the $I_1/I_2$ computed by the MC algorithm of \cite{atay2005monte}, with over 100 replications. The $15$ different graphs are indicated on the horizontal axis. Each sequence of four digits indicates the number of paths of lengths 2, 3, 4, and 5 in the graph. For example, ``$3110$" represents a graph with $3$ disjoints paths of length $2$, $1$ of length $3$, $1$ of length $4$, and $0$ of length $5$. \label{fig:boxplot disjoint} }
\end{figure}

To illustrate the results in Theorem \ref{theorem:disjoint}, we report the ratio $I_1/I_2$ \eqref{eq:i1i2} following the MC approach of \cite{atay2005monte} as well as the lower bound $B$ in Equation \ref{eq:B}.
We note that, if $I_1/I_2\approx 1$ our approximation is good, without any additional conditions. %, MC method should reflect that by being close to 1.  
Note that, $1-I_1/I_2$ reflects the error rate of our approximation \eqref{eq:approx} for the prior normalizing constant of $G$-Wishart. 
Since $I_1/I_2$ and $B$ are functions of $\delta$ and type of disjoint paths ($d$ and $\ell_{\lambda}$), 
our simulation is based on graphs with different types of disjoint paths as well as different values of $\delta$.
%Since $I_1/I_2$ and $B(\delta,d, \ell_{\lambda})$ are function of $\delta$ and type of disjoint paths ($d$ and $\ell_{\lambda}$) in graphs, here report these values for graphs with different type of disjoint paths and then for different values of $\delta$.
%We consider $15$ different scenarios for graphs having five different paths between $q$ and $p$. 
We consider $15$ different types of graphs with five different paths between $q$ and $p$. 
These graphs are indicated on the horizontal axis in Figure \ref{fig:boxplot disjoint}. 
Each sequence of four digits denotes the number of paths of length $2$, $3$, $4$, and $5$ in the graphs. 
For example, ``$3110$" indicates a graph configuration with $3$ disjoints paths of length $2$, $1$ of length $3$, $1$ of length $4$, and $0$ of length $5$.

Figure \ref{fig:boxplot disjoint} represents the values of  $I_1/I_2$ (over $100$ replications) as well as the lower bound $B$ \eqref{eq:B} for two values of $\delta = 3$ and  $\delta = 10$. 
%The boxplot of the numerical values of $I_1/I_2$ obtained over $100$ replications. % here for the particular graphs with disjoint paths.
%We see that  $I_1/I_2$ are close to 1, always greater than 0.88. The relative error is always less than $0.12$.
The worst-case scenarios are for the case $\delta=3$ and no paths of length two ($d=0$), likes the graph ``$0500$" which has $5$ paths of length $3$ and no other type of paths; 
These types of graphs are highly unlikely cases. 
Even for this case, the relative error is around $0.12$.
For the case $\delta = 10$, we see that our approximation has pretty good performance with the maximum relative error $1-I_1/I_2$ around $0.025$.

\begin{figure} [ht]
\begin{center}
	  \includegraphics[width=6.5in]{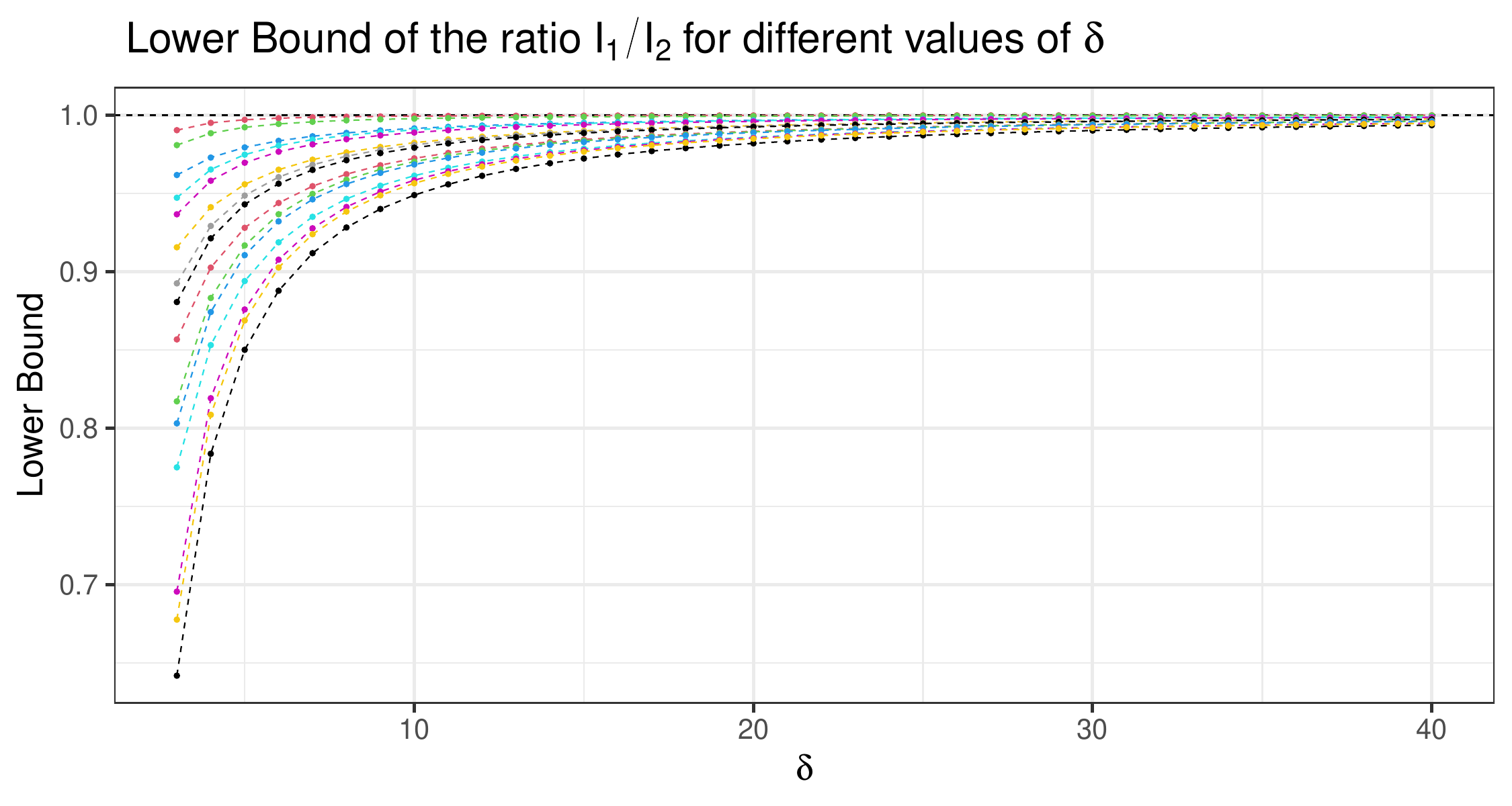}
\end{center}
\caption{ Plot visualization for the lower bound $B$ \eqref{eq:B} for $\delta = \{ 3, 4, ..., 40 \}$ and for the $15$ different graphs which are indicated on the horizontal axis in Figure \ref{fig:boxplot disjoint}. Each dotted line represents the $B$ values for a graph with specific types of paths. For example, the black bottom line is for the case ``$0500$" which means a graph with $0$ disjoints paths of length $2$, $5$ of length $3$, $0$ of length $4$, and $0$ of length $5$. \label{fig:Boundary} }
\end{figure}

Figure \ref{fig:Boundary} reports the values of the lower bound $B$ for different values of $\delta$ ($\delta = \{ 3, 4, ..., 40 \}$) and for the $15$ different graphs which are indicated on the horizontal axis in Figure \ref{fig:boxplot disjoint}.
Each dotted line represents the $B$ values for a specific graph with different type of paths. 
For instance, the black bottom line is for the configuration ``$0500$''. % which indicates a graph with $5$ disjoints paths of length three and no paths of length two, four, and five.
In general, this plot indicates that the accuracy of our approximation is increased by increasing the value of $\delta$.
As we can see the worst-case scenario is for the minimum value of $\delta$($=3$), while for the cases $\delta > 10$ the lower bound $B$ for our approximation is cloth to 1. % has a pretty good performance. 

% - - - - - - - - - - - - - - - - - - - - - - - - - - - - - - - - - - - - - - - - - - - - - - - - - - - - - - - - - - - -|
\section{The ratio in general case}
\label{sec:joint}

When the paths between $q$ and $p$ are not disjoint, the expression of $b_1$ \eqref{eq:b} becomes more complicated. 
It can be expressed in terms of variables $\psi_{jp}\sim N(0,1) , j<p$ and variables of the type
\begin{equation*}
%\label{beta}
X_{ij}=\frac{\psi_{ij}}{\psi_{jj}},\;\;i<j,
\end{equation*}
where $\psi_{ij}\sim N(0,1)$ and $\psi_{jj}^2\sim \chi^2_{\delta+\nu_j}$.
%These are variables which, as in the previous section, are distributed like $ Z/\sqrt{Q}$ with $Z\sim N(0,1)$ independent of $Q^2\sim \chi^2_{\delta+\nu_j}$. 
As a toy example, for the graph of Figure \ref{fig:paths} (right) with tedious computations yield
\[
b_1 =  \psi_{26} X_{23}^2 X_{24} \psi_{47} 
                        + \psi_{26} X_{24} \psi_{47} 
                        + \psi_{16} X_{14} \psi_{47} 
                        + \psi_{26} X_{23}\psi_{37}.
\]
For the details, see Example \ref{ex:smallgraph} in Section \ref{subsec:bilinear} of the Supplementary File. 
We see that $b_1$ is the sum of polynomials in $X_{ij}, (i,j)\in E$ multiplied by the product of two independent $N(0,1)$. %, as in the expression  \eqref{b1lambda}  of $b_{1\lambda}$. 
But, unlike in the case of disjoint paths between $q$ and $p$, the polynomials here are not linear in each $X_{ij}$;
We see in our simple example that one of them has degree 2, and larger graphs would lead to polynomials of higher linear degree. 
So, we could not find a lower bound, similar to Theorem \ref{theorem:disjoint}.
We therefore should find another argument to prove that $I_1/I_2$ is close to $1$. 
This result is given in the following Theorem as the second main result of the paper.

%If we follow the argument parallel to that of the previous section, assuming that we can find an inequality of the type \eqref{eq:bl}, instead of having to bound $\E \left( \left| N_pX_1\ldots X_{\ell-1} \right| \right) $, we will have to find an upper bound for $\E \left( \left| N_p \prod X_{ij}^{s_{ij}} \right| \right) $. 
%But, as proven in part 1 of Proposition \ref{proposition:lots}, the expected value of $E(X_{ij}^{s_{ij}})$ is finite only for $-1<s_{ij}<\delta+1$. In practice $\delta$ is taken to be equal to $3$. 
%So, in the preceding example, this expected value would be infinite and we could not find an upper bound.  We therefore have to find another argument to prove that $I_1/I_2$ is close to $1$. This result is given in Theorem \ref{theorem:joint} as the second main result of the paper.

\begin{theorem}
\label{theorem:joint}
Under the approximation $b_1 \sim N(0, v_D)$, the ratio $I_1/I_2$ \eqref{eq:i1i2} can be written
\begin{equation}
\label{eq:ratioJoint}
\frac{I_1}{I_2} = 
\frac{ \E \left( e^{-\frac{D}{2} } g\left( \delta^*, v_D \right) \right) }{ \E(e^{-\frac{D}{2}}) },
\end{equation}
where 
\begin{equation*}%\label{eq:fun-g}
g\left( \delta^*, v_D \right) = \frac{  \left( \frac{v_D}{2} \right)^{\delta^*}  }{ \Gamma( \delta^* ) }  \int_0^{\infty} t^{\delta^* - \frac{1}{2} } (1+t)^{\frac{-1}{2}} e^{-\frac{v_D t}{2}} dt,
\end{equation*}
in which $\delta^* = \frac{\delta+d}{2}$. 
Moreover, when $v_D$  is small, we have
\begin{equation*}
%\label{eq:jointEQ2}
 g\left( \delta^*, v_D \right) = 1 - \frac{\Gamma \left( \delta^* + \frac{1}{2} \right) }{ \Gamma \left( \delta^* \right)  } \left( \frac{v_D}{2} \right)^{\delta^*}  {\mathcal O} \left( \left|\frac{v_D}{2} \right|^{ \delta^* -1} \right).
\end{equation*}
And when, for all $D$, $v_D$ is uniformly bounded by a small quantity, we have
\begin{equation*}
%\label{eq:approxvsmall}
\frac{I_1}{I_2} = 1 - \frac{\Gamma \left( \delta^* + \frac{1}{2} \right) }{ \Gamma \left( \delta^* \right)  } \frac{\E \left(e^{-\frac{D}{2}} \left( \frac{v_D}{2} \right)^{\delta^*}  {\mathcal O}(|\frac{v_D}{2}|^{ \delta^* -1})\right)}{\E(e^{-\frac{D}{2}})}\
\approx 1\;
\end{equation*}
and it leads that our approximation \eqref{eq:approx} holds. 
%It also holds that $I_1/I_2$  always satisfies $\frac{I_1}{I_2} \leq 1$. 
\end{theorem}

\textbf{Proof.} %This theorem follows immediately form the Lemma \ref{lemma:main}.
%The proof is given in Section \ref{sup:proof-theorem2} of the Supplementary file.
The proof is in three steps. First, we show $b_1$ can be expressed as a bilinear form. %, condition on $\psi_{\tilde{D}_1}$ in \eqref{eq:newi1i2}. 
Then, using the bilinear expression, we prove $b_1$ is distributed as the continuous scale mixture of centered Gaussian variables. 
Finally, this allows us to deduce that there exists a unique $v_D$ so that the normal $N(0, v_D)$ distribution best approximates the $b_1$ distribution. 
For detailed proof see Section \ref{sup:proof-theorem2} of the Supplementary file.

In Theorem \ref{theorem:joint}, we prove that $I_1/I_2$ can accurately be approximated by 1, under the assumption that $v_D$ is small, 
%or equivalently the ratio in \eqref{eq:ratio_pro} can accurately be approximated by \eqref{eq:approx}. 
or equivalently our approximation in Equation \ref{eq:approx} is accurate.
The validity of the assumption that $v_D$ is small and the accuracy of the approximation is demonstrated numerically in the following subsection.

%In Subsection \ref{subsec:bilinear}, we show in Proposition \ref{proposition:quadratic} that, conditional on a quantity $\psi_{\tilde{D}_1}$ as defined in \eqref{eq:newi1i2}, $b_1$ can be expressed as a bilinear form. 
%In Proposition \ref{proposition:mixture} of Subsection \ref{subsec:mixture}, using its expression as a bilinear form, we show that $b_1$ is distributed like the continuous scale mixture of centered Gaussian variables. 
%This allows us to deduce that there exists a unique $v_D$ such that the normal $N(0, v_D)$ distribution best approximates the distribution of $b_1$. 
%Finally, in Subsection \ref{subsec:intunderapprox}, we state and prove our  second main result, Theorem \ref{theorem:joint}: under the assumption that $v_D$ is small, $I_1/I_2$ can accurately be approximated by 1 or equivalently the ratio in \eqref{eq:ratio_pro} can 
%accurately be approximated by \eqref{eq:approx}, which is what we want to prove. The validity of the assumption that $v_D$ is small and the accuracy of the approximation for different kinds of graphs of different sizes will be demonstrated numerically in Section  \ref{sec:simulation}.

% - - - - - - - - - - - - - - - - - - - - - - - - - - - - - - - - - - - - - - - - - - - - - - - - - - - - - - - - - - - -|
\subsection{Simulated experiments for the general case}
\label{subsec:sim-joint}

We compute the ratio $I_1/I_2$ in two different ways, first following the MC approach of \cite{atay2005monte} and second using our approximation in Theorem \ref{theorem:joint}; 
We call these values  $I_{1/2,\mbox{MC}}$ and $I_{1/2,\mbox{Gamm}}$, respectively. 
We note that, if our approximation $I_1/I_2\approx 1$ is good, without any additional conditions, $I_{1/2,\mbox{MC}}$ should reflect that by being close to 1. 
However, if our approximation $I_1/I_2\approx 1$ is good, according to Theorem \ref{theorem:joint}, $I_{1/2,\mbox{Gamm}}$ will be close to 1 if the assumption of $v_D$ small is satisfied.

%We will give the values of both  $I_{1/2,\mbox{MC}}$ and $I_{1/2,\mbox{Gamm}}$ for $9$ different graph structures: % (see Figure \ref{fig:plot_graphs_p30}):
%\begin{itemize}
%\item[1.] \textbf{Random\_1:}  A graph in which edges are randomly generated from independent Bernoulli distributions with probability $0.1$. It representative of sparse random graph.
%\item[2.] \textbf{Random\_2:}  The same as the Random\_1 graph for Bernoulli distributions with probability $0.2$. It representative of relatively sparse radnom graph.
%\item[3.] \textbf{Random\_5:}  The same as the Random\_1 graph for Bernoulli distributions with probability $0.5$. It representative of dense radnom graph.
%\item[4.] \textbf{Cluster:} A graph in which the number of clusters is $3$. Each cluster has the same structure as the \textbf{Random\_5} graph.
%\item[5.] \textbf{Scale-free:} A graph which has a power-law degree distribution generated by the Barab{\'a}si-Albert algorithm \citep{albert2002statistical}.
%\item[6.] \textbf{Star:} A graph in which every node is connected to the one node and that node is selected by random.
%\item[7.] \textbf{Hub:} A graph in which every node is connected to the $3$ nodes and those nodes are selected by random.
%\item[8.] \textbf{Circle:} See Figure \ref{fig:plot_graphs_p30} for the case $p=30$.
%\item[9.] \textbf{Lattice:} See Figure \ref{fig:plot_graphs_p30} for the case $p=25$.
%\end{itemize}
% - - - - - - - - - - - - - - - - - - - - - - - - - - - - - - - - - - - - - - - - - - - - - - - - - - - - - - - - - - - -|
\begin{figure} [!ht]
\begin{center}
    \includegraphics[width=11.5cm,height=8.5cm]{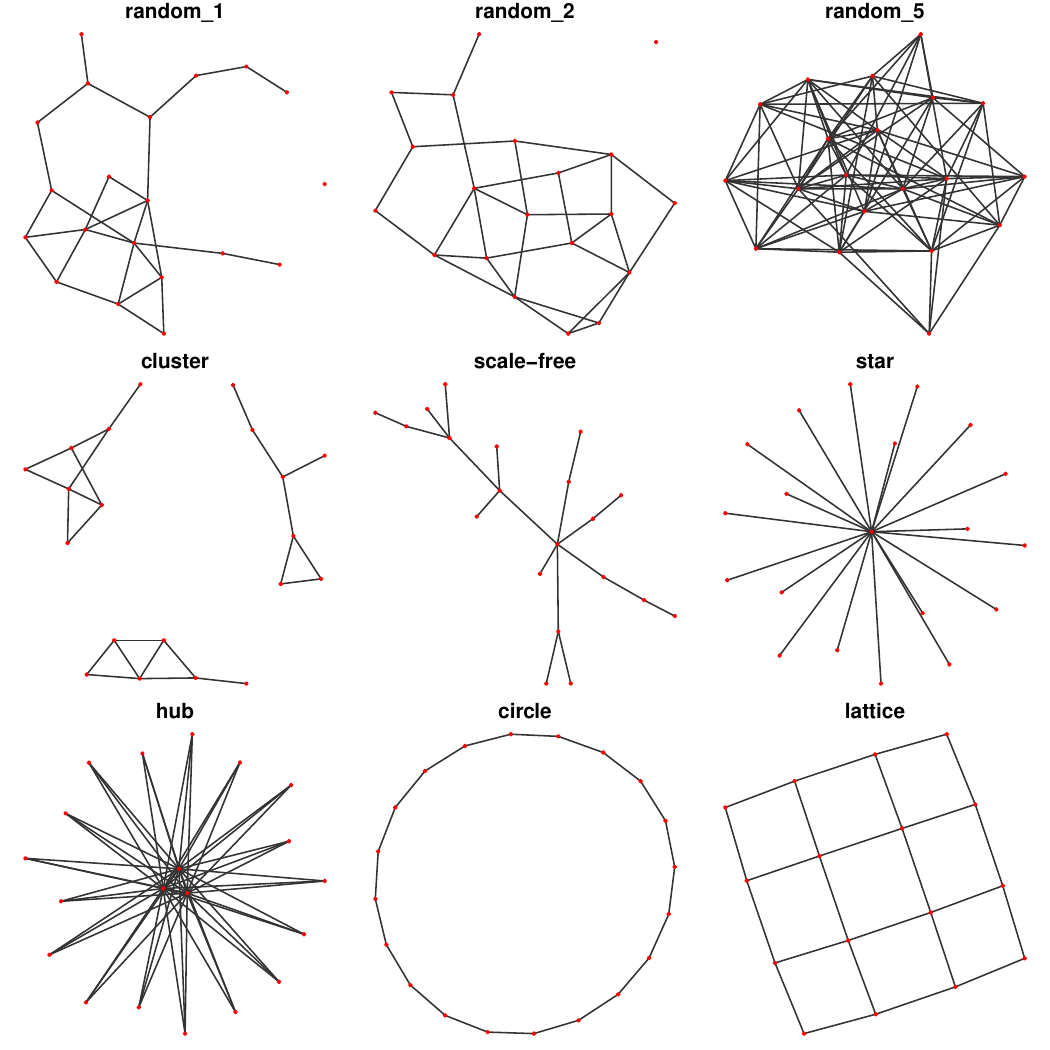}
\end{center}
\caption{ The $9$ different types of undirected graphs for $p=20$, as a number of nodes. For the case graph is \textit{Lattice} $p=16$. The graphs \textit{Random\_1},  \textit{Random\_2},  \textit{Random\_5} are random graphs with  edge probabilites equal to $0.1$, $0.2$, and $0.5$ respectively ranging from sparse to dense graphs. \label{fig:plot_graphs_p30} }
\end{figure}
% - - - - - - - - - - - - - - - - - - - - - - - - - - - - - - - - - - - - - - - - - - - - - - - - - - - - - - - - - - - -|

While it is straightforward to evaluate $I_{1/2,\mbox{MC}}$, it is less obvious how to compute $I_{1/2,\mbox{Gamm}}$ using Equation \ref{eq:ratioJoint}. 
% We proceed as follows.
%We have just proved that  $v_D$  close to $0$ is a sufficient condition for   approximation \eqref{approx} to hold with  accuracy proportional to $|\frac{v_D}{2}|^{\frac{\delta+d}{2}-1}$.
% This is not to say, though, that it is a necessary condition. Indeed, the gray boxes  in the top graph of Figure \ref{fig:boxplot disjoint} give the boxplots for the evaluations of $1-I_1/I_2$ obtained using its expression in \eqref{eq:newi1i2}.
The pseudo-code for evaluating the $I_{1/2,\mbox{Gamm}}$ is given in Section \ref{sec:pseudo-code} of the Supplementary file. 
%We see that the value of the relative error $1-I_1/I_2$ using  \eqref{eq:newi1i2} is rather large while we know from Theorem \ref{theorem:disjoint} and from the boxplots of Figure \ref{fig:boxplot} (also the white boxes in the top graph of Figure \ref{fig:boxplot disjoint}), 
%where the simple Monte-Carlo method of \cite{atay2005monte} has been used, that  in fact approximation \eqref{approx} is good.
%From the bottom graph in Figure \ref{fig:boxplot disjoint}, we also verify that, following \eqref{eq:ufunction}, when $v_D$ is relatively large, for example for graphs with 4 paths of length 2 between $q$ and $p$, the approximation to $I_1/I_2$ given by \eqref{eq:newi1i2} 
%is poor while when $v_D$ is smaller, for example when there are no paths of length 2, then this approximation becomes good.
We represent the boxplot of the numerical values of $I_{1/2,\mbox{MC}}$ and $I_{1/2,\mbox{Gamm}}$ obtained over $100$ replications for nine different types of graphs (Figure \ref{fig:plot_graphs_p30}) along with three different numbers of nodes $p=\{10,20,30\}$ and two different values for $\delta = \{3,10\}$. 
Besides, we report the corresponding values of $v_D$ so that one can see the variation of the accuracy of $I_{1/2,\mbox{Gamm}}$, as $v_D$ varies, as predicted by Theorem \ref{theorem:joint}, but also that of  $I_{1/2,\mbox{MC}}$.
 
\begin{figure} % [!ht]
\begin{center}
     \includegraphics[width=18cm]{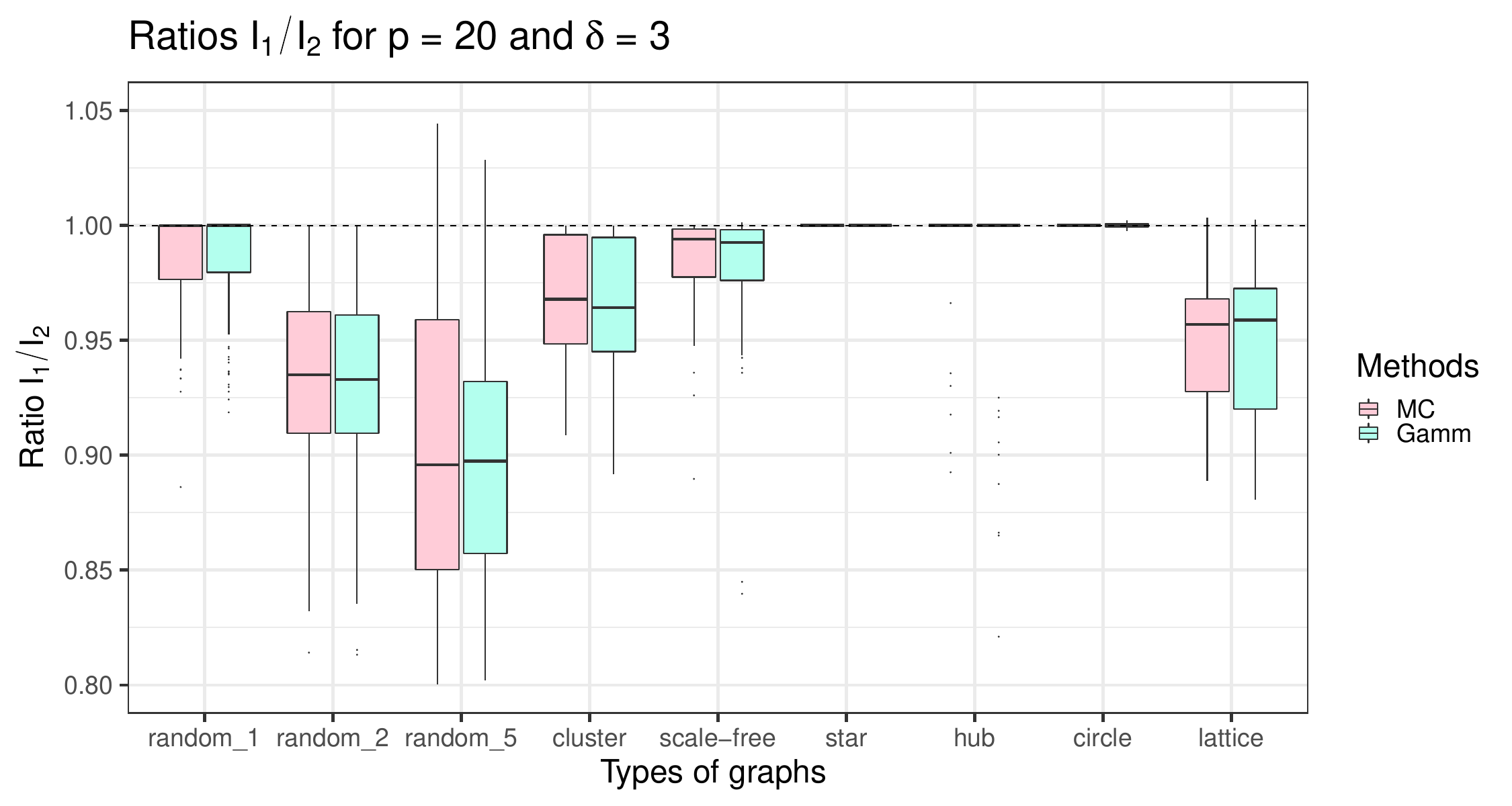}
     \includegraphics[width=15cm]{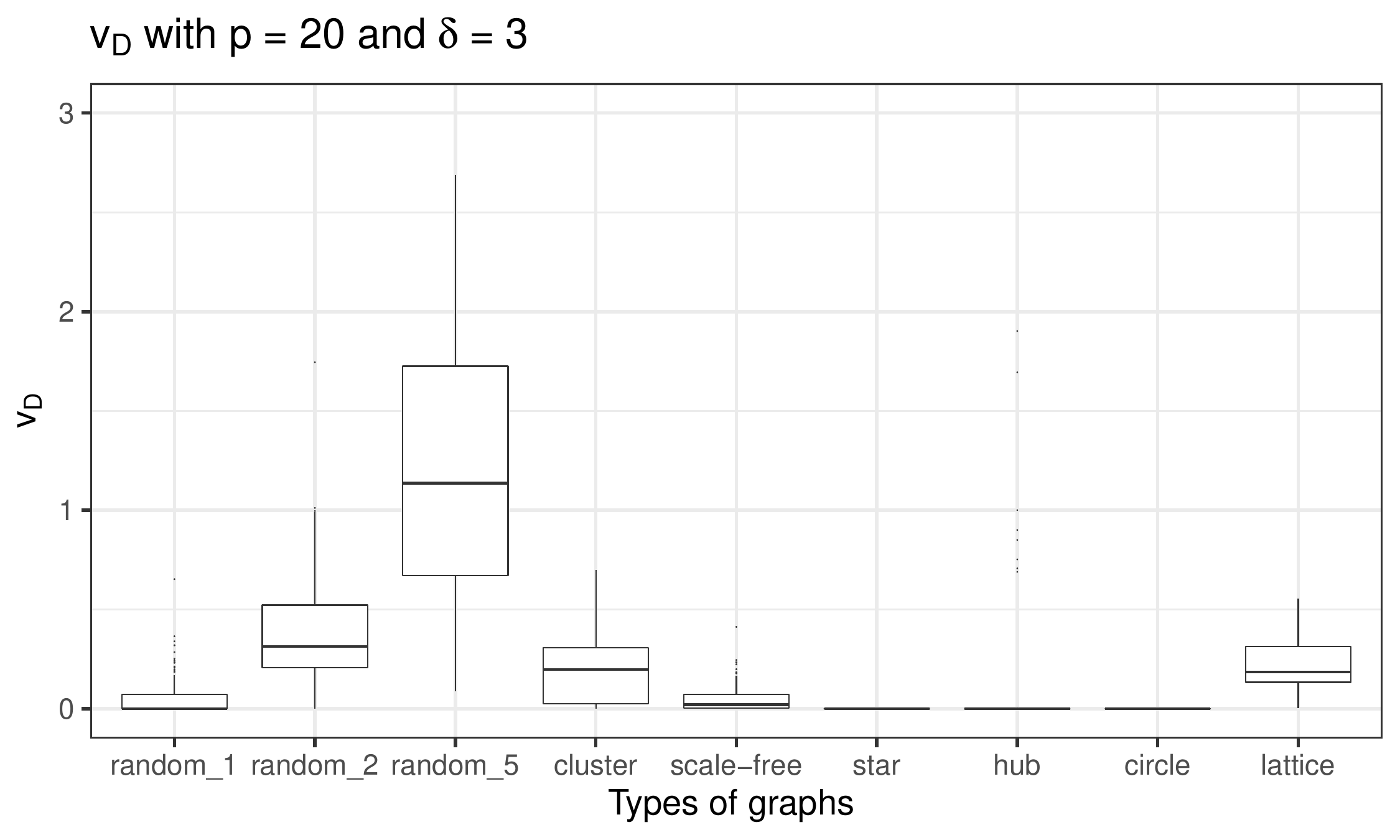}
\end{center}
\caption{ (Top) The boxplot for the ratio $I_1/I_2$ computed by the MC approach of \cite{atay2005monte} (in red) and our approximation \eqref{eq:ratioJoint} (in blue). (Bottom) The boxplot of the variance $v_D$ of $b_1$ for the corresponding graphs. These computations are done over 100 replications for nine different graphs (Figure \ref{fig:plot_graphs_p30}) with $p=20$ nodes and $\delta=3$. \label{fig:boxplot_p20} }
\end{figure}

%For the general case of  the nine different types of graphs,  
For the case $\delta = 3$,
the values of $I_{1/2,\mbox{MC}}$ and $I_{1/2,\mbox{Gamm}}$ are represented in Figure \ref{fig:boxplot_p20} for $p=20$, and Figures \ref{fig:boxplot_p10} and \ref{fig:boxplot_p30} in Section \ref{sec:plots} of the Supplementary File  for $p=\{10, 30 \}$. 
We see that the values of $I_1/I_2$ slightly move away from 1 as $v_D$ moves away from 0. 
But in all cases, we see that $I_{1/2,\mbox{MC}}$ and $I_{1/2,\mbox{Gamm}}$ cover the same range of values and their medians are between $0.9$ and 1, giving relative errors less than $0.10$. 
While, from these facts, we cannot immediately conclude that the assumption of $v_D$ small is always satisfied, it is a strong indication that it is satisfied enough to ensure that our approximation is acceptable.
% - - - - - - - - - - - - - - - - - - - - - - - - - - - - - - - - - - - - - - - - - - - - - - - - - - - - - - - - - - - -|
% - - - - - - - - - - - - - - - - - - - - - - - - - - - - - - - - - - - - - - - - - - - - - - - - - - - - - - - - - - - -|

\begin{figure} %[!ht]
\begin{center}
   \includegraphics[width=18cm]{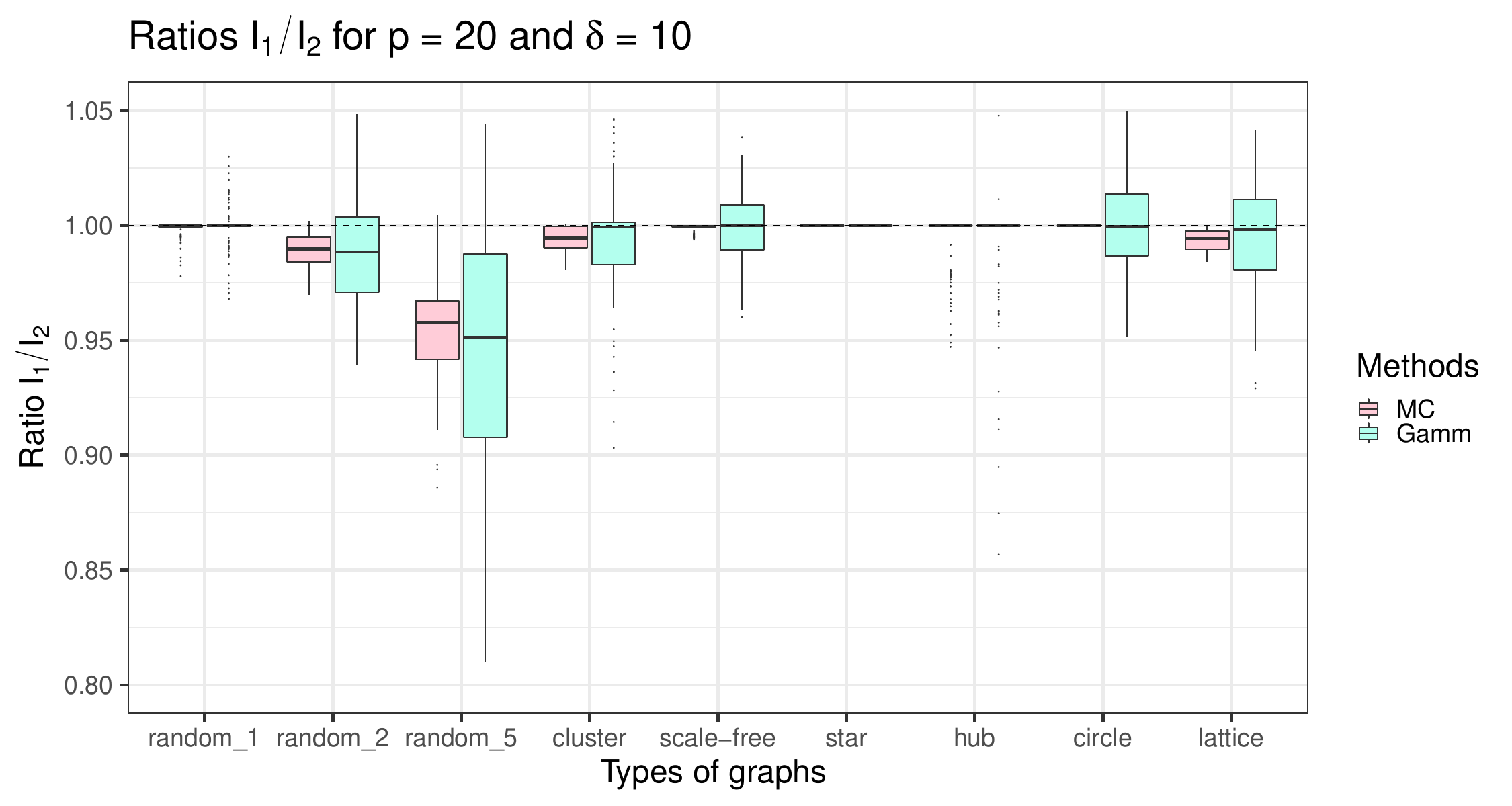}
    \includegraphics[width=15cm]{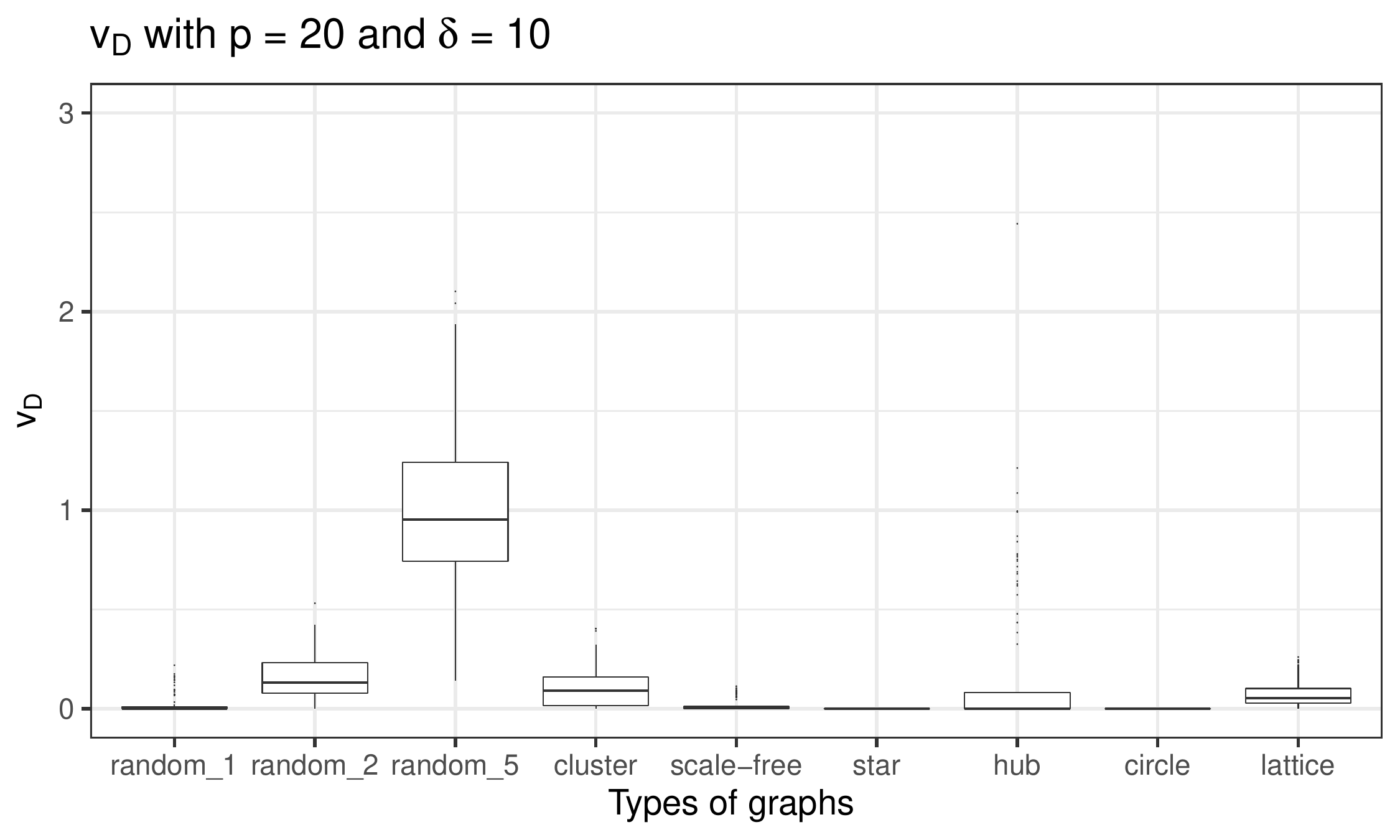}
\end{center}
 \caption{ (Top) The boxplot for the ratio $I_1/I_2$ computed by the MC approach of \cite{atay2005monte} (in red) and our approximation \eqref{eq:ratioJoint} (in blue). (Bottom) The boxplot of the variance $v_D$ of $b_1$ for the corresponding graphs. These computations are done over 100 replications for nine different graphs (Figure \ref{fig:plot_graphs_p30}) with $p=20$ nodes and $\delta=10$. \label{fig:boxplot_p20_delta10}}
\end{figure}

%For the general case of  the nine different types of graphs,  
For the case $\delta = 10$,
the values of $I_{1/2,\mbox{MC}}$ and $I_{1/2,\mbox{Gamm}}$ are represented in Figure \ref{fig:boxplot_p20_delta10} for $p=20$, and Figures \ref{fig:boxplot_p10_delta10} and \ref{fig:boxplot_p30_delta10} for $p=\{10, 30 \}$ in Section \ref{sec:plots} of the Supplementary File. 
%We see that the values of $I_1/I_2$ slightly move away from 1 as $v_D$ moves away from 0. 
In all cases, we see that $I_{1/2,\mbox{MC}}$ and $I_{1/2,\mbox{Gamm}}$ cover the same range of values and their medians are between $0.995$ and 1, giving pretty low relative errors of less than $0.005$. 
%While, from these facts, we cannot immediately conclude that the assumption of $v_D$ small is always satisfied, it is a strong indication that it is satisfied enough to ensure that our approximation is acceptable.
% - - - - - - - - - - - - - - - - - - - - - - - - - - - - - - - - - - - - - - - - - - - - - - - - - - - - - - - - - - - -|
% - - - - - - - - - - - - - - - - - - - - - - - - - - - - - - - - - - - - - - - - - - - - - - - - - - - - - - - - - - - -|

% to check
We verify this result numerically. 
Besides, our numerical results show our approximation $I_1/I_2$ (given by $I_{1/2,\mbox{Gamm}}$) is accurate (close to 1) even for the cases that $v_D$'s are not close to $0$.
%%If $v_D$ is not close to $0$, we see that numerically still, the approximation $I_1/I_2$, given by the values of $I_{1/2,\mbox{Gamm}}$  is close to 1. 
In fact, both set of values for $I_{1/2,\mbox{MC}}$ and $I_{1/2,\mbox{Gamm}}$ seem to be affected by the size of $v_D$ but are reasonably close to 1, whatever the value of $v_D$.

We should mention that our simulations indicate that our approximation is more accurate for the sparser graphs. 
For example, in Figure \ref{fig:boxplot_p20} (top) consider the graphs \textit{Random\_1},  \textit{Random\_2},  and \textit{Random\_5} which are respectively ranging from sparse to dense graphs.
This figure as well as the other figures in this section indicate that our approximation is more accurate for the sparser graphs. 
%Other figures in this Section confirm the same pattern. 

%So, for general graphs with up to $p=30$ nodes, we see that approximation \eqref{approx} is such that the relative error is of the order 0.10, which given the randomness of alternative methods such as double Metropolis-Hastings, is a reasonable accurate approximation.
For $p>30$, we cannot verify the accuracy of our approximation directly by computing $I_{1/2,\mbox{MC}}$ and $I_{1/2,\mbox{Gamm}}$
because of the limitations of the Monte Carlo method of \cite{atay2005monte}. 
So, in the next section, for graphs with up to $150$ nodes, we will verify the performance of our approximation in the search algorithm that represents in Section \ref{subsec:BD}. 

% - - - - - - - - - - - - - - - - - - - - - - - - - - - - - - - - - - - - - - - - - - - - - - - - - - - - - - - - - - - -|
\section{Simulation study for high-dimensional graphs}
\label{sec:simulation}

We perform Bayesian structure learning on simulated data from high-dimensional graphs using the BDMCMC search algorithm, %developed by \cite{mohammadi2015bayesianStructure} 
 represented in Algorithm \ref{alg:BDMCMC}. 
%We use our approximation ratio \eqref{eq:approx}, thus bypassing the lengthy Monte Carlo computations, we call it here BDMCMC-Gamm.
We use our approximation \eqref{eq:approx} within Algorithm \ref{alg:BDMCMC} and we call it BDMCMC-Gamm.
For the sake of comparison, we also evaluate the ratio of normalizing constants, within the BDMCMC search algorithm, using the exchange algorithm which is represented in Algorithm \ref{alg:BDMCMC-DMHJ}, we call it BDMCMC-DMH; This algorithm can be considered as the state-of-the-art. 
Both approaches are implemented in the \textbf{BDgraph} \textit{R} package \citep{BDgraph, mohammadi2019bdgraph} in the function \textit{bdgraph()}.

%paper we approximate the ratio of prior normalizing constants  using  ratio \eqref{approx} which is simple and fast to compute. 
%In the next subsection we show the perfomance of using our approximation within the BDMCMC algorithm.

%\subsection{Simulation experiments}
%\label{comparison}

%Here, we will perform Bayesian structure learning on high-dimensional simulated data using the BDMCMC approach  developed by \cite{mohammadi2015bayesianStructure}. Ratio \eqref{eq:ratio}  is central to all the computations in this method and an approximation that bypasses the lengthy Monte Carlo computations is essential. 
%We bypass the computation of this ratio by replacing ratio \eqref{eq:ratio}  by approximation \eqref{approx} or by using the double Metropolis-Hastings.
%Using the ratio \eqref{eq:ratio}  within the BDMCMC aglorithm is implemented in the \textbf{BDgraph} \textit{R} package \citep{BDgraph} in the function \textit{bdgraph()}. 

%We consider four different graph structures (see Figure \ref{fig:plot_graphs} of the Supplementary file) as follows:
We consider four following graph structures: % (see Figure \ref{fig:plot_graphs} of the Supplementary file):
\begin{itemize}
\item[1.] \textit{Scale-free:} A graph which has a power-law degree distribution generated by the Barab{\'a}si-Albert algorithm 
\citep{albert2002statistical}.
\item[2.] \textit{Random\_p:}  A graph in which edges are randomly generated from independent Bernoulli distributions with mean equal to $p$.
\item[3.] \textit{Random\_2p:} The same as the Random\_p graph with mean equal to $2p$.
\item[4.] \textit{Cluster:} A graph in which the number of clusters is $\left| p/20 \right| $. Each cluster has the same structure as the \textit{Random\_p} graph.
\end{itemize}
For each graph, we consider various scenarios based on the number of nodes $p \in \{ 50, 100, 150 \}$ and the sample size $n \in \{ p, 2p \}$.
We draw $n$ independent samples from the normal $N_p(0,K)$ distribution.
We consider $\delta=3$, which is the worst-case value for our approximation (see subsections \ref{subsec:sim-disjoint} and \ref{subsec:sim-joint}). 
%All simulations are performed using the \textbf{BDgraph} \textit{R} package \citep{BDgraph}, where approximation \eqref{approx} and double Metropolis-Hastings have been implemented. See also \citet{mohammadi2017bdgraph}.

For each scenario, we run Algorithm \ref{alg:BDMCMC} by using our approximation \eqref{eq:approx} as well as Algorithm \ref{alg:BDMCMC-DMHJ} which is based on an exchange algorithm.
The number of iterations is $100,000$ with $60,000$ iterations as burn-in.
To evaluate the performance of  both algorithms we use  ROC curves, based on model averaging, by computing true and false-positive rates for each of $50$ replicated data sets and then by averaging over the $50$ replicates.

\begin{figure} %[!ht]
\begin{center}
    \includegraphics[width=10cm]{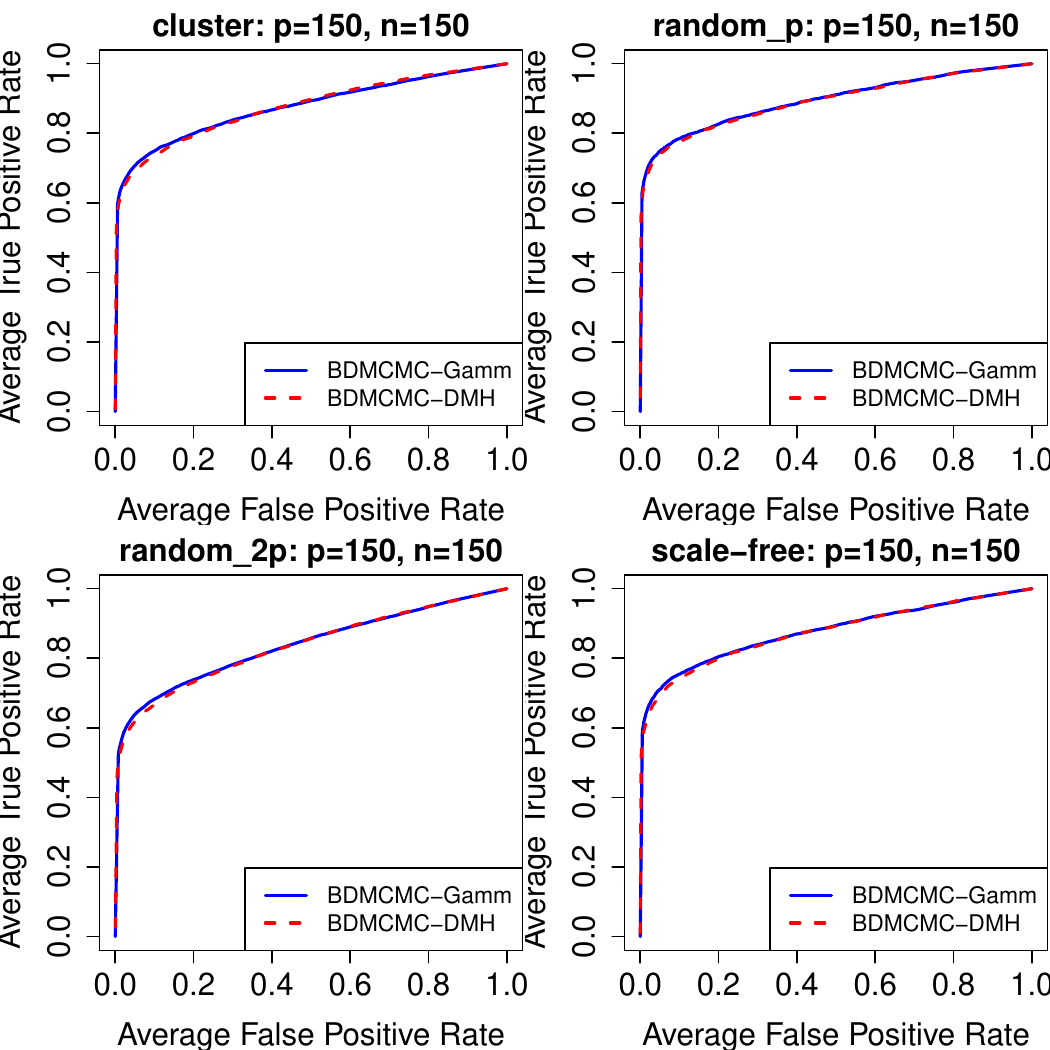} 
    \includegraphics[width=10cm]{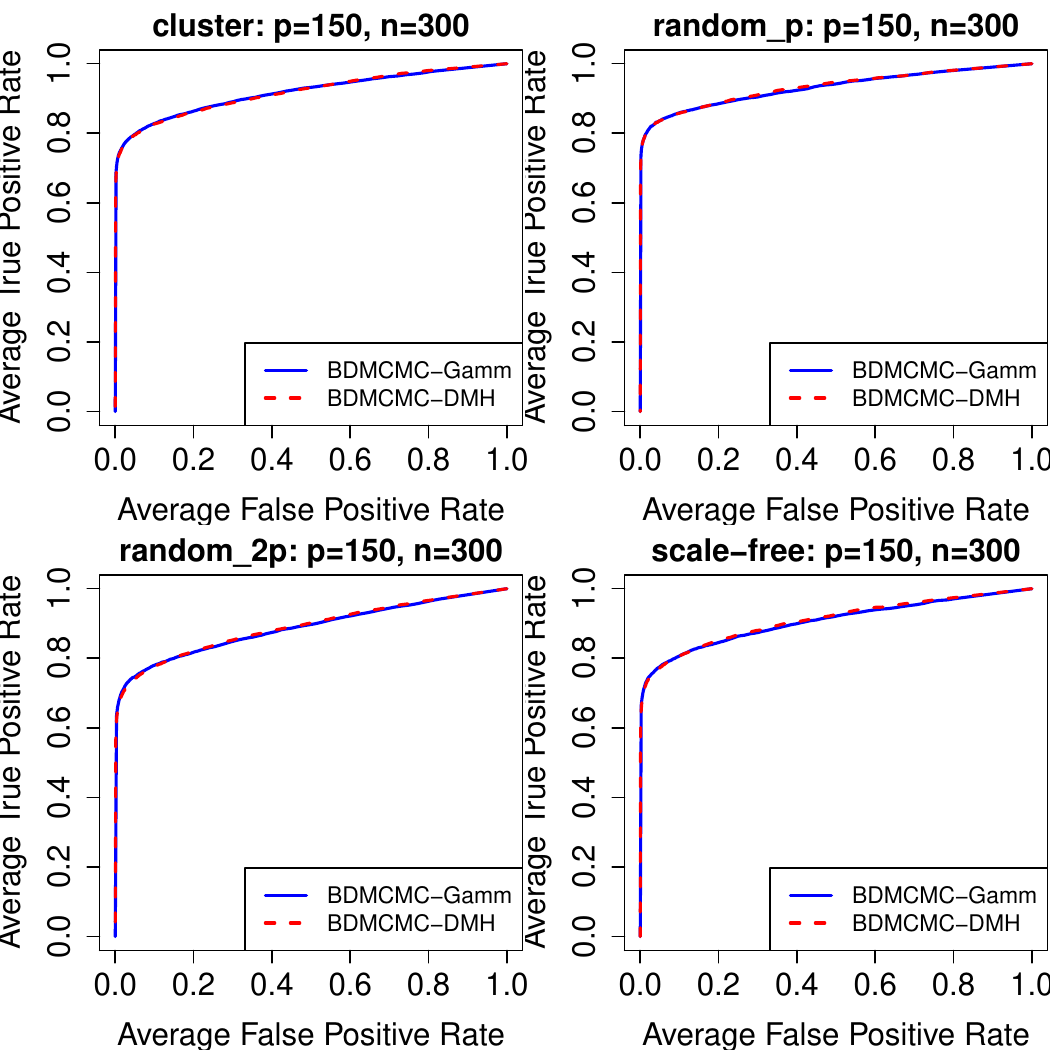}
\end{center}
\caption{ ROC curves for the BDMCMC algorithm with our approximation \eqref{eq:approx} (BDMCMC-Gamm) and BDMCMC algorithm with exchange algorithm (BDMCMC-DMH), over $50$ replications. 
Here, $p=150$, $n \in \{150, 300 \}$, and $4$ different graph structures. \label{fig:Rocplot_p=150} }
\end{figure}

Figure \ref{fig:Rocplot_p=150} represents the ROC curves for the cases $p=150$ with $n \in \{150, 300\}$. 
The ROC curves for $p=50$ and $p=100$ are, respectively, in Figures \ref{fig:Rocplot_p=50} and \ref{fig:Rocplot_p=100} in Section \ref{sec:plots} of the Supplementary File.
As we can see, in almost all cases, the performance of the BDMCMC algorithm based on both approximations is the same. 
In a few cases, the BDMCMC algorithm using our approximation \eqref{eq:approx} performs slightly better than the BDMCMC algorithm using the exchange algorithm: this happens especially when $p$ is large, for example, when $p=150$ and $n=150$. 
This discrepancy can be due to the convergence issue of the exchange algorithm in high-dimensional graphs. %(see \citealt{liang2010double}). 

The execution times for both algorithms are represented on the right-hand side of Figure \ref{fig:Rocplot_time}. 
It indicates the computational gain of using our approximation within the search algorithm. 
For example, in the case $p=150$, the BDMCMC algorithm using our approximation is more than $3$ times faster than the BDMCMC algorithm using the exchange algorithm. 

In summary, our simulation study shows that, from an accuracy point of view, the BDMCMC algorithm using our approximation \eqref{eq:approx}, performs well especially for high-dimensional sparse graphs, which is the case for many real-world applications. % as, for example in genetics. 
From a computational point of view, using our approximation speeds up the BDMCMC search algorithm for the models with high-dimensional graphs.

% - - - - - - - - - - - - - - - - - - - - - - - - - - - - - - - - - - - - - - - - - - - - - - - - - - - - - - - - - - - -|
% - - - - - - - - - - - - - - - - - - - - - - - - - - - - - - - - - - - - - - - - - - - - - - - - - - - - - - - - - - - -|

% - - - - - - - - - - - - - - - - - - - - - - - - - - - - - - - - - - - - - - - - - - - - - - - - - - - - - - - - - - - -|

\section{Conclusion}
\label{sec:conclusion}

%In this paper, we propose approximation \eqref{eq:approx} to the ratio \eqref{eq:ratio_pro} of prior normalizing constants.  This allows for model selection without the heavy computational burden of Monte Carlo methods such as the double Metropolis-Hastings algorithm. 
%This approximation consists in ignoring the paths of length strictly greater than two between the end points of the edge $e$ in ratio \eqref{eq:ratio_pro}. As a consequence, the approximation is exact when,  in $G^{-e}$,  there are only paths of length two between these endpoints. 
%As recalled earlier, \citet{uhler2014exact} have proved that if $G^{-e}$ is such that $G$ is decomposable, then equality holds in \eqref{approx}, that is the approximation is exact. 

%In this paper, we propose the approximation in Equation \ref{eq:approx} for the prior normalizing constants of G-Wishart. 
%This allows for Bayesian structure learning to avoid the sampling-based methods as computationally expensive updates within the search algorithm.
In this paper, we represent a search algorithm in which the the prior normalizing constants of G-Wishart is carried out by our approximation in Equation \ref{eq:approx}.
Using our approximation allows for Bayesian structure learning to avoid the sampling-based methods as computationally expensive updates within the search algorithm.
% Monte Carlo approximation such as the double Metropolis-Hastings algorithm.
%Hitherto, authors have primarily used exchange algorithm within the search algorithm to avoid evaluating the prior normalizing constants. 
%Since these types of algorithms can be computationally expansive and probably slow in convergence, for the graph with higher than 100 nodes ($p=100$) require super computers. 
%In this paper, with our approximation we make it possible to apply Bayesian structure learning for the graph with $p>100$ without requiring supper computers with the same accuracy as the state-of-the-art. 
%Our approximation consists in ignoring the length paths strictly greater than two between the end points of the edge $e$. % in ratio \eqref{eq:ratio_pro}. 
%As a consequence, our approximation \eqref{eq:approx} is exact when,  in $G^{-e}$,  there are only paths with length less than three between these endpoints. 
We give theoretical results to justify this approximation when certain assumptions are satisfied. 
Then, as importantly, we show, through numerical experiments that the assumptions are reasonably satisfied and yield a good accuracy of the approximation. 
In Theorem \ref{theorem:disjoint}, we consider the specific case where the paths between the endpoints are disjoint. 
Though this case is unrealistic in practice of course, it is interesting because we can obtain an analytic lower bound to the ratio $I_1/I_2$, which is a function of $\delta$ and the number of paths and their length. 
We see that the actual accuracy is much better than  that given by the lower bound. % $B(\delta, d, \ell_{\lambda})$. 

In the realistic and general case where the paths  are not necessarily disjoint, we give an alternative expression in Theorem \ref{theorem:joint} for the ratio $I_1/I_2$, then an approximation to this expression. 
We show that when the variance $v_D$ is small, then the accuracy is good. 
When performing structure learning in practice, one will not verify this assumption any more than one would verify that the paths are disjoint. 
But we do examine a large array of standard graphs and verify numerically that the assumption of $v_D$ small is satisfied in most cases.
Whatever the value of $v_D$, the accuracy of the approximation $I_1/I_2\approx 1$, or equivalently of the approximation in Equation \ref{eq:approx}, is very good. 
We do so by direct computation for graphs of size $p \leq 30$. 
Due to the limitations of Monte Carlo method to compute $I_1/I_2$, we cannot perform these direct computations for $p>30$. 
In that case, we perform structure learning on graphical models with up to $150$ variables and obtain the good results of Section \ref{sec:simulation}. 
We should emphasize here that we stopped at $p=150$ because, beyond this size, the state-of-the-art algorithms become computationally expensive but the BDMCMC search algorithm with our approximation \eqref{eq:approx} can scale up to higher dimensions still.

The accuracy of our approximation \eqref{eq:approx} depends on (i) the value of $\delta$ (scale parameter of G-Wishart) (ii) the structure of the graphs, more specifically its sparsity.
We illustrate it in the simulations of Sections \ref{subsec:sim-disjoint} and \ref{subsec:sim-joint}. 
It also can be interpreted from Theorems \ref{theorem:disjoint} and \ref{theorem:joint}.
The accuracy of our approximation is increased by increasing the values of $\delta$. 
Thus, our recommendation in practice is to choose, preferably, the value of $\delta$ higher than $10$, to be on the save side.
For the case of graph structure, the accuracy of our approximation depends on the sparsity of the graphs.
Our results indicated that our approximation is more accurate for the sparser graphs, as is indicated in the simulations of Section \ref{subsec:sim-joint}.
Since in real-life applications the underlying graphs are not dense (mainly sparse) is safe to use our approximation in practice. 

In conclusion, we think that our approximation can be safely adopted in the search algorithm to replace the sampling-based methods such as the exchange algorithm.
Finally, we also  proved that $I_1/I_2 \leq 1$.  
It shows that our approximation \eqref{eq:approx} yields a Bayes factor which favours $G^{-e}$ compared to $G$ so that we know that a model search using our approximation might lead to  sparser graphs.
%, a point we may want to correct with the help of an adequate prior on the space of graphs. 

% - - - - - - - - - - - - - - - - - - - - - - - - - - - - - - - - - - - - - - - - - - - - - - - - - - - - - - - - - - - -|
%\section*{Acknowledgments}

%The authors have contributed equally to this work and are listed in alphabetical order. 
%The authors thank the Editor, Associate Editor and two referees for their helpful comments that have resulted in significant improvements of the article.
%We are grateful to two referees and an associate editor for their constructive comments and helpful suggestions.

%We are grateful to the Editor, Dr. Wang, and two referees for their insightful comments that have resulted in significant improvements to the article.

% - - - - - - - - - - - - - - - - - - - - - - - - - - - - - - - - - - - - - - - - - - - - - - - - - - - - - - - - - - - -|
\bigskip
% - - - - - - - - - - - - - - - - - - - - - - - - - - - - - - - - - - - - - - - - - - - - - - - - - - - - - - - - - - - -|
\begin{center}
{\large\bf SUPPLEMENTARY MATERIALS}
\end{center}

%Supplementary materials contains the proofs of Theorem \ref{theorem:disjoint} and \ref{theorem:joint}, Lemma \ref{lemma:I_1}, Lemma \ref{lemma:psie}, some general results, and additional simulation results. 

The supplementary materials contain technical proofs for all the theorems from the main article as well as additional simulation results.
%Simulation results and detailed proofs are provided in the supplementary materials.

%The following seven sections comprise the online supplementary materials: A. Proposition for our approximation, B. Proof of Lemma \ref{lemma:I_1}, C. Proof of Lemma \ref{lemma:psie}, D. Proof of Theorem \ref{theorem:disjoint}, E. Proof of Theorem \ref{theorem:joint}, F. Pseudo-code to evaluate expression in Theorem \ref{theorem:joint}, H. Additional simulation results.

% - - - - - - - - - - - - - - - - - - - - - - - - - - - - - - - - - - - - - - - - - - - - - - - - - - - - - - - - - - - -|
%\begin{description}
%\item[Appendix:] Theorems, detailed proofs, and additional simulation results. (.zip file)

%\item[R-package:] R-package ÒBDgraphÓ containing code to perform the diagnostic methods described in the article. To download, please visit: https://CRAN.R-project.org/package=BDgraph. (R)
%\end{description}
% - - - - - - - - - - - - - - - - - - - - - - - - - - - - - - - - - - - - - - - - - - - - - - - - - - - - - - - - - - - -|

\appendix

% - - - - - - - - - - - - - - - - - - - - - - - - - - - - - - - - - - - - - - - - - - - - - - - - - - - - - - - - - - - -|
%\section{Appendix }

% - - - - - - - - - - - - - - - - - - - - - - - - - - - - - - - - - - - - - - - - - - - - - - - - - - - - - - - - - - - -|
\section{Proposition}
\label{sec:proposition EA}

The following proposition is used to compute $\E \left( e^{-\frac{A^2}{2}} \right)$ where $A$ defined in Equation \ref{eq:A} of the manuscript.
\begin{proposition}
\label{proposition:EA}
Let $U_1,\ldots,U_k$, $V_1\dots,V_k$, and $Q$ be independent random variables such that $U_i $ and $V_i$ are $N(0,1)$ and $Q \sim \chi^2_{\delta} $. 
Then
\begin{equation*}
%\label{eq:EA}
E \left( e^{-\frac{1}{2Q} \left( \sum_{i=1}^k U_iV_i \right )^2 } \right) =
\frac{\Gamma \left( \frac{\delta+k}{2} \right) \Gamma \left( \frac{\delta+1}{2} \right) }{ \Gamma \left( \frac{\delta}{2} \right) \Gamma \left( \frac{\delta+k+1}{2} \right) }.
\end{equation*}
\end{proposition}

\textit{Proof}. We have 
\begin{equation}
\label{eq:E3}
\sum_{i=1}^k U_i V_i \sim \left( U_1^2 + \cdots + U_k^2 \right)^{\frac{1}{2}} V_1 = X^{\frac{1}{2}} V_1,
\end{equation}
where $X \sim \Gamma( \frac{k}{2}, \frac{1}{2} )$;
To see this, we compute the Laplace transforms of both sides % of \eqref{eq:E3}. Thus we have
as follow
\begin{equation*}
E \left( e^{-\frac{1}{2Q} \left( \sum_{i=1}^k U_iV_i \right )^2 } \right) 
= E \left( e^{ \frac{X}{Q} \frac{-V_1}{2} } \right) 
= E \left( \left( 1 + \frac{X}{Q} \right)^{\frac{-1}{2}} \right),
\end{equation*}
where the last equality is due to integrating with regard to variable $V_1$.  
Since $X \sim \Gamma( \frac{k}{2}, \frac{1}{2} )$ and $Q \sim \Gamma( \frac{d}{2}, \frac{1}{2} )$, we have $U = \frac{X}{Q} \sim B_2 \left( \frac{k}{2}, \frac{d}{2} \right)$ which is a Beta distribution of second kind. Thus
\begin{equation*}
E \left( \left( 1 + \frac{X}{Q} \right)^{\frac{-1}{2}} \right) = \frac{\Gamma \left( \frac{\delta+k}{2} \right) \Gamma \left( \frac{\delta+1}{2} \right) }{ \Gamma \left( \frac{\delta}{2} \right) \Gamma \left( \frac{\delta+k+1}{2} \right) }.
\end{equation*}

% - - - - - - - - - - - - - - - - - - - - - - - - - - - - - - - - - - - - - - - - - - - - - - - - - - - - - - - - - - - -|

\section{Proof of Lemma \ref{lemma:I_1} }
\label{sec:lemma:I_1}

%We first define the conditioning set
%\begin{equation*} 
% \label{firstD}
%\Psi_{E_{qp}^-} = \Psi_{E_{q \cap p}} = 
%\psi_{\tilde{D}} = \left\{ \psi_{ij}: (i, j) \in E \setminus \left( N_q \cap N_p \right) \right \},
%\end{equation*}
%We have

% $\psi_{\tilde{D}}$ --> $\Psi_{\cap}^-$
By considering $\Psi_{\cap}^- = \left\{ \psi_{ij}: (i, j) \in E \setminus \left( E_q \cap E_p \right) \right \}$, we have
\begin{eqnarray*}
\E \left( e^{-\frac{ D}{2} - \frac{(A+b)^2}{2}} \right) = \E \left( e^{-\frac{ D}{2}} \E \left( e^{-\frac{(A+b)^2}{2}} \big| \Psi_{\cap}^- \right) \right).
\end{eqnarray*}
Note that $D$ and $b_1$ are $\Psi_{\cap}^-$-measurable, that is, are functions of the elements of $\Psi_{\cap}^-$ only. 
Due to Equation \ref{eq:E3}, we have 
\begin{equation*} %\label{JD}
\E \left( e^{ \frac{-(A_1+b_1)^2 }{ 2 \psi_{qq} } } \big| \Psi_{\cap}^- \right) = \E \left( e^{-\d \frac{(UV+b_1)^2 }{ \psi_{qq} } } \big| \Psi_{\cap}^- \right)
\end{equation*}
where $U$, $V$, and $\psi_{qq}$ are independent random variables such that $U^2 \sim \chi^2_{d}$, $V \sim N(0,1)$, and $\psi_{qq} \sim \chi^2_{\delta}$.
Integrating with respect to $V\sim N(0,1)$, we obtain
\begin{align*}
\E \left( e^{-\d \frac{(UV+b_1)^2 }{ \psi_{qq} } } \big| \Psi_{\cap}^- \right) 
   & = \E \left( e^{-\d\frac{b_1^2 }{ U^2 + \psi_{qq} } } \left( \frac{ \psi_{qq} }{ U^2 + \psi_{qq} } \right)^{ \frac{1}{2} } \Big| \Psi_{\cap}^- \right) \\
   & = \E \left( \sqrt{B} \right) \E \left( e^{ \frac{-b_1^2}{2 Y}} \big| \Psi_{\cap}^- \right),   
\end{align*}
where $B = \frac{ \psi_{qq} }{ U^2 + \psi_{qq} } \sim Beta( {\frac{\delta}{2}, \frac{d}{2}} )$ and is independent of $Y = U^2 + \psi_{qq} \sim \chi^2_{\delta + d}$. Thus
\[
\E \left( \sqrt{B} \right) = \frac{\Gamma(\frac{\delta+1}{2})\Gamma(\frac{\delta+d}{2})}{\Gamma(\frac{\delta+d+1}{2})\Gamma(\frac{\delta}{2})},
\]
which is equal to $\E \left( e^{-\frac{A^2}{2}} \right)$. 
Since $Y \sim \chi_{\delta + d}^2$, we have
\[
\E \left( e^{\frac{-b_1^2}{2Y}} \big| \Psi_{\cap}^- \right) = \E \left( h(b_1, \delta^*) \mid \Psi_{\cap}^- \right).   
\]
Regarding that $\Psi_{\cup}^- \subset \Psi_{\cap}^-$ we have
\[
\E \left( h(b_1, \delta^*) \mid \Psi_{\cap}^- \right) = \E \left( h(b_1, \delta^*) \mid \Psi_{\cup}^- \right).   
\]
Note that, while $b_1$ is $\Psi_{\cap}^-$-measurable, it is not $\Psi_{\cup}^-$-measurable. 
%Next, we will show that the distribution of $b_1$ is a scale mixture of normal distributions which can be approximated by another $N(0, v_D)$ distribution where the variance $v_D$ depends on $\psi_{\tilde{D}_1}$. But, to do so, we must first express $b_1$ as a bilinear form in two standard normal random vectors.

% - - - - - - - - - - - - - - - - - - - - - - - - - - - - - - - - - - - - - - - - - - - - - - - - - - - - - - - - - - - -|

\section{Proof of Lemma \ref{lemma:psie} }
 \label{sec:lemma:psie}
 
We note three important facts.
First, the elements of the first row of the matrix $\psi$ are all zero except for those corresponding to the edges of the path $\lambda_1$, i.e.
\begin{equation*}
%\label{eq:fact1}
\psi_{1v} = 0, \;\; \text{for} \;\; v \in \cup_{\lambda \in \Lambda} V_{\lambda} \;\; \text{and} \;\; v \not = \{ 1, 2, q \}.
\end{equation*}
Second, based on the above fact and Equation \ref{eq:aliye} of the manuscript, the remaining non-free entries in all the columns of $\psi$ except for the columns $q$ and $p$, are equal to zero. %%, i.e., for $\lambda, \lambda'\in \Lambda$ and for $\lambda$ ranked before $\lambda'$,
%%\begin{equation}
%%\label{fact2}
%%\psi_{i_{\lambda},j_{\lambda'}}=0, \;1_{\lambda}<i_{\lambda}<j_{\lambda'},\;,j_{\lambda'}\not=q,p\;.
%%\end{equation}
Third, due to the first entry $\psi_{1q}$ of column $q$ being free, none of the entries of column $q$ are necessarily zero. 
However, for each $\lambda\in \Lambda$, using iteratively Equation \ref{eq:aliye} of the manuscript, we see the non-free entries of column $p$ are zero except for the last one $\psi_{qp}$. %% which is a free variable, i.e.
%%\begin{equation}
%%\label{fact3}
%%\psi_{j_{\lambda}, p}=0,\;\; j_{\lambda}<\ell_{\lambda}\;, \lambda\in \Lambda\;.
%%\end{equation}
Considering these important facts and applying Equation \ref{eq:aliye}  of the manuscript  yields
\begin{equation}
\label{eq:psia2}
\psi_{qp} = \frac{-1}{\psi_{qq}} \sum_{\lambda \in \Lambda} \sum_{ i_{\lambda \in V_{\Lambda}}} \psi_{i_{\lambda}q} \psi_{i_{\lambda}p} =  \frac{-1}{\psi_{qq}} \sum_{\lambda\in \Lambda}\psi_{\ell_{\lambda}q} \psi_{\ell_{\lambda}p}.
\end{equation}
The entries $ \psi_{\ell_{\lambda}p},\;\lambda\in \Lambda$ are free. The entries $\psi_{\ell_{\lambda}q}$ are obtained by successively applying Equation \ref{eq:aliye} of the manuscript and the fact that $\psi_{(j-1)_{\lambda} j_{\lambda}}$, $j =\{1,\ldots, (l-1) \}$ are free and the non-free entries of $\Psi$ are equal to zero expect for the columns $q$ and $p$. 
That is
\begin{equation*} %\label{psia3}
\begin{split}
\psi_{\ell_{\lambda}q} & = -\frac{\psi_{(\ell-1)_{\lambda} \ell_{\lambda}} \psi_{(\ell-1)_{\lambda}q} }{ \psi_{(\ell-1)_{\lambda} (\ell-1)_{\lambda}}} 
= +\frac{\psi_{(\ell-1)_{\lambda} \ell_{\lambda}}\psi_{(\ell-2)_{\lambda} \ell_{\lambda}}\psi_{(\ell-2)_{\lambda}q}}{\psi_{(l-1)_{\lambda} (l-1)_{\lambda}}\psi_{(l-2)_{\lambda} (l-2)_{\lambda}}} \nonumber \\
& = \ldots \\
 & = (-1)^{\ell_{\lambda}-1} \frac{\psi_{1_{\lambda}q}\prod_{j=1}^{\ell-1}\psi_{j_{\lambda} (j+1)_{\lambda}}}{\prod_{j=2}^{\ell}\psi_{j_{\lambda} j_{\lambda}}}
= (-1)^{\ell_{\lambda}-1} \frac{\psi_{1_{\lambda}q}\prod_{j=1}^{\ell-1}\psi_{j_{\lambda} (j+1)_{\lambda}}}{\prod_{j=1}^{\ell-1}\psi_{(j+1)_{\lambda}(j+1)_{\lambda}}}, \\
& = (-1)^{\ell_{\lambda}-1}\psi_{1_{\lambda}q}\prod_{j=1}^{l-1}\frac{\psi_{j_{\lambda} (j+1)_{\lambda}}}{\psi_{(j+1)_{\lambda}(j+1)_{\lambda}}}.
\end{split}
\end{equation*}
Above equality and Equation \ref{eq:psia2} together yield
\begin{equation*} % \label{psia4}
\psi_{qp} =  \frac{1}{\psi_{qq}} \sum_{\lambda\in \Lambda} (-1)^{\ell_{\lambda}} \frac{\psi_{1_{\lambda},q} \psi_{\ell_{\lambda},p} \prod_{j=1}^{\ell-1} \psi_{j_{\lambda}, (j+1)_{\lambda}}}{ \prod_{j=2}^{l}\psi_{j_{\lambda},j_{\lambda}}} 
= \frac{1}{\psi_{qq}} \sum_{\lambda\in \Lambda} (-1)^{\ell_{\lambda}}\psi_{\ell_{\lambda},p} \psi_{1_{\lambda},q} \prod_{j=1}^{\ell-1}\frac{\psi_{j_{\lambda}, (j+1)_{\lambda}}}{\psi_{(j+1)_{\lambda},(j+1)_{\lambda}}},
\end{equation*}
which is identical to Equation \ref{eq:psie} of the manuscript.

% - - - - - - - - - - - - - - - - - - - - - - - - - - - - - - - - - - - - - - - - - - - - - - - - - - - - - - - - - - - -|
%\section{Appendix }
\section{Proof of Theorem \ref{theorem:disjoint} }
\label{sup:proof-theorem1}
 
The proof relies on a fact that for the case where the paths between $q$ and $p$ are disjoint, we rewrite $A_1$ and $b_1$ in Equations \ref{eq:A} and \ref{eq:b} of the manuscript as 
\begin{align*}
A_1 = \sum_{\lambda\in \Lambda, l_{\lambda=1}}\psi_{\ell_{\lambda},q}\psi_{\ell_{\lambda},p},  \;\;\;
b_1 = \sum_{\lambda \in \Lambda, \ell_{\lambda} \geq 2} b_{1\lambda},
\end{align*}
where
\[
b_{1\lambda} = (-1)^{\ell_{\lambda}} \psi_{1_{\lambda},q}\psi_{\ell_{\lambda},p}\prod_{j_{\lambda}=1}^{\ell_{\lambda}-1}\frac{\psi_{j_{\lambda}, (j+1)_{\lambda}}}{\psi_{(j+1)_{\lambda},(j+1)_{\lambda}}}
\]
with the convention that $b_{1\lambda}=0$ if $\ell_{\lambda}=1$,
and
\begin{equation*}
D = \sum_{\lambda \in \Lambda}  D_{\lambda} \; \text{ where } \;\; 
D_{\lambda} = 
\sum_{k=2}^{\ell_{\lambda}} \left( (-1)^{k-1}
\psi_{1_{\lambda},q}\prod_{j_{\lambda}=1}^{k-1}\frac{\psi_{j_{\lambda}, (j+1)_{\lambda}}}{\psi_{(j+1)_{\lambda},(j+1)_{\lambda}}} \right) ^2.
\end{equation*}
All the entries  appearing in the expression for $A_1$ and $b_1$ are free variables independent of each other and those appearing in $b_{1\lambda}, \lambda\in \Lambda, \ell_{\lambda}\geq 2$ are different from those appearing in $A_1$. 
Thus, $A_1$ and $\sum_{\lambda\in \Lambda, \ell_{\lambda}\geq 2}b_{1\lambda}$ are stochastically independent. 
Moreover, according to Proposition \ref{roposition:atay}, all $\psi_{ij}, i\not= j$ are $N(0,1)$ random variables while $\psi_{ii}^2$ follow a $\chi^2_{\delta+\nu_i}$ distribution. In particular $\psi^2_{qq}\sim\chi^2_{\delta}$.  
  
To find a lower bound for the $I_1/I_2$, we use the Gaussian equality as follows:
%To prove \eqref{ineq:approx}, we thus have to find a lower bound for the $I_1/I_2$. In the sequel, we will often use the Gaussian equality: 
if $Z\sim N(0, \sigma^2)$, then
\begin{equation*} %\label{fourier}
\E \left( e^{ i t Z } \right) = \int_{-\infty}^{+\infty} e^{i t z} e^{-\frac{z^2}{2\sigma^2}}\frac{dz}{\sigma\sqrt{2\pi}}=e^{-\frac{\sigma^2t^2}{2}}.
\end{equation*}
Applying above equality with $t = b_1$ and $ \sigma^2 = \frac{1}{y}$, we have
\begin{eqnarray*}
h(b_1, \delta^*)
&=& \frac{ 2^{-\delta^*} }{ \Gamma( \delta^* )} \int_0^{+\infty} y^{\delta^* - 1} e^{ \frac{-y}{2} } e^{ \frac{-b_1^2}{2 y} }  dy \\
&=& \frac{ 2^{-\delta^*} }{ \Gamma( \delta^* )}  \int_0^{+\infty} y^{\delta^* - 1} e^{ \frac{-y}{2} } \left( \int_{-\infty}^{+\infty} e^{i b_1 z} e^{-\frac{y z^2}{2}} \sqrt{y} \frac{dz }{ \sqrt{2 \pi}} \right)  dy \\
&=& \frac{ 2^{-\delta^*} }{ \Gamma( \delta^* ) \sqrt{2 \pi}}  \int_{+\infty}^{+\infty} e^{i b_1 z} \left( \int_0^{+\infty} y^{\delta^* - 1/2 } e^{  \frac{-(1+z^2)y}{2} }  dy \right) dz \\
&=& \frac{ 2^{-\delta^*} }{ \Gamma( \delta^* ) \sqrt{2 \pi}} \int_{+\infty}^{+\infty} e^{i b_1 z} \left( \frac{\Gamma( \delta^* + 1/2 ) }{ \left( \frac{1+z^2}{2} \right)^{\delta^* + 1/2} } \right) dz \\
&=& \int_{+\infty}^{+\infty} e^{i b_1 z} f(z) dz,
\end{eqnarray*}
where 
\[
f(z) = \frac{ \Gamma( \delta^* + 1/2 )}{ \sqrt{\pi} \Gamma( \delta^*) } \left( 1+z^2 \right)^{ -\delta^* + 1/2 }.
\]
Thus
\begin{equation*}
\begin{split}
\E \left( e^{-\frac{D}{2}} \E \left( h(b_1, \delta^*) \mid \Psi_{\cup}^- \right)  \right) & = \int_{+\infty}^{+\infty} \E \left( e^{ -\frac{D}{2} + i b_1 z } \right) f(z) dz \\
 & = \prod_{\lambda\in \Lambda} \int_{+\infty}^{+\infty} \E \left( e^{ -\frac{D_\lambda}{2} + i b_{ 1 \lambda} z } \right) f(z) dz
\end{split}.
\end{equation*}
Similarly, we have
\begin{eqnarray*}
\E \left( e^{-\frac{D}{2}}  \right) =  \prod_{\lambda\in \Lambda} \E \left( e^{-\frac{D_{\lambda}}{2}} \right).
\end{eqnarray*}

Consider independent identically distributed random variables $X_1,\ldots,X_n,\ldots$  such that $X_1\sim Z/\sqrt{Q}$ with $Z\sim N(0,1)$ independent of $Q \sim \chi^2_{\delta+1}$. 
For $\ell=\ell_{\lambda}, \lambda\in \Lambda$, we define
\begin{equation} 
\label{eq:S_l}
S_{\ell}=X_1^2+X_1^2X_2^2+\cdots+(X_1\ldots X_{\ell-1})^2,\ \ B_{\ell}=X_1X_2\ldots X_{\ell-1}.
\end{equation} 
We see that for $\lambda\in \Lambda$, we have 
\begin{eqnarray*}
%\label{b1lambda}
D_{\lambda}&\sim&N_{1_{\lambda}q} S_{\ell_{\lambda}}, \nonumber\\
b_{1\lambda}&\sim&N_{1_{\lambda}q}N_{\ell_{\lambda}p}B_{\ell_{\lambda}},
\end{eqnarray*}
where $N_{1_{\lambda}q}$ and $N_{\ell_{\lambda}p}$ are  independent $N(0,1)$  random variables, independent of $X_1,\ldots,X_{\ell},\ldots$.
We note that, from the independence of the entries of $\psi_E$, we have
\[
(b_{1\lambda}, D_{\lambda}, N_{1_{\lambda}q}, N_{\ell_{\lambda}p}),\;\;\;\lambda\in \Lambda
\]
are mutually independent.

Omitting the index $\lambda$ on $\ell_{\lambda}$, and simplifying $N_{1_{\lambda}q}$ to $N_q$ and $N_{\ell_{\lambda}p}$ to $N_p$, we define
 \begin{equation*}
 %\label{eq:E6}
 g_{\ell}(x) = \E \left( e^{-\frac{N_q^2S_{\ell} }{ 2} + iN_pN_qB_{\ell}x} \right).
 \end{equation*}
 Then
\begin{eqnarray*}
I_1 = \prod_{ \lambda\in \Lambda } \int_{-\infty}^{\infty} g_{\ell_{\lambda}}(x)f(x)dx
\;\; \text{ and } \;\;
I_2 = \prod_{\lambda\in \Lambda} g_{\ell_{\lambda}}(0).
\end{eqnarray*}
Therefore we can write
\begin{equation} 
\label{eq:i2-i1}
\begin{split}
\frac{I_2 - I_1}{I_2} 
&= \int_{-\infty}^{\infty} \frac{\prod_{\lambda \in \Lambda} g_{\ell_{\lambda}}(0) - \prod_{\lambda\in \Lambda} g_{\ell_{\lambda}}(x)}{ \prod_{\lambda\in \Lambda} g_{\ell_{\lambda}}(0) } f(x) dx \\
& \leq \int_{-\infty}^{\infty} \sum_{\lambda\in \Lambda} \frac{ g_{\ell_{\lambda}}(0)-g_{\ell_{\lambda}}(x) }{ g_{\ell_{\lambda}}(0)} f(x) dx,
\end{split}
\end{equation}
where the last inequality is based on Lemma \ref{lemma:sumdiff} applied to $a_{\lambda}=g_{\ell_{\lambda}}(0)$ and $b_{\lambda}=g_{\ell_{\lambda}}(x)$. 
Writing $\ell$ for $\ell_{\lambda}$, we have
\begin{eqnarray*}
%\label{eq:bl}
g_{\ell}(0) - g_{\ell}(x) &=& \E \left( e^{-\frac{N_q^2S_{\ell}}{2}} \left( 1 - e^{iN_pN_qB_{\ell}x} \right) \right) \nonumber \\
&\leq&\E\left(e^{-\frac{N_q^2S_{\ell}}{2}}|N_pN_qX_1\ldots X_{\ell-1}|\right)|x|\nonumber\\
&\leq& \E\left(e^{-\frac{N_q^2X_1^2}{2}}|N_qX_1|\right)\E\left(|N_pX_2\ldots X_{\ell-1}|\right)|x|\\
&=& \frac{2}{\pi} \frac{\delta}{\delta+2} r(\delta)^{\ell}\;|x|, \nonumber
\end{eqnarray*}
where the first inequality is due to the fact that $|1-e^{iN_pN_qB_{\ell}x}|\leq |N_pN_qB_{\ell}|\;|x|$, 
the second inequality is due to the fact that $ X_1^2 \leq S_{\ell}$ and the independence of $(N_q,X_1)$ and $(N_p, X_2,\ldots, X_{\ell-1})$, 
and the last equality is obtained using Equations \ref{eq:nqx1} and \ref{eq:npx1xl}. 
Moreover, by using Equation \ref{eq:E7}, we have 
%that $r(\delta)\leq g_{\ell}(0)$. Thus
\begin{eqnarray*}
\frac{g_{\ell}(0) - g_{\ell}(x)}{g_{\ell}(0)} &\leq & \frac{\delta^2}{\pi (\delta + 2)} \left[ \frac{\Gamma( \frac{\delta}{2} )}{\Gamma( \frac{\delta+1}{2} )} \right]^2 r(\delta)^{\ell} \;|x|,
\end{eqnarray*}
and Equation \ref{eq:i2-i1} yields
\begin{equation*} 
\begin{split}
0 \leq \frac{I_2-I_1}{I_2} 
& \leq \frac{\delta^2}{\pi (\delta + 2)} \left( \frac{\Gamma( \frac{\delta}{2} )}{\Gamma( \frac{\delta+1}{2} )} \right)^2 \left( \sum_{\lambda\in \Lambda} r(\delta)^{\ell_{\lambda}}\right) \int_{-\infty}^{\infty} |x|f(x)dx  \\
&\leq \frac{\delta^2}{\pi (\delta + 2)} \left( \frac{\Gamma( \frac{\delta}{2} )}{\Gamma( \frac{\delta+1}{2} )} \right)^2 \left( \sum_{\lambda\in \Lambda} r(\delta)^{\ell_{\lambda}} \right) \frac{\Gamma( \delta^* - 1/2 )  }{ \sqrt{\pi} \Gamma( \delta^* ) } \\
&= \frac{\delta^2}{\pi (\delta + 2)} \left( \frac{\Gamma( \frac{\delta}{2} )}{\Gamma( \frac{\delta+1}{2} )} \right)^2 \left(\sum_{\lambda\in \Lambda} r(\delta)^{\ell_{\lambda}}\right)
r(\delta+d-1)
\end{split}
\end{equation*}
which leads to Equation \ref{ineq:approx} of the manuscript.
%where the last equality is due to \eqref{ouf}.

The following lemmas are used in the proof of Theorem \ref{theorem:disjoint}.
\begin{lemma}
Let $X_1,\ldots,X_{\ell-1}$ be independent identically distributed random variables such that $X_1 \sim Z / \sqrt{Q}$ with $Z \sim N(0,1)$ independent of $Q \sim \chi^2_{\delta+1}$ where $\delta \geq 3$. 
Let $N_p$ and $N_q$ also be standard normal $N(0,1)$ random variables, mutually independent and independent of $X_1,\ldots,X_{\ell-1}$. 
We then have
\begin{equation}
\label{eq:nqx1}
\E \left( e^{-\frac{N_q^2X_1^2}{2}} \left| N_qX_1 \right| \right) = \sqrt{\frac{2}{\pi}}\frac{\delta}{\delta+2} r(\delta)
\end{equation}
and
\begin{equation}
\label{eq:npx1xl}
\E \left( \left| N_pX_1 \ldots X_{\ell - 1} \right| \right) = \sqrt{\frac{2}{\pi}}r(\delta)^{\ell-1},
\end{equation}
where $r(\delta) = \frac{\Gamma(\frac{\delta}{2})}{\sqrt{\pi}\Gamma(\frac{\delta+1}{2})}$. 
Let $S_{\ell}$ as defined in Equation \ref{eq:S_l} then 
\begin{equation}
\label{eq:E7}
 \E \left( e^{-\frac{N_q^2 S_{\ell} }{ 2}} \right) > \frac{2}{\delta} \left( \frac{\Gamma \left( \frac{\delta+1}{2} \right) }{ \Gamma \left( \frac{\delta}{2} \right)} \right)^2.
\end{equation}
\end{lemma}

\begin{proof}
For Equation \ref{eq:nqx1}, the variable $Y = X_1^2 \sim B_2 \left( \frac{1}{2}, \frac{\delta+1}{2} \right)$ which is a Beta distribution of the second kind. 
%The density of such a variable is $f(u)=\frac{1}{B(\alpha, \beta)}\frac{u^{\alpha-1}}{(1+u)^{\alpha+\beta}}{\bf 1}_{(0, +\infty)}(u)$. 
Thus
\begin{eqnarray*}
\E \left( e^{-\frac{N_q^2X_1^2}{2}}|N_qX_1| \right) &=& \E \left( \E \left( e^{-\frac{N_q^2X_1^2}{2}}|N_qX_1| \; \Big| X_1 \right) \right) 
\end{eqnarray*}
where 
\begin{eqnarray*}
\E \left( e^{-\frac{N_q^2X_1^2}{2}}|N_qX_1| \; \Big| X_1 \right) &=& \frac{ |X_1| }{ \sqrt{ 2 \pi } }  \int_{-\infty}^{+\infty} |u| e^{-\frac{(1+X_1^2)}{2} u^2 } du =\sqrt{\frac{2}{\pi}} \frac{ |X_1| }{ 1 + X_1^2 }.
\end{eqnarray*}
Thus
\begin{eqnarray*}
\E \left( e^{-\frac{N_q^2X_1^2}{2}}|N_qX_1| \right) = \sqrt{ \frac{2}{\pi} } \E \left( \frac{ |X_1| }{ 1 + X_1^2 } \right) = \sqrt{ \frac{2}{\pi} } \E \left( \frac{ \sqrt{Y} }{ 1 + Y } \right) = \sqrt{ \frac{2}{\pi} } \frac{ \delta }{ \delta+2} r(\delta). 
\end{eqnarray*}

For Equation \ref{eq:npx1xl}, by using the mutual independence of $N_p$ and $X_1,\ldots, X_{\ell-1}$ and the fact that $X_i^2 \sim B_2 \left( \frac{1}{2}, \frac{\delta+1}{2} \right)$ for $ i = \{ 1, \ldots, \ell-1 \}$, we have
\begin{eqnarray*}
\E \left( |N_pX_1\ldots X_{\ell-1}| \right) = \E \left( |N_p| \right) \E \left( |X_1| \right)^{\ell-1} = \frac{\sqrt{2}}{\sqrt{\pi}}\left(\frac{\Gamma(\frac{\delta}{2})}{\sqrt{\pi} \Gamma(\frac{\delta+1}{2})}\right)^{\ell-1}=\sqrt{\frac{2}{\pi}} r(\delta)^{\ell-1}
\end{eqnarray*}

For Equation \ref{eq:E7}, from \citet[Example 9]{chamayou1991explicit}, if $U \sim B_2( a, b )$ for $b>a$ and $V$ are independent. Then $U(1+V) \sim V$ if and only if $V \sim B_2( a, b-a )$. 
One applies this to $U = X_1^2$, $V = S$, $a = \frac{1}{2}$, and $b = \frac{\delta+1}{2}$ since $S' = \sum_{i=2}^{\infty}X_2^2\ldots X_i^2\sim S$ and $X_1^2(1+S')=S$. 
Thus $S \sim B_2 \left( \frac{1}{2}, \frac{\delta}{2} \right)$. 
Since $S_{\ell} < S$, we have

\begin{equation*}
% g_{\ell}(0) = 
 \E \left( e^{-\frac{N_q^2 S_{\ell} }{ 2}} \right) > \E \left( e^{-\frac{N_q^2 S }{ 2}} \right) = \frac{2}{\delta} \left( \frac{\Gamma \left( \frac{\delta+1}{2} \right)}{\Gamma \left( \frac{\delta}{2} \right)} \right)^2.
\end{equation*}
\end{proof}

\begin{lemma}
\label{lemma:sumdiff}
Let $a_1, \ldots, a_n$ and $b_1,\ldots,b_n$ be complex numbers such that $\left| b_i \right| \leq \left|a_i \right|$, for $i = 1, \ldots, n$. Then
\begin{equation*}
%\label{eq:sumdiff}
\left| \prod_{i=1}^na_i - \prod_{i=1}^nb_i \right| \leq \prod_{i=1}^n \left| a_i \right| \sum_{j=1}^n\frac{|a_j-b_j|}{|a_j|}.
\end{equation*}
\end{lemma}

\begin{proof}
\begin{eqnarray*}
\left| a_1 \ldots a_n - b_1 \ldots b_n \right| &=& | \left( a_1 \ldots a_n - b_1 a_2 \ldots a_n \right) + \left( b_1 a_2 \ldots a_n - b_1 b_2 a_3 \ldots a_n \right)   \\
                                                                         && \hspace{3cm} + \ldots \ldots + \left( b_1 \ldots b_{n-1} a_n - b_1 b_2 \ldots b_n \right) | \\
&\leq& \sum_{j=1}^n \left| \prod_{i=1}^{j-1} b_i \prod_{i=j}^n a_i - \prod_{i=1}^{j} b_i \prod_{i=j+1}^n a_i \right| =
\sum_{j=1}^n \left| \prod_{i=1}^{j-1} b_i \prod_{i=j+1}^n a_i \right| \left| a_j - b_j \right| \\
&\leq& \sum_{j=1}^n \left| \prod_{i=1}^{j-1} a_i \prod_{i=j+1}^n a_i \right| \left| a_j-b_j \right| = \prod_{i=1}^{n} \left| a_i \right| \sum_{j=1}^n \left| 1 - \frac{b_j}{a_j} \right|
\end{eqnarray*}
\end{proof}

% - - - - - - - - - - - - - - - - - - - - - - - - - - - - - - - - - - - - - - - - - - - - - - - - - - - - - - - - - - - -|
%\section{Appendix }
\section{Proof of Theorem \ref{theorem:joint} }
\label{sup:proof-theorem2}

The proof of Theorem \ref{theorem:joint} is long and is done in the next three subsections. 
%$ In Subsection \ref{subsec: newintegral}, we give a new integral expression \eqref{eq:newi1i2} for $I_1/I_2$, in terms of $b_1$. 
In Subsection \ref{subsec:bilinear}, we show in Proposition \ref{proposition:quadratic} that, conditional on a quantity $\Psi_{\cup}^-$ as defined in Equation \ref{eq:smallD}  of the manuscript, $b_1$ can be expressed as a bilinear form. 
In Proposition \ref{proposition:mixture} of Subsection \ref{subsec:mixture}, using its expression as a bilinear form, we show that $b_1$ is distributed like the continuous scale mixture of centered Gaussian variables. 
This allows us to deduce that there exists a unique $v_D$ such that the normal $N(0, v_D)$ distribution best approximates the distribution of $b_1$. 
Finally, in Subsection \ref{subsec:intunderapprox}, we prove Theorem \ref{theorem:joint}:
under the assumption that $v_D$ is small, $I_1/I_2$ can accurately be approximated by 1 or equivalently the ratio of the normalizing constants can accurately be approximated by Equation \ref{eq:approx} of the manuscript, which is what we want to prove. 

% - - - - - - - - - - - - - - - - - - - - - - - - - - - - - - - - - - - - - - - - - - - - - - - - - - - - - - - - - - - -|
\subsection{Expression of $b_1$ as a bilinear form}
\label{subsec:bilinear}

Here we want to show that the distribution of $b_1$ is a scale mixture of normal distributions which can be approximated by another $N(0, v_D)$ distribution where the variance $v_D$ depends on $\Psi_{\cup}^-$. 
But, to do so, we must first express $b_1$ as a bilinear form in two standard normal random vectors.

Regarding $E_q = \{ (i,j): (i,q) \in E \}$ and $E_p = \{ (i,j): (i,p) \in E \}$, we define 
\begin{equation*}
%\label{eq:v+-}
\Psi_{E_q^-} = \{ \psi_{ij}: (i,j) \in E_q \setminus E_p \}, \;\;\; \Psi_{E_p^-} = \{ \psi_{ij}: (i,j) \in  E_p \setminus E_q \}.
\end{equation*}
The $\Psi_{E_q^-}$ represents the free elements of matrix $\Psi$ regarding to just the neighbor of $q$ and the same for $p$.
In the following, we express $D$ and $b_1$ as polynomials in $\Psi_{E_q^-}$ and $\Psi_{E_p^-}$. 

\begin{proposition}
\label{proposition:quadratic}
Let $N_q^{-}$ denote the set of nodes that are neighbours of $q$ but not of $p$ and $N_p^{-}$ denote the set of nodes that are neighbours of $p$ but not of $q$. 
There exist  vectors $M_i^q \in \R^{N_q^{-}}$ and $M_i^p \in \R^{N_p^{-}}$, $i = \{ 1, \ldots, q-d-1\}$, functions of $\Psi_{\cup}^-$, such that 
%% if $R^q$ and $R^p$ are the two symmetric matrices defined by
%% \[ R^q = \sum_{i=1}^{p-d-2} M^q_i (M^q_i)^t, \ \ R^p = \sum_{i=1}^{p-d-2} M^p_i (M^p_i)^t \] and 
if $C$ is the $| N_q^{-} | \times | N_p^{-} |$-dimensional matrix
\[
C = \sum_{i=1}^{q-d-1} M^q_i (M^p_i)^t
\]
then we have 
\begin{eqnarray}
%\label{LEBOND}
%D &=& \sum_{(i,j)\in \overline{E}, \ j<q} \psi_{ij}^2 + \tr \left( \Psi_{N_q^-} R^q \Psi_{N_q^-} \right) + \tr \left( \Psi_{N_p^-} R^p \Psi_{N_p^-} \right),  \\
\label{eq:LEBONb}
b_1 = \tr \left(  \Psi_{E_q^-} C  \Psi_{E_p^-} \right).
\end{eqnarray} 
%where $ \Psi_{E_q^-}$ and $ \Psi_{E_p^-}$ are defined in \eqref{eq:v+-}. 
Furthermore, $\left\{ M^q_i, M^p_i \right\}_{i=1}^{q-d-1}$, $ \Psi_{E_q^-}$ and $ \Psi_{E_p^-}$ are  independent.
\end{proposition}
%The proof is given in Section \ref{sec:lemma:psie} of the Supplementary File.

\begin{proof}
%% From the expression of $D$ in \eqref{D}  of the manuscript we have
%% \begin{eqnarray*}
%% D = \sum_{(i,j)\in \overline{E}} \psi_{ij}^2 =  \sum_{(i,j)\in \overline{E}, \ j<q} \psi_{ij}^2  + \sum_{(i,q)\in \overline{E}, i\le p-d-1}\psi_{iq}^2+\sum_{(i,p)\in \overline{E}, i\le p-d-1}\psi_{ip}^2.
%% \end{eqnarray*}
From the expression of $b_1$ in Equation \ref{eq:b}  of the manuscript we have
\begin{eqnarray*}
b_1 = \sum^{q - d - 1 }_{i=1} \psi_{iq}\psi_{ip},
\end{eqnarray*}
which is based on assume the nodes which are neighbours to both $q$ and $p$, are numbered $q-d, q-d+1, \ldots, p$. 
By Equation \ref{eq:aliye} of the manuscript, each $\psi_{iq}, (i,q)\in \overline{E}$ is equal to the sum of products 
\begin{equation}
\label{eq:sum}
\psi_{iq} = \frac{-1}{\psi_{ii}}\sum_{l=1}^{i-1} \psi_{l i}\psi_{l q},
\end{equation}
where each of these $\psi_{l i}$ or $\psi_{l q}$, $l = \{1, \ldots, q-d-1 \}$ may be free or not free. 
If $\psi_{l q}$ is free, $l$ necessarily belongs to $N_q^{-}$ because it is a neighbour of $q$ and, since $i \le q-d$ and $l \le i$, it cannot be a neighbour of both $q$ and $p$. 
If it is not free, then, we write the expression of $\psi_{l q}$ according to Equation \ref{eq:aliye} of the manuscript and we repeat this process until $\psi_{l q}$ has been expressed in terms of a ratio of a product of $\psi_{uv}, (u,v)\in E, u\leq l, v\leq q$, 
one of which  is necessarily (since the sum in Equation \ref{eq:sum} is finite) equal to $\psi_{u_{l} q}$ for some $u_{l}\leq l, u_{l} \in N_q^{-}$, and a product of  $\psi_{vv}, v\leq l$. 
Similarly $\psi_{l i}$ is free or not free. If not free, it will be expressed as a ratio of products of elements of $\psi_E$, none of which, in the numerator, can be equal to $\psi_{u_{l}q}$ since it is the product of entries $\psi_{uv}$ of $\psi$ with $u\leq l, v\leq i<q$. 
Thus, from Equation \ref{eq:sum}, we can write, for each $i = \{ 1, \ldots, q-d-1 \}$ 
\begin{equation}
\label{eq:psiiq}
\psi_{iq} = \frac{-1}{\psi_{ii}} \sum_{l=1}^{i-1} \psi_{l i}\psi_{l q} = \sum_{l \in N_q^{-}}(M_i^q)_{l} \; \psi_{u_{l},q} = \tr \left( M_i^q \Psi_{E_q^-}  \right),
\end{equation}
where $(M_i^q)_{l}, l \in N_q^{-}$ are the components of $M_i^q$, some of which can be equal to 0, if the Cholesky equations (in Equation \ref{eq:aliye} of the manuscript) do not lead to that particular $l \in N_q^{-}$, or 1, if $l \in N_q^{-}$. 
%% Now, we have
%% \begin{eqnarray*}
%% \sum_{(i,q) \in \overline{E}, i\le p-d-1} \psi_{iq}^2 
%%   &=& \sum_{(i,q) \in \overline{E}, i\le p-d-1}  \tr \left( M_i^q \Psi_{N_q^-}  \right) \\
%%   &=& \Psi_{N_q^-}^t \left( \sum_{i=1}^{p-d-2} M_i^q (M_i^q)^t \right) \Psi_{N_q^-} \\
%%   &=& \Psi_{N_q^-}^t R^q \Psi_{N_q^-}. 
%% \end{eqnarray*}
%A similar argument holds for $\psi_{ip}$ and thus \eqref{LEBOND} of the manuscript is proved.
Similarly, we have
\begin{equation}
\label{eq:psiip}
\psi_{ip} = \frac{-1}{\psi_{ii}} \sum_{l=1}^{i-1} \psi_{l i} \psi_{l p} = \sum_{l \in N_p^{-} } (M_i^p)_{l} \; \psi_{v_{l},p} = \tr \left( M_i^p \Psi_{E_p^-}  \right),
\end{equation}
for some $v_{l} < l, v_{l} \in N_p^{-}$ and Equation \ref{eq:LEBONb} follows from Equations \ref{eq:psiiq} and \ref{eq:psiip}.
\end{proof}

\begin{example} 
\label{ex:smallgraph}
Consider the graph in Figure \ref{fig:paths} (right)  of the manuscript where $q=6$,
$p=7$, $N_q = \{ 1, 2, 5 \}$, $ N_p = \{ 3, 4, 5 \}$, $N_q^{-} = \{ 1, 2 \}$, $ N_p^{-} = \{ 3, 4 \}$, $d=1$, and 
$\psi_E = \{ \psi_{14}, \psi_{16}, \psi_{23},\psi_{24},\psi_{26}$,  $\psi_{37}, \psi_{47}, \psi_{56},\psi_{57} \}$. 
Thus
%\[\Psi_{E_q^-} =  \left( \psi_{16}, \psi_{26} \right) \;\; \text{ and } \;\; \Psi_{E_p^-} =  \left( \psi_{37}, \psi_{47} \right).\]
$\Psi_{E_q^-} =  \left( \psi_{16}, \psi_{26} \right)$ and $\Psi_{E_p^-} =  \left( \psi_{37}, \psi_{47} \right)$.
Using the notation $X_{ij} = \psi_{ij} / \psi_{jj}$ for convenience, the non-free entries are
\begin{eqnarray*}
\psi_{34} &=& - \psi_{24} X_{23},   \\
\psi_{36} &=& - \psi_{26} X_{23} = - \left( \psi_{16}, \psi_{26} \right) \left( 0, X_{23} \right)^t = \tr( \Psi_{E_q^-} M^q_3 ), \\
M^q_3 &=& \left( 0, X_{23} \right),   \\
\psi_{46} &=& - \psi_{26} X_{23}^2 X_{24} - \psi_{26} X_{24} - \psi_{16} X_{14} 
                     = -\left( \psi_{16}, \psi_{26} \right) \left( X_{14}, X_{24} + X_{23}^2 X_{24} \right)^t 
                     = \tr \left( \Psi_{E_q^-} M^q_4 \right),   \\
M^q_4 &=& - \left( X_{14}, X_{24} + X_{23}^2 X_{24} \right)^t,   \\
\psi_{67} &=& \frac{1}{ \psi_{66} } \left( - \psi_{57} \psi_{56} 
                        + \psi_{26} X_{23}^2X_{24}\psi_{47} 
                        + \psi_{26} X_{24} \psi_{47} 
                        + \psi_{16} X_{14} \psi_{47} 
                        + \psi_{26} X_{23}\psi_{37} \right) \\
                      &=& \frac{1}{ \psi_{66} } \left( A_1 + b_1 \right),
\end{eqnarray*}
where
\begin{eqnarray*}
A_1 &=&  - \psi_{57} \psi_{56}, \\ 
b_1 &=&  \psi_{26} X_{23}^2X_{24} \psi_{47} 
                        + \psi_{26} X_{24} \psi_{47} 
                        + \psi_{16} X_{14} \psi_{47} 
                        + \psi_{26} X_{23} \psi_{37}.
 %            =  \Psi_{E_q^-}^t C \Psi_{E_p^-},
\end{eqnarray*}
It leads $b_1 = \Psi_{E_q^-}^t C \Psi_{E_p^-}$ where 
\[
C = \left[ 
\begin{array}{cc}
        0          &X_{14} \\ 
        X_{23} & X_{24} + X_{23}^2 X_{24} 
\end{array} \right].
\]
It also follows from the definition of $\Psi_{E_q^-}$ that $M_1^q = (1,0)^t$ and $M_2^q = (0,1)^t$. 
From the definition of $\Psi_{E_p^-}$, we have $M_1^p = M_2^p = (0,0)^t$ and $M_3^p = (1,0)^t$ and $M_4^p = (0,1)^t$. 
We can then verify
\[
\tr \left( \Psi_{E_q^-}^t C \Psi_{E_p^-} \right) = \sum_{i=1}^4 \tr \left( \Psi_{E_q^-} M_i^q \right) \tr \left( \Psi_{E_p^-} M_i^p \right).
\]
\end{example}

% - - - - - - - - - - - - - - - - - - - - - - - - - - - - - - - - - - - - - - - - - - - - - - - - - - - - - - - - - - - -|
\subsection{A normal approximation to the distribution of $b_1$}
\label{subsec:mixture}

%We saw from the expression of $D$ and $b$ in \eqref{LEBOND} and \eqref{LEBONb} respectively, that if we condition on $\psi_{\tilde{D}_1}$, as defined in \eqref{eq:smallD}, then $b_1$ can be expressed as the bilinear form in $\Psi_{E_q^-}$ and $\Psi_{E_p^-}$,
We see in Equation \ref{eq:LEBONb}, that if we condition on $\Psi_{\cup}^-$, as defined in Equation \ref{eq:smallD} of the manuscript, then $b_1$ can be expressed as the bilinear form of $\Psi_{E_q^-}$ and $\Psi_{E_p^-}$ as 
\[
b_1 = \Psi_{E_q^-}^t C \Psi_{E_p^-},
\]
where $C = \sum_{i=1}^{q-d-1} (M^q_i)^t M_i^p$ is a matrix of rank $m \leq \min \left( | N_q^{-} |, | N_p^{-} | \right)$. 
Once $\Psi_{\cup}^-$ is known, $C$ is fixed.
We are now going to show that, conditional on $\Psi_{\cup}^-$,  the distribution of $b_1$ has the following property.

\begin{proposition} 
\label{proposition:mixture}
When conditioned by $\Psi_{\cup}^-$,   the distribution of $b_1$ is a continuous scale mixture of centered normal distributions. 
More precisely, $b_1$  follows the same distribution as $X\sqrt{Y/2}$ where $X$ and $Y$ are independent with $X\sim N(0,1)$ and  
\[
Y = \sum_{i=1}^m \frac{Y_i}{\lambda_i} \;\; \text{ where } \;\; Y_i  \stackrel{iid}{\sim} \chi_2^2
\]
and where $1/\lambda_1,\ldots,1/\lambda_m$ are the non-zero eigenvalues of $C^TC.$
\end{proposition}

\begin{proof} 
We first show that $b_1$ follows the same distribution as $X\sqrt{Y/2}$. To do so, it suffices to show the two Laplace transforms $ \E \left( e^{s b_1} \right)$ and $\E \left( e^{sX\sqrt{Y/2}} \right)$ coincide. 
Integrating this last expected value first with respect to $X$, holding $Y$ fixed, and then with respect to $Y$, we obtain
\begin{equation*}
%\label{lp}
\E \left( e^{sX\sqrt{Y/2}} \right) = \E \left( e^{s^2Y/2} \right) = \prod_{j=1}^m \left( 1 - \frac{\lambda_j s^2}{2} \right)^{\frac{-1}{2}}.
\end{equation*}
Next, from Equation \ref{eq:LEBONb} and then integrating with respect to $\Psi_{E_q^-}$, we have 
\[
\E \left( e^{s b_1} \right) = \E \left( e^{ s \tr \left( \Psi_{E_q^-} C \Psi_{E_p^-} \right) } \right) = \E \left( e^{ \frac{s^2}{2} \tr \left( \Psi_{E_p^-} C^t C \Psi_{E_p^-} \right) } \right).
\]
Now we have 
\[
\tr \left( \Psi_{E_p^{-}} C^t C \Psi_{E_p^{-}} \right) \sim Z^t \mathrm{diag}(1/\lambda_1,\ldots,1/\lambda_m) Z,
\]
where $Z = (Z_1,\ldots,Z_m)$ are independent  $N(0,1)$ random variables. Thus  $b_1$ and $X\sqrt{Y/2}$ have the same Laplace transform.
\end{proof}

To show the distribution of $b_1$ is a scale mixture of centered normals, we note that if $X \sim N(0,1)$ and $V = Y/2$ is any positive random variable with distribution  $\mu(dv)$, then if $U=X\sqrt{Y/2}$,  the density of $U$ is
\begin{eqnarray*}
f_U(u) = \int_0^{+\infty} \frac{e^{-\frac{u^2}{2v}} }{ \sqrt{2 \pi v}} \mu(dv).
\end{eqnarray*}
So, the distribution of $U$, that is the distribution of $b_1$, is a mixture of normal $N(0,v)$ distributions. 
Following \citet[Theorem 3.1.]{letac2020gaussian}, with the distribution of $b_1$ for $f$, we deduce that there exists a unique $v_D$ such that the normal $N(0, v_D)$ distribution best approximates the distribution of $b_1$.

% - - - - - - - - - - - - - - - - - - - - - - - - - - - - - - - - - - - - - - - - - - - - - - - - - - - - - - - - - - - -|
\subsection{Expression $I_1/I_2$ regarding $b_1 \sim N(0, v_D)$}
\label{subsec:intunderapprox}

We will now derive an expression for $I_1/I_2$ when we approximate the distribution of $b_1$ by the $N(0, v_D)$ distribution.
We start with the following lemma.

\begin{lemma}
\label{lemma:main}
Under the approximation $b_1\sim N(0, v_D)$, we have
\begin{equation*}
%\label{eq:bige}
\E \left( h(b_1, \delta^*) \big| \Psi_{\cup}^- \right) = \frac{ v_D^{\delta^* } }{ 2^{\delta^*} \Gamma( \delta^* )} \int_0^{\infty} t^{\delta^* - \frac{1}{2} } (1+t)^{\frac{-1}{2}} e^{-\frac{v_D t}{2}} dt,
\end{equation*}
where $\delta^* = \frac{\delta+d}{2}$ and $h(b_1, \delta^*)$ defined in Equation \ref{eq:int-b1} of the manuscript. 
Moreover, when $v_D$ is small, we have
\begin{equation}
\label{eq:ufunction}
\E \left( h(b_1, \delta^*) \big| \Psi_{\cup}^- \right) = 1 - \frac{ \Gamma( \delta^* + \frac{1}{2} ) }{ \Gamma( \delta^* ) } \left( \frac{v_D}{2} \right)^{\delta^*} {\mathcal O} \left( \Big| \frac{v_D}{2} \Big|^{\delta^*-1} \right).
\end{equation}
%where$$f( \delta^*, v_D ) = \frac{ \Gamma( \delta^* + \frac{1}{2} ) }{ \Gamma( \delta^* ) } \left( \frac{v_D}{2} \right)^{\delta^*}. $$
\end{lemma}

\begin{proof}
For $b_1\sim N(0, v_D)$, we have  
\begin{eqnarray*}
\E \left( h(b_1, \delta^*) \big| \Psi_{\cup}^- \right)
 &=& \frac{ 2^{-\delta^*} }{ \Gamma( \delta^* )} \int_{-\infty}^{+\infty} \left( \int_0^{+\infty} y^{\delta^* - 1} e^{ \frac{-1}{2} \left( y + \frac{b_1^2}{y} \right) }  dy \right) \frac{e^{-\frac{b_1^2}{2v_D}}}{\sqrt{2\pi v_D}} d b_1 \\
 &=& \frac{ 2^{-\delta^*} }{ \Gamma( \delta^* )} \int_0^{+\infty} \left( \int_{-\infty}^{+\infty} \frac{e^{-\frac{y + v_D }{ 2v_D y} b_1^2} }{ \sqrt{2 \pi v_D}} db_1 \right) y^{\delta^* - 1} e^{-\frac{y}{2}} dy \\
&=& \frac{ 2^{-\delta^*} }{ \Gamma( \delta^* )} \int_0^{+\infty} \left( \frac{y}{y + v_D} \right)^\frac{1}{2} e^{-\frac{y}{2}} y^{ \delta^* - 1} dy \\
 &=& \frac{ 2^{-\delta^*} v_D^{\delta^* } }{ \Gamma( \delta^* )} \int_0^{\infty} t^{\delta^* - \frac{1}{2} } (1+t)^{\frac{-1}{2}} e^{-\frac{v_D t}{2}} dt. 
\end{eqnarray*}
We note that the above integral is a confluent hypergeometric function of the form 
\[
\Gamma(a)U(a,b,z)=\int_0^{+\infty}e^{-zt}t^{a-1}(1+t)^{b-a-1}dt,
\]
%Abramovitz&Stegun , p.505, formula 13.2.5
with $z = \frac{v_D}{2}$, $a = \delta^* + \frac{1}{2}$, and $b = \delta^* + 1$ (see \cite{abramovitzStegun}, p.505, formula 13.2.5) and from p.508, formula 13.5.6 of the same, we know when $|z|\to 0$ and $b>2$, then 
\[
U(a,b,z)=\frac{\Gamma(b-1)}{\Gamma(a)}z^{1-b}+{\mathcal O}(|z|^{b-2}).
\]
This yields Equation \ref{eq:ufunction}.
%, and from \eqref{eq:newi1i2}, we obtain \eqref{eq:approxvsmall} and \eqref{approx} holds.
%Finally, the fact that $I_1/I_2 \leqslant 1$ follows immediately from \eqref{eq:relationship2} since $b_1^2 / Y \geqslant 0$ is always positive and $e^{-b_1^2/Y} \leqslant 1.$
\end{proof}

%As already mentioned at the end of Section \ref{subsec:reformation}, we will give a different expression for $I_1/I_2$,  use an approximation to the distribution of $b_1$ and then find a sufficient condition for $I_1/I_2$ under that approximation to be close to 1.
%Of course, this does not mean that $I_1/I_2$ is not well approximated by 1 when the sufficient condition is not satisfied. It simply means that this particular approximation is not always good. 

% - - - - - - - - - - - - - - - - - - - - - - - - - - - - - - - - - - - - - - - - - - - - - - - - - - - - - - - - - - - -|
\section{Pseudo-code for $I_1/I_2$ in Theorem \ref{theorem:joint}}
%\eqref{eq:newi1i2}
\label{sec:pseudo-code}

 Since we cannot evaluate the approximation in Theorem \ref{theorem:joint} directly, we write instead
%To evaluate approximation \eqref{eq:ratioJoint} in Theorem \ref{theorem:joint}, we write
\begin{eqnarray*}
%\frac{  2^{\delta^*} \Gamma( \delta^* ) }{ \Gamma( \delta^* + \frac{1}{2} ) } \times
\frac{I_1}{I_2} = \frac{{\mathbb E}\left(e^{-\frac{D}{2}} g\left( \delta^*, v_D \right) \right)}{{\mathbb E}\left(e^{-\frac{D}{2}}\right)}=\frac{\int g\left( \delta^*, v_D \right) e^{-\frac{D}{2}}\pi(D) dD}{\int e^{-\frac{D}{2}} \pi(D)dD}=\int  g\left( \delta^*, v_D \right)\pi_1(D) dD
\end{eqnarray*}
where $\pi(D)$ is the unknown density of $D$, and %$\pi_1(D)=e^{-\frac{D}{2}}\pi(D)(\int e^{-\frac{D}{2}} \pi(D)dD)^{-1}$. % and $I_3(D)$ is equal to the left-hand side of \eqref{eq:jointEQ2}. 
\[
\pi_1(D) = \frac{e^{-\frac{D}{2}}\pi(D) }{ \int e^{-\frac{D}{2}} \pi(D)dD }.
\] 
% and $I_3(D)$ is equal to the left-hand side of \eqref{eq:jointEQ2}. 
We then approximate $I_1/I_2$  %$\int g\left( \delta^*, v_D \right) \pi_1(D) dD$ 
%following  the pseudo-code given in Section \ref{sec:pseudo-code} of the Supplementary file. 
%Let $N$ be the number of iterations such as $10^3$ or $10^6$ for example. 
%The value of \eqref{eq:ratioJoint} of the manuscript is obtained through the following sequence of steps.
by the following sequence of steps.
\begin{enumerate}
\item[\bf{1.}] Generate $D_i, i=1,\ldots,N$ the usual way. Divide the range of $D$ into appropriate small intervals $int^{(q)}, q=1,\ldots, Q$, and for each interval $int^{(q)}$, compute the relative frequency $f^{(q)}$.
\item[\bf{2.}]Compute 
 $D^{(q)} = \frac{1}{Nf^{(q)}} \sum_{D_i \in int^{(q)}}D_i$ and
 $r^{(q)}=\frac{e^{-\frac{D^{(q)}}{2}}f^{(q)}}{\sum_{j=1}^Q e^{-\frac{D^{(j)}}{2}}f^{(j)}}$, $q=1,\ldots,Q.$
\item[\bf{3.}] Sample $M$ values of $D^{(m)}, m=1,\ldots,M$ with probabilities given by the empirical distribution of the $r^{(q)}, q=1,\ldots,Q$.
\item[\bf{4.}] For each $D^{(m)}$,
generate $b_1^{(m,k)}, k=1,\ldots,K$ the usual way. Compute $v_{D^{(m)}} =\frac{1}{K} \sum_{k=1}^K(b_1^{(m,k)}-\overline{b_1}^{(m)})^2 $ where 
$\overline{b_1}^{(m)} =\frac{1} {K}b_1^{(lmk)}$.
\item[\bf{5.}] Compute 
\[
I_3( v_{D^{(m)}} ) = {\mathbb E} \left(  \sqrt{ \frac{t}{t+1} } \right)
\]
by simulating from the $\Gamma( \delta^*, \frac{v_{D^{(m)}} }{2} )$ distribution for $t$.
\item[\bf{6.}] Take the average $\frac{1}{M} \sum_{m=1}^M I_3( v_{D^{(m)}} )$ as the estimate of $I_1/I_2$.
\end{enumerate}

% - - - - - - - - - - - - - - - - - - - - - - - - - - - - - - - - - - - - - - - - - - - - - - - - - - - - - - - - - - - -|
\section{Additional simulation results}
\label{sec:plots}

\begin{figure} %[!ht]
\begin{center}
    \includegraphics[width=18cm]{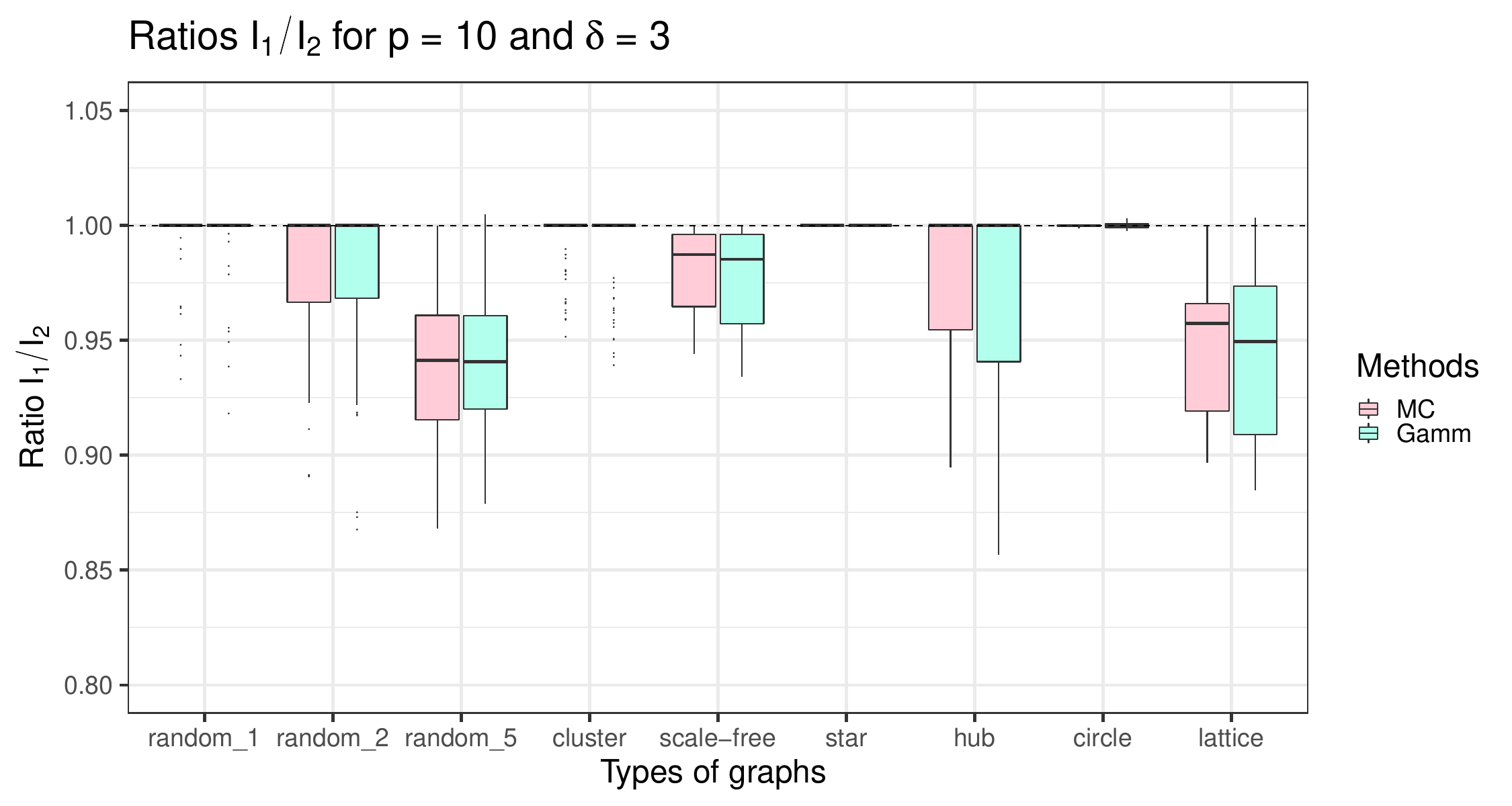}
    \includegraphics[width=15cm]{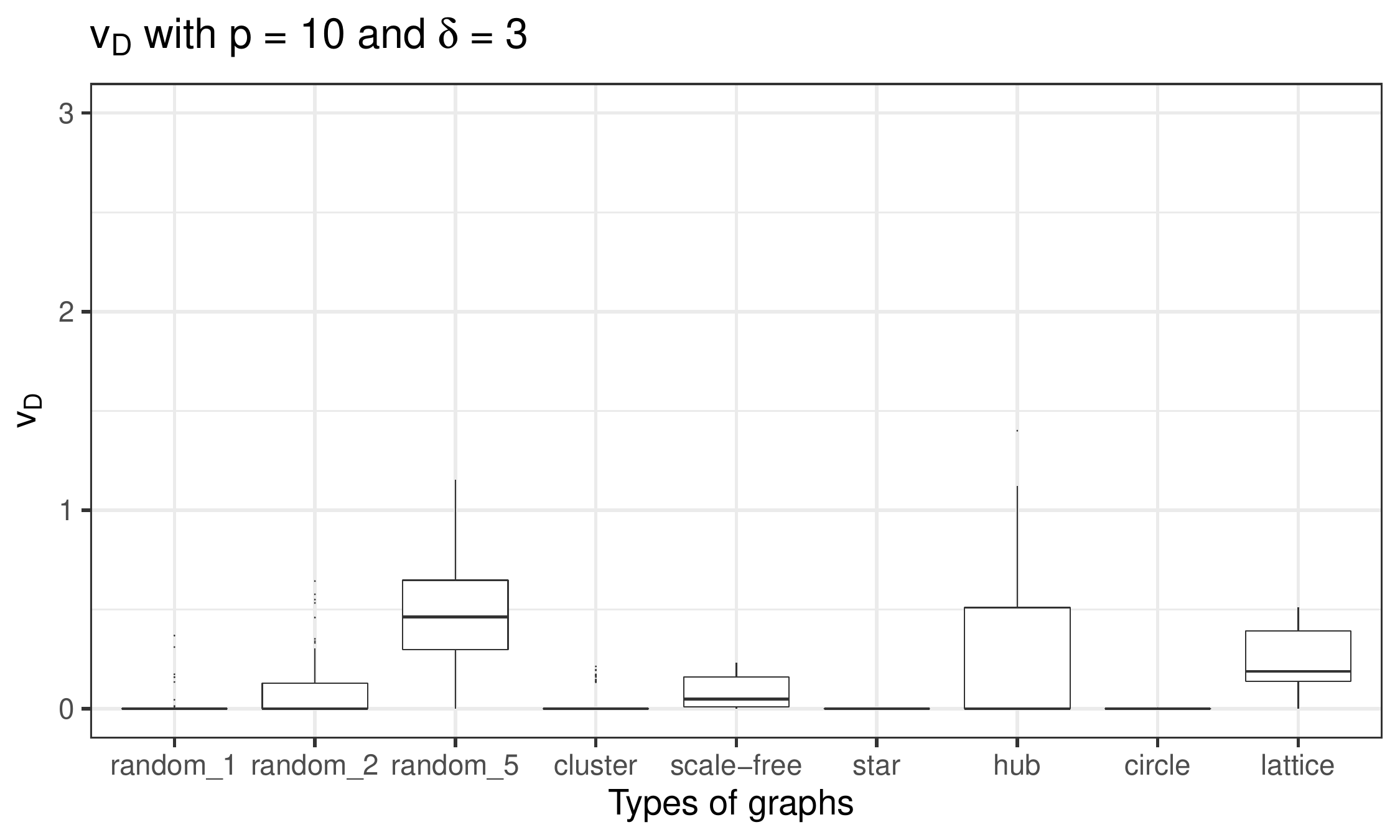}
\end{center}
\caption{ 
(Top) The boxplot for the ratio $I_1/I_2$ computed by the MC approach of \cite{atay2005monte} (in red) and approximation \eqref{eq:ratioJoint} (in blue). (Bottom) The boxplot of the variance $v_D$ of $b_1$ for the corresponding graphs. These computations are done over 100 replications for nine different graphs (Figure \ref{fig:plot_graphs_p30}) with $p=10$ nodes and $\delta=3$. \label{fig:boxplot_p10} }
\end{figure}

\begin{figure} %[!ht]
\begin{center}
    \includegraphics[width=18cm]{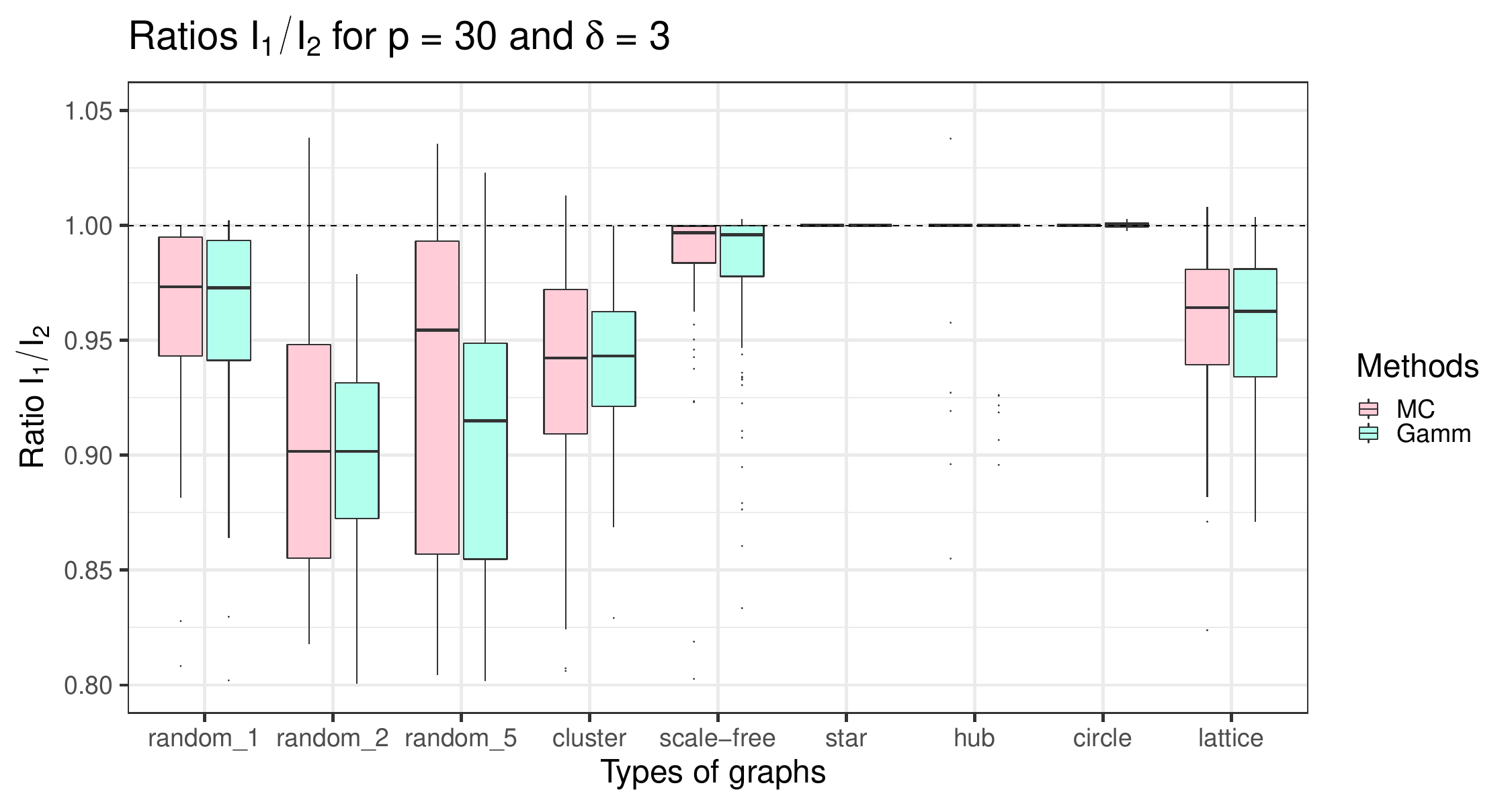}
    \includegraphics[width=15cm]{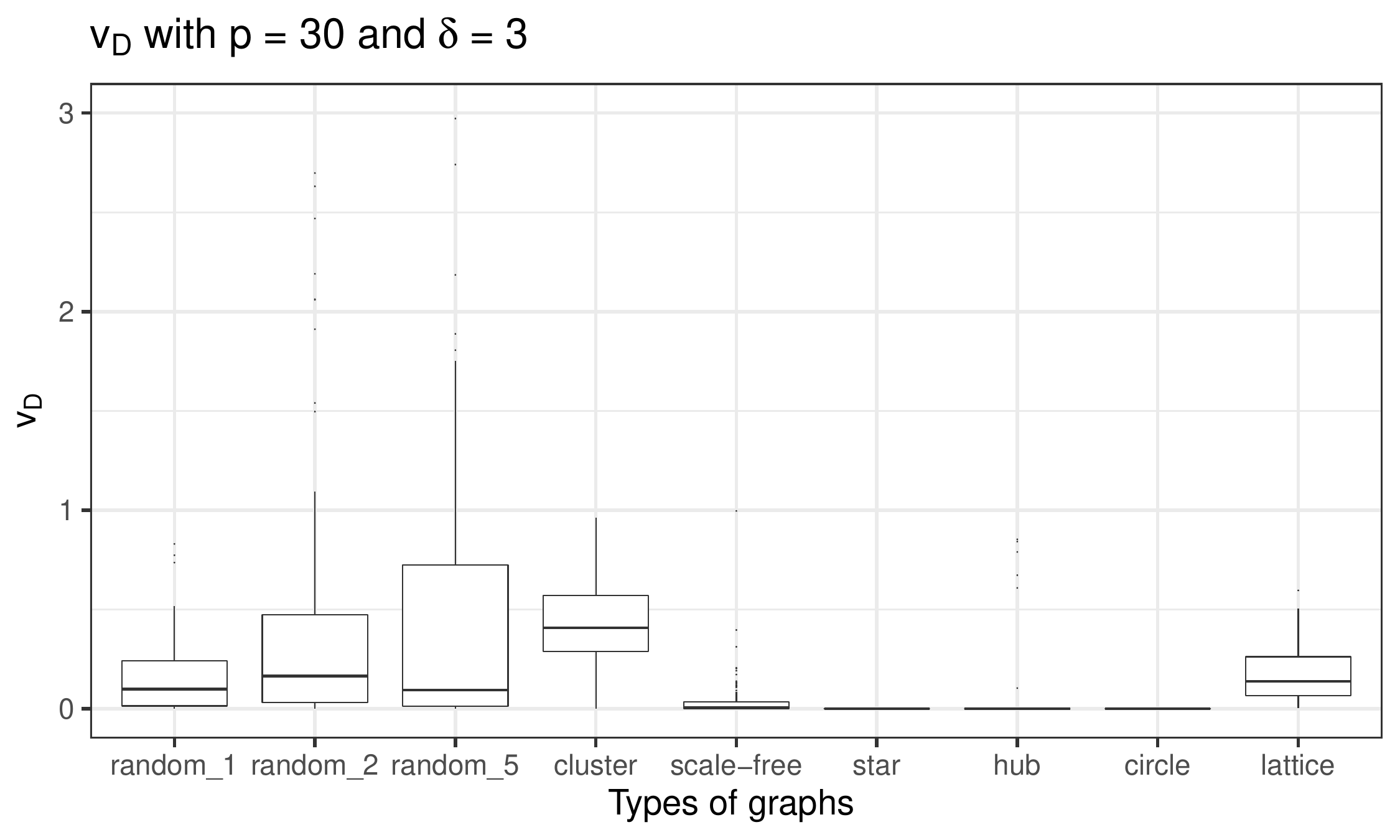}
\end{center}
\caption{ (Top) The boxplot for the ratio $I_1/I_2$ computed by the MC approach of \cite{atay2005monte} (in red) and our approximation \eqref{eq:ratioJoint} (in blue). (Bottom) The boxplot of the variance $v_D$ of $b_1$ for the corresponding graphs. These computations are done over 100 replications for nine different graphs (Figure \ref{fig:plot_graphs_p30}) with $p=30$ nodes and $\delta=3$. \label{fig:boxplot_p30} }
\end{figure}

\begin{figure} %[!ht]
\begin{center}
    \includegraphics[width=18cm]{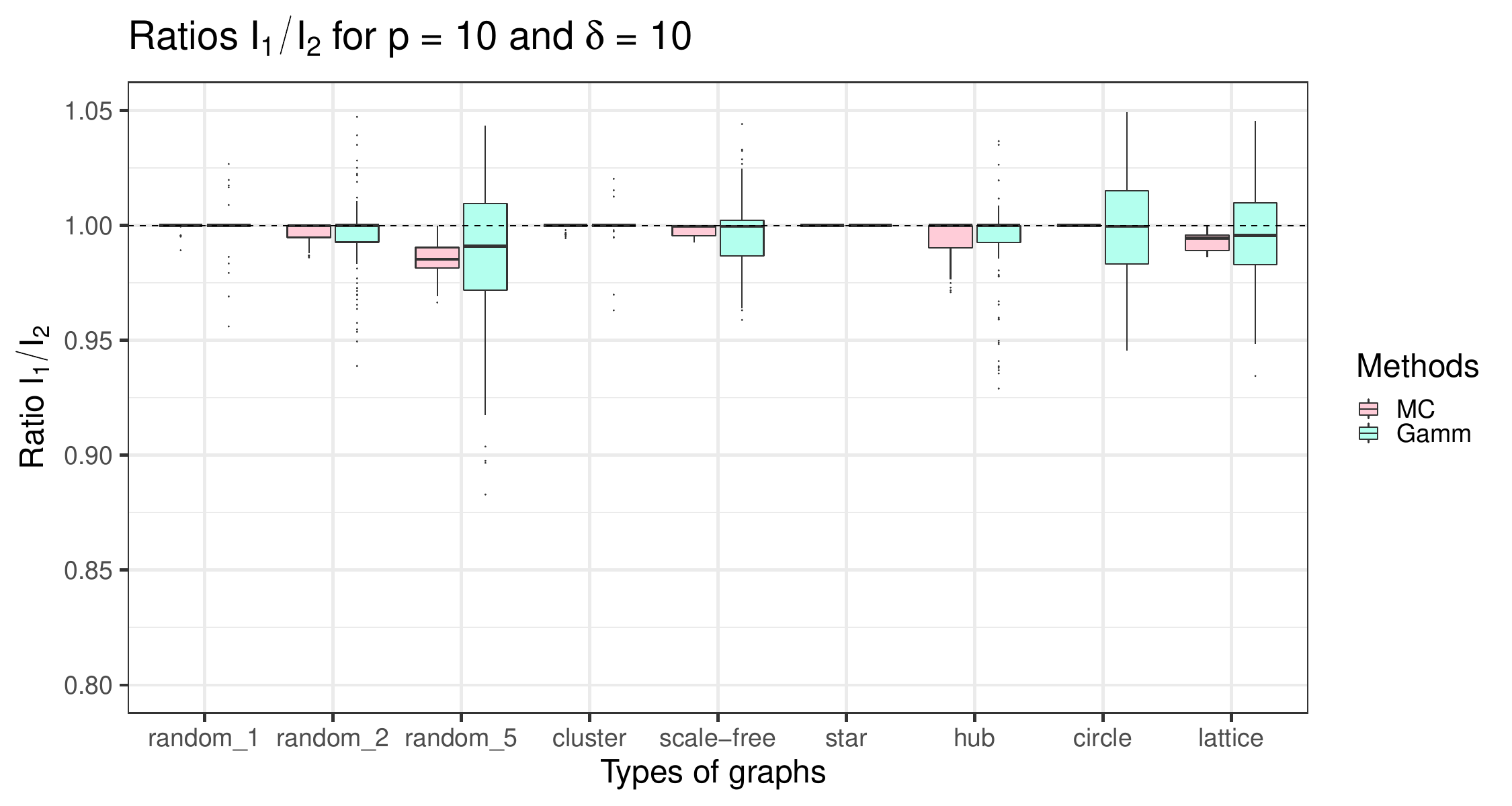}
    \includegraphics[width=15cm]{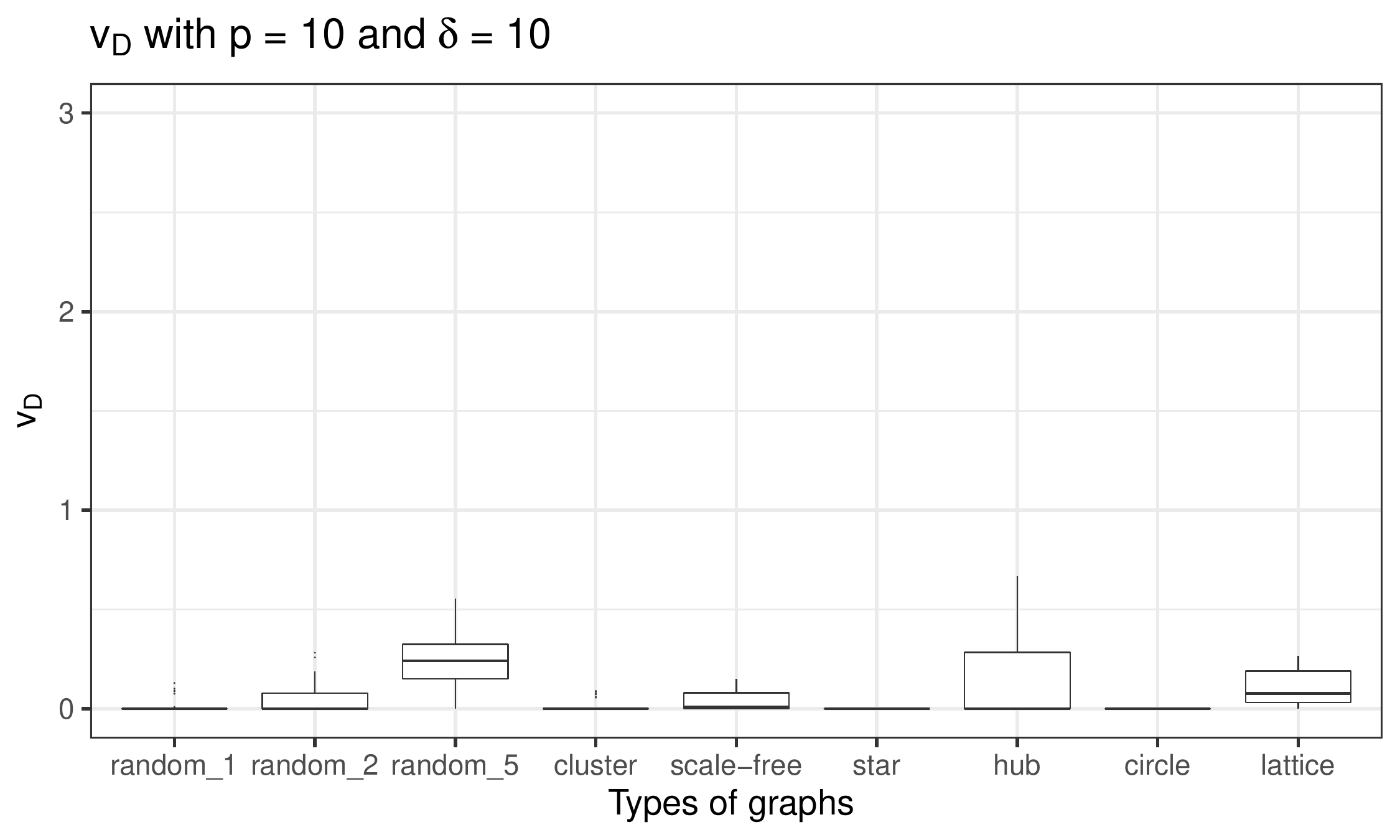}
\end{center}
\caption{ (Top) The boxplot for the ratio $I_1/I_2$ computed by the MC approach of \cite{atay2005monte} (in red) and our approximation \eqref{eq:ratioJoint} (in blue). (Bottom) The boxplot of the variance $v_D$ of $b_1$ for the corresponding graphs. These computations are done over 100 replications for nine different graphs (Figure \ref{fig:plot_graphs_p30}) with $p=10$ nodes and $\delta=10$. \label{fig:boxplot_p10_delta10} }
\end{figure}

\begin{figure} %[!ht]
\begin{center}
    \includegraphics[width=18cm]{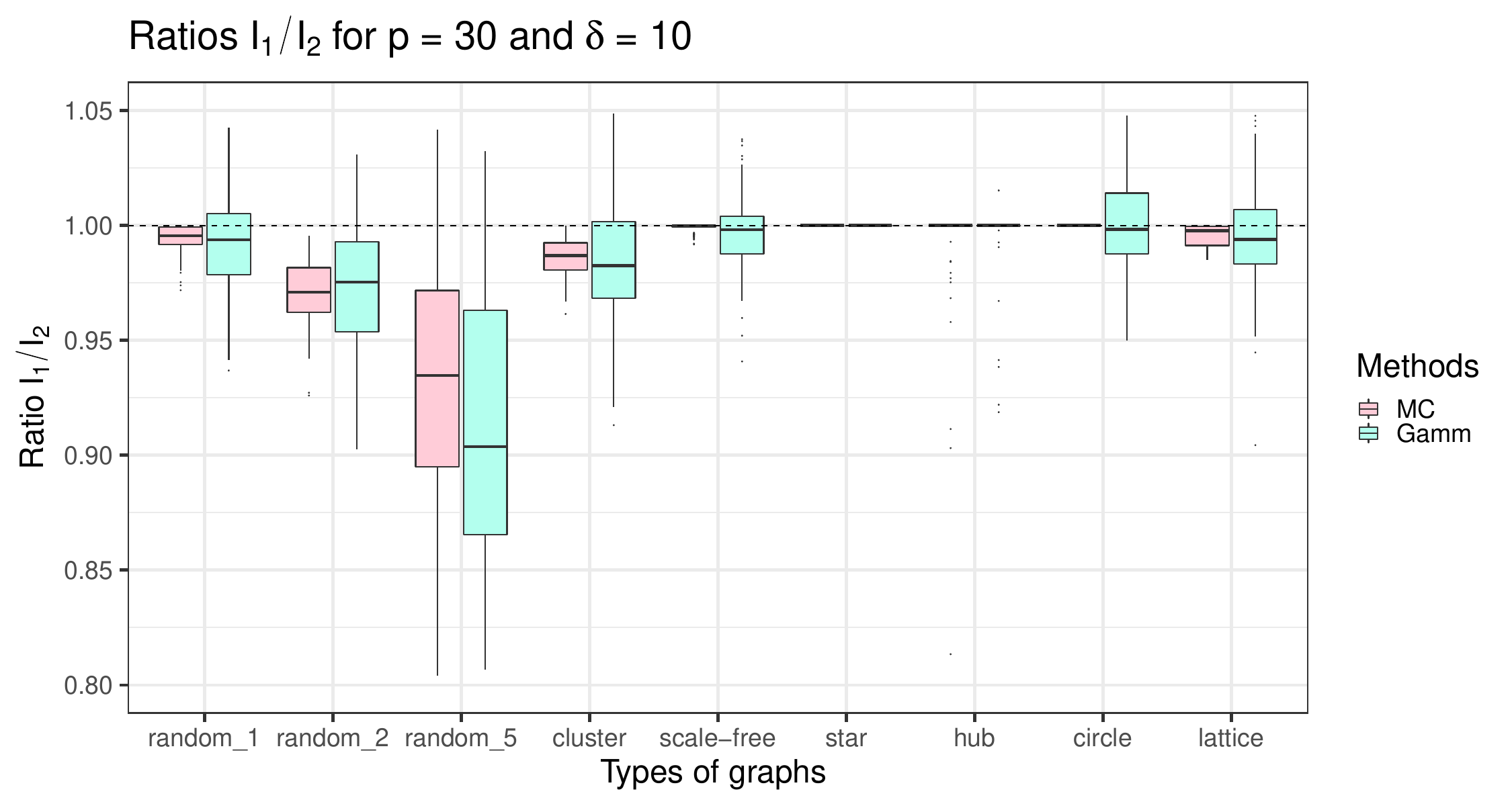}
    \includegraphics[width=15cm]{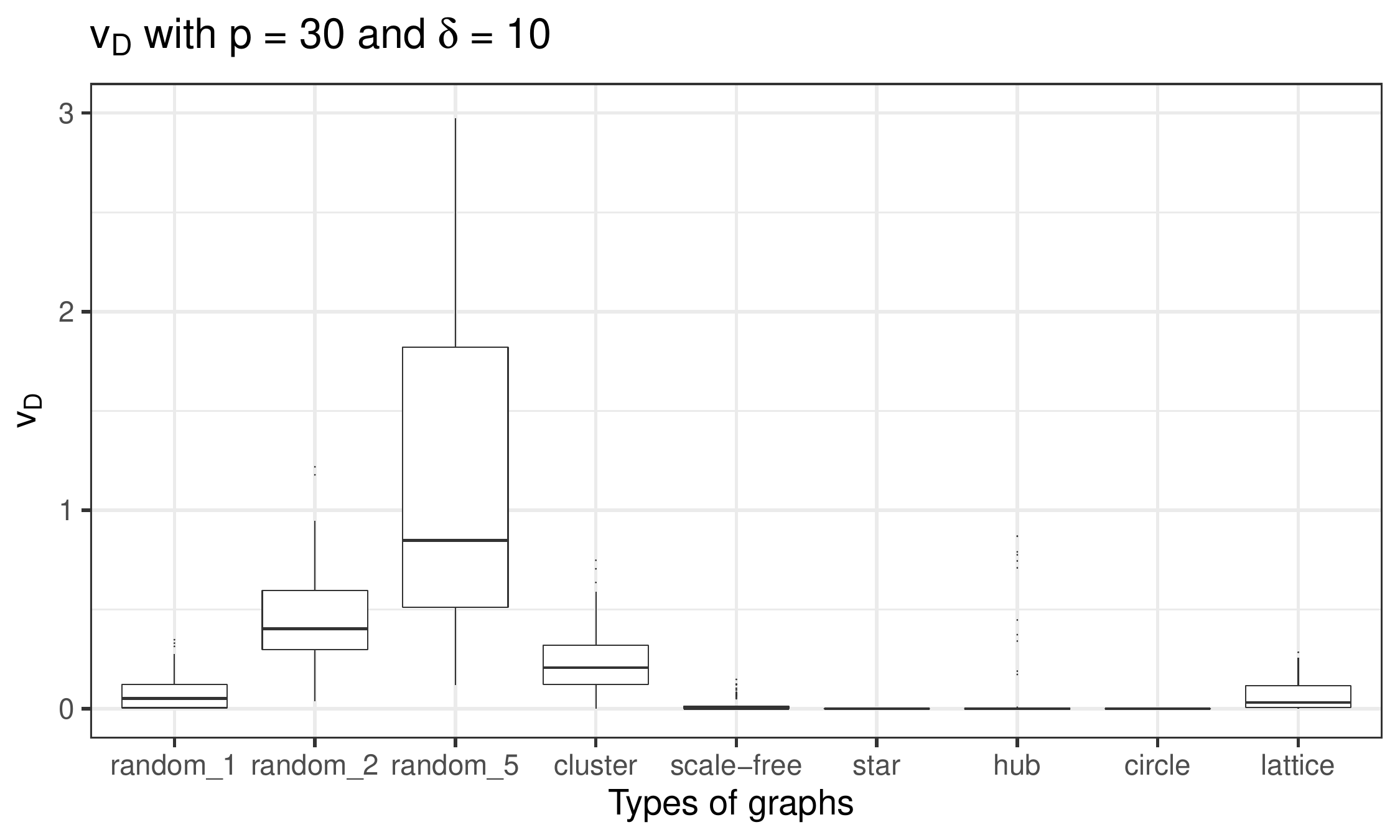}
\end{center}
\caption{ (Top) The boxplot for the ratio $I_1/I_2$ computed by the MC approach of \cite{atay2005monte} (in red) and our approximation \eqref{eq:ratioJoint} (in blue). (Bottom) The boxplot of the variance $v_D$ of $b_1$ for the corresponding graphs. These computations are done over 100 replications for nine different graphs (Figure \ref{fig:plot_graphs_p30}) with $p=30$ nodes and $\delta=10$. \label{fig:boxplot_p30_delta10} }
\end{figure}

\begin{figure} %[!ht]
\begin{center}
    \includegraphics[width=10cm]{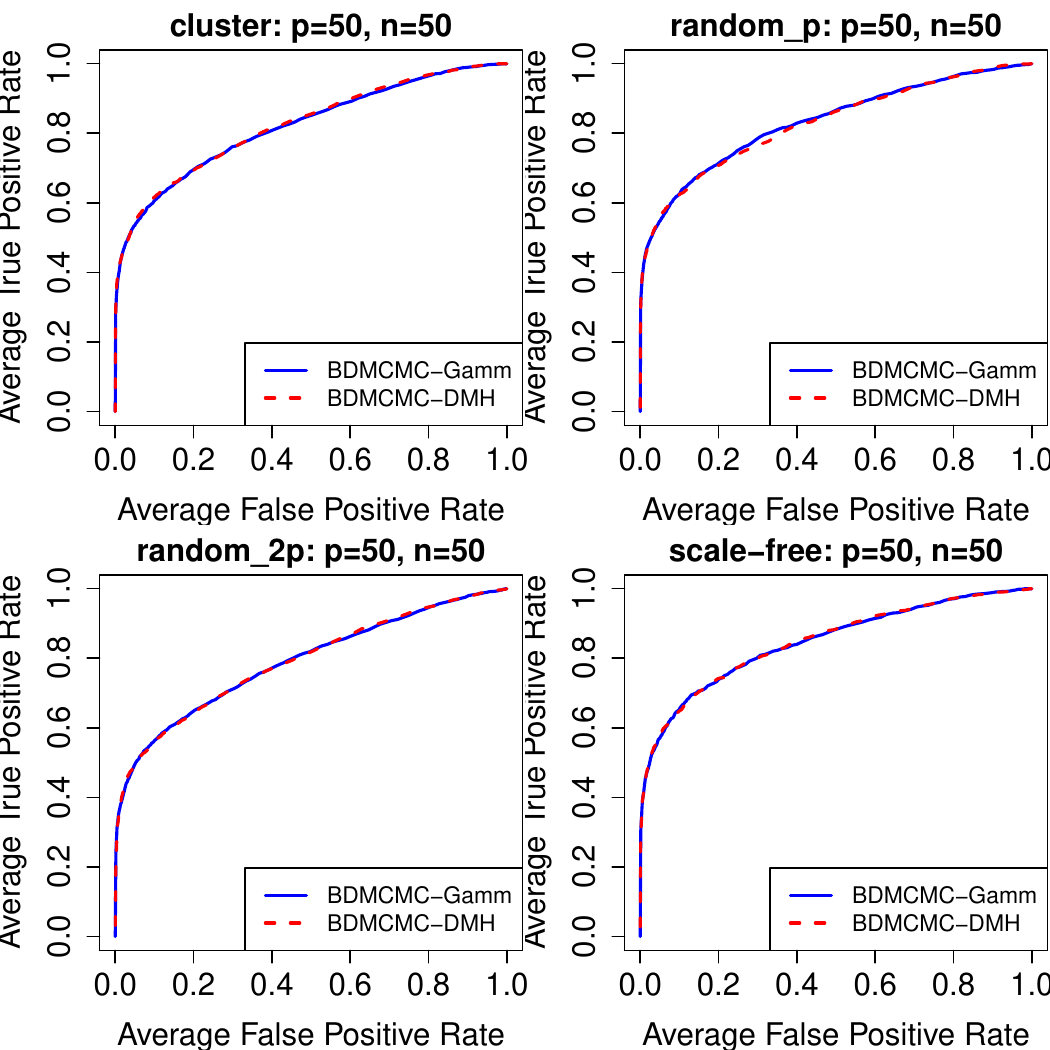} 
    \includegraphics[width=10cm]{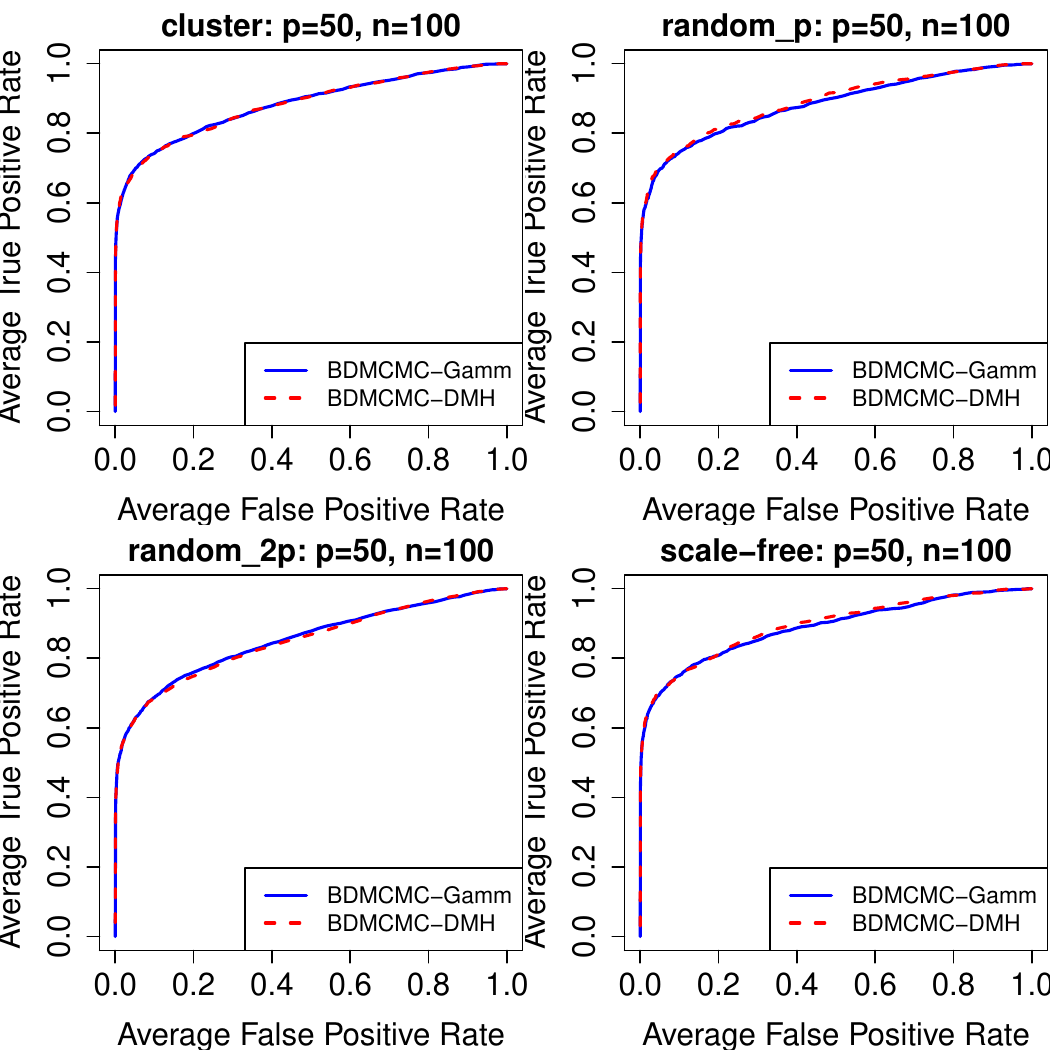}
\end{center}
\caption{ ROC curves for the BDMCMC algorithm with our approximation in Equation \ref{eq:approx} (BDMCMC-Gamm) and BDMCMC algorithm with exchange algorithm (BDMCMC-DMH), over $50$ replications. Here, $p=50$, $n \in \{50, 100 \}$, and $4$ different graph structures. \label{fig:Rocplot_p=50} }
\end{figure}

\begin{figure} %[!ht]
\begin{center}
    \includegraphics[width=10cm]{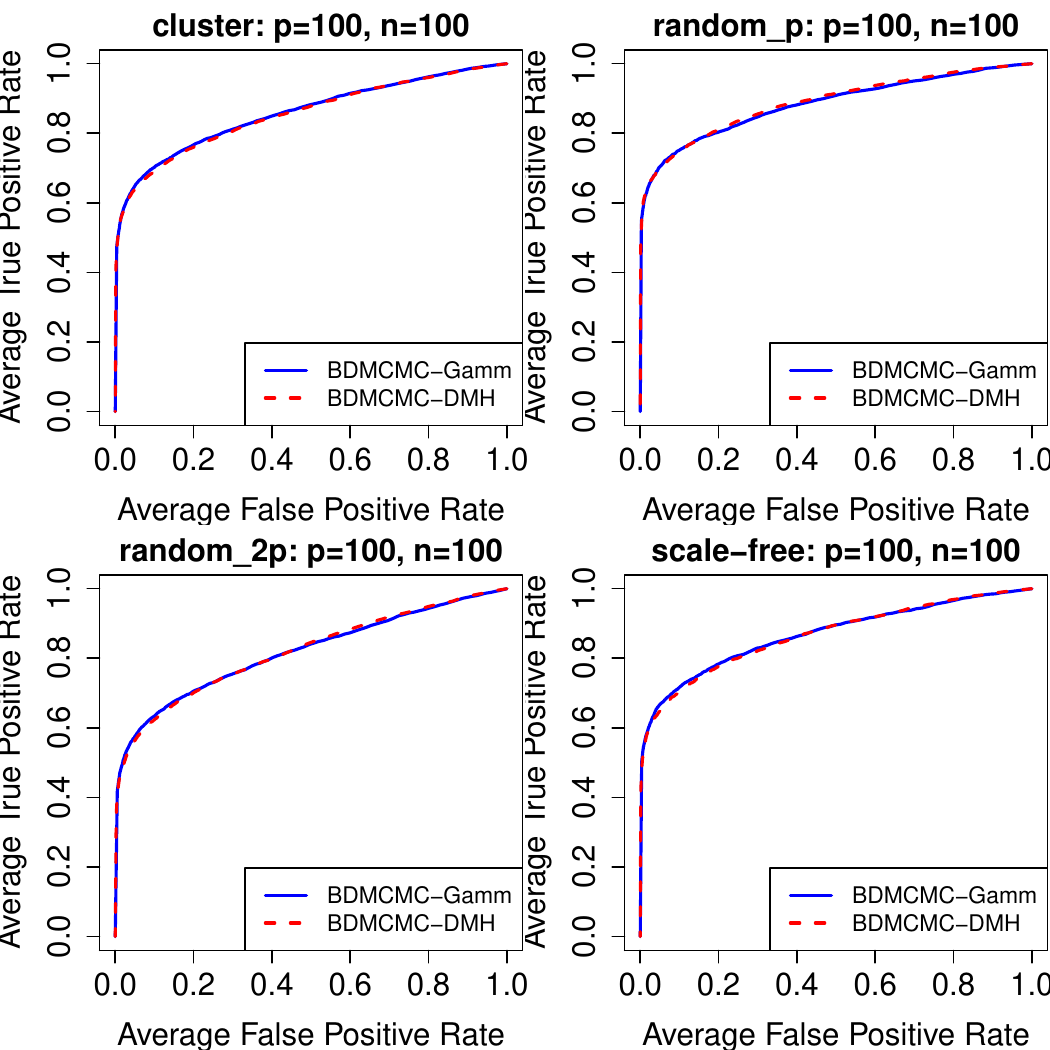} 
    \includegraphics[width=10cm]{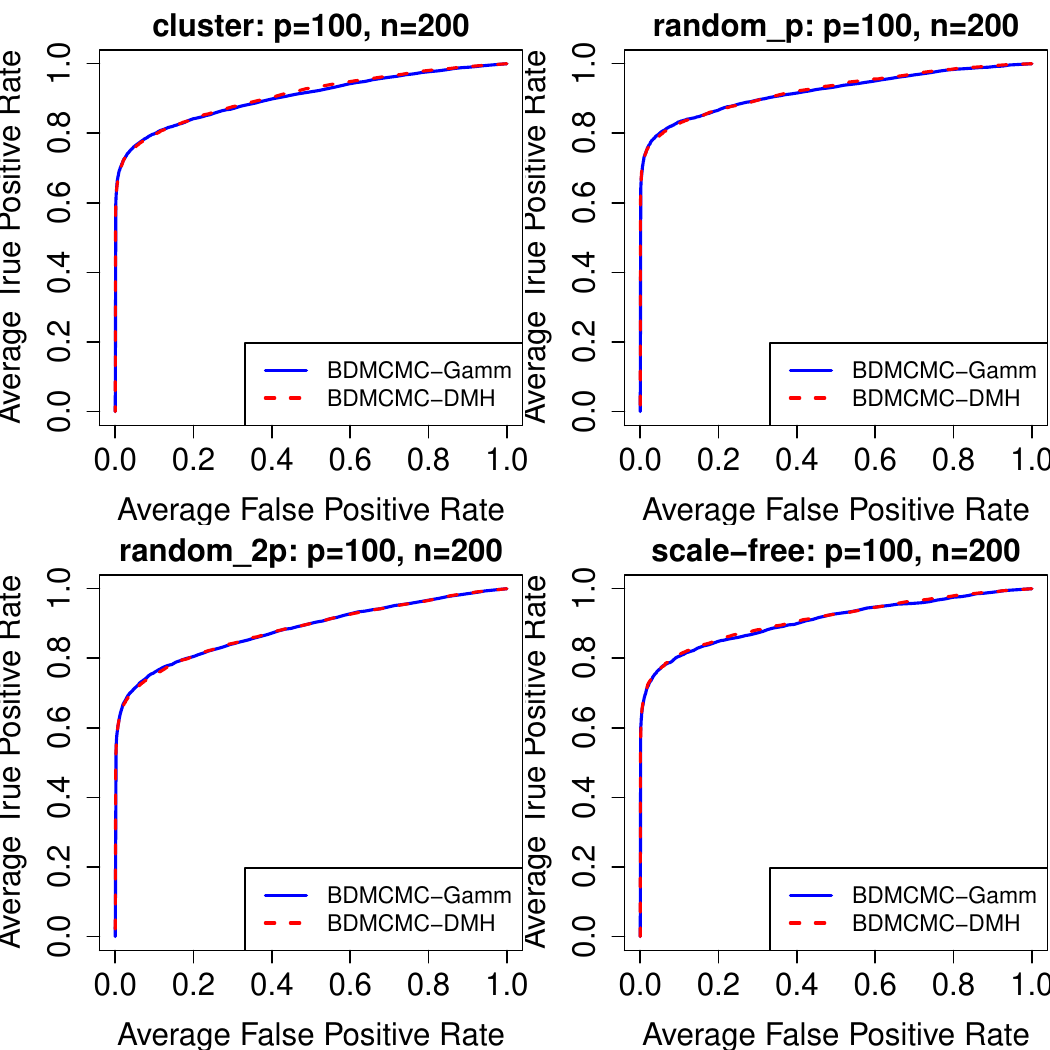}
\end{center}
\caption{ ROC curves for the BDMCMC algorithm with our approximation in Equation \ref{eq:approx} (BDMCMC-Gamm) and BDMCMC algorithm with exchange algorithm (BDMCMC-DMH), over $50$ replications. Here, $p=100$, $n \in \{100, 200 \}$, and $4$ different graph structures. \label{fig:Rocplot_p=100} }
\end{figure}

Figures \ref{fig:boxplot_p10}, \ref{fig:boxplot_p30}, \ref{fig:boxplot_p10_delta10}, and \ref{fig:boxplot_p30_delta10} are the results for the simulation in Section \ref{subsec:sim-joint} of the manuscript.
Figures \ref{fig:Rocplot_p=50} and \ref{fig:Rocplot_p=100} are the ROC plots for Section \ref{sec:simulation} of the manuscript.

% - - - - - - - - - - - - - - - - - - - - - - - - - - - - - - - - - - - - - - - - - - - - - - - - - - - - - - - - - - - -|

% - - - - - - - - - - - - - - - - - - - - - - - - - - - - - - - - - - - - - - - - - - - - - - - - - - - - - - - - - - - -|

% - - - - - - - - - - - - - - - - - - - - - - - - - - - - - - - - - - - - - - - - - - - - - - - - - - - - - - - - - - - -|

% - - - - - - - - - - - - - - - - - - - - - - - - - - - - - - - - - - - - - - - - - - - - - - - - - - - - - - - - - - - -|

% - - - - - - - - - - - - - - - - - - - - - - - - - - - - - - - - - - - - - - - - - - - - - - - - - - - - - - - - - - - -|
% ROC plots
% - - - - - - - - - - - - - - - - - - - - - - - - - - - - - - - - - - - - - - - - - - - - - - - - - - - - - - - - - - - -|

% - - - - - - - - - - - - - - - - - - - - - - - - - - - - - - - - - - - - - - - - - - - - - - - - - - - - - - - - - - - -|
%Plots for p = 100

% - - - - - - - - - - - - - - - - - - - - - - - - - - - - - - - - - - - - - - - - - - - - - - - - - - - - - - - - - - - -|
\bibliographystyle{Chicago}
\bibliography{ref}

\begin{thebibliography}{}

\bibitem[\protect\citeauthoryear{Abramovitz and Stegun}{Abramovitz and
  Stegun}{1972}]{abramovitzStegun}
Abramovitz, M. and I.~Stegun (1972).
\newblock {\em Handbook of Mathematical Functions}, Volume~55.
\newblock National Bureau of Standards, Applied Mathematics.

\bibitem[\protect\citeauthoryear{Albert and Barab{\'a}si}{Albert and
  Barab{\'a}si}{2002}]{albert2002statistical}
Albert, R. and A.~Barab{\'a}si (2002).
\newblock Statistical mechanics of complex networks.
\newblock {\em Reviews of modern physics\/}~{\em 74\/}(1), 47.

\bibitem[\protect\citeauthoryear{Atay-Kayis and Massam}{Atay-Kayis and
  Massam}{2005}]{atay2005monte}
Atay-Kayis, A. and H.~Massam (2005).
\newblock A {M}onte {C}arlo method for computing the marginal likelihood in
  nondecomposable {G}aussian graphical models.
\newblock {\em Biometrika\/}~{\em 92\/}(2), 317--335.

\bibitem[\protect\citeauthoryear{Chamayou and Letac}{Chamayou and
  Letac}{1991}]{chamayou1991explicit}
Chamayou, J.-F. and G.~Letac (1991).
\newblock Explicit stationary distributions for compositions of random
  functions and products of random matrices.
\newblock {\em Journal of Theoretical Probability\/}~{\em 4\/}(1), 3--36.

\bibitem[\protect\citeauthoryear{Cheng and Lenkoski}{Cheng and
  Lenkoski}{2012}]{cheng2012hierarchical}
Cheng, Y. and A.~Lenkoski (2012).
\newblock Hierarchical {G}aussian graphical models: Beyond reversible jump.
\newblock {\em Electronic Journal of Statistics\/}~{\em 6}, 2309--2331.

\bibitem[\protect\citeauthoryear{Dempster}{Dempster}{1972}]{dempster1972covariance}
Dempster, A.~P. (1972).
\newblock Covariance selection.
\newblock {\em Biometrics\/}, 157--175.

\bibitem[\protect\citeauthoryear{Dobra and Lenkoski}{Dobra and
  Lenkoski}{2011}]{dobra2011copula}
Dobra, A. and A.~Lenkoski (2011).
\newblock Copula {G}aussian graphical models and their application to modeling
  functional disability data.
\newblock {\em The Annals of Applied Statistics\/}~{\em 5\/}(2A), 969--993.

\bibitem[\protect\citeauthoryear{Dobra, Lenkoski, and Rodriguez}{Dobra
  et~al.}{2011}]{dobra2011bayesian}
Dobra, A., A.~Lenkoski, and A.~Rodriguez (2011).
\newblock Bayesian inference for general {G}aussian graphical models with
  application to multivariate lattice data.
\newblock {\em Journal of the American Statistical Association\/}~{\em
  106\/}(496), 1418--1433.

\bibitem[\protect\citeauthoryear{Friedman, Hastie, and Tibshirani}{Friedman
  et~al.}{2008}]{friedman2008sparse}
Friedman, J., T.~Hastie, and R.~Tibshirani (2008).
\newblock Sparse inverse covariance estimation with the graphical lasso.
\newblock {\em Biostatistics\/}~{\em 9\/}(3), 432--441.

\bibitem[\protect\citeauthoryear{Green}{Green}{1995}]{green1995reversible}
Green, P. (1995).
\newblock Reversible jump markov chain {M}onte carlo computation and {B}ayesian
  model determination.
\newblock {\em Biometrika\/}~{\em 82\/}(4), 711--732.

\bibitem[\protect\citeauthoryear{Green}{Green}{2003}]{green2003trans}
Green, P. (2003).
\newblock Trans-dimensional {M}arkov chain {M}onte {C}arlo.
\newblock {\em Oxford Statistical Science Series\/}, 179--198.

\bibitem[\protect\citeauthoryear{Hinne, Lenkoski, Heskes, and van Gerven}{Hinne
  et~al.}{2014}]{hinne2014efficient}
Hinne, M., A.~Lenkoski, T.~Heskes, and M.~van Gerven (2014).
\newblock Efficient sampling of gaussian graphical models using conditional
  bayes factors.
\newblock {\em Stat\/}~{\em 3\/}(1), 326--336.

\bibitem[\protect\citeauthoryear{Hinoveanu, Leisen, and Villa}{Hinoveanu
  et~al.}{2018}]{hinoveanu2018loss}
Hinoveanu, L.~C., F.~Leisen, and C.~Villa (2018).
\newblock A loss-based prior for gaussian graphical models.
\newblock {\em arXiv preprint arXiv:1812.05531\/}.

\bibitem[\protect\citeauthoryear{Lauritzen}{Lauritzen}{1996}]{lauritzen1996graphical}
Lauritzen, S. (1996).
\newblock {\em Graphical Models}, Volume~17.
\newblock Oxford University Press, USA.

\bibitem[\protect\citeauthoryear{Lenkoski}{Lenkoski}{2013}]{lenkoski2013direct}
Lenkoski, A. (2013).
\newblock A direct sampler for {G}-{W}ishart variates.
\newblock {\em Stat\/}~{\em 2\/}(1), 119--128.

\bibitem[\protect\citeauthoryear{Lenkoski and Dobra}{Lenkoski and
  Dobra}{2011}]{lenkoski2011computational}
Lenkoski, A. and A.~Dobra (2011).
\newblock Computational aspects related to inference in {G}aussian graphical
  models with the g-wishart prior.
\newblock {\em Journal of Computational and Graphical Statistics\/}~{\em
  20\/}(1), 140--157.

\bibitem[\protect\citeauthoryear{Letac and Massam}{Letac and
  Massam}{2020}]{letac2020gaussian}
Letac, G. and H.~Massam (2020).
\newblock Gaussian approximation of gaussian scale mixture.
\newblock {\em kybersetika\/}~{\em 56\/}(6), 1063--1080.

\bibitem[\protect\citeauthoryear{Letac, Massam, et~al.}{Letac
  et~al.}{2007}]{letac2007wishart}
Letac, G., H.~Massam, et~al. (2007).
\newblock Wishart distributions for decomposable graphs.
\newblock {\em The Annals of Statistics\/}~{\em 35\/}(3), 1278--1323.

\bibitem[\protect\citeauthoryear{Liang}{Liang}{2010}]{liang2010double}
Liang, F. (2010).
\newblock A double {M}etropolis--{H}astings sampler for spatial models with
  intractable normalizing constants.
\newblock {\em Journal of Statistical Computation and Simulation\/}~{\em
  80\/}(9), 1007--1022.

\bibitem[\protect\citeauthoryear{Mohammadi and Wit}{Mohammadi and
  Wit}{2015}]{mohammadi2015bayesianStructure}
Mohammadi, A. and E.~Wit (2015).
\newblock {B}ayesian structure learning in sparse {G}aussian graphical models.
\newblock {\em Bayesian Analysis\/}~{\em 10\/}(1), 109--138.

\bibitem[\protect\citeauthoryear{Mohammadi and Wit}{Mohammadi and
  Wit}{2019a}]{BDgraph}
Mohammadi, R. and E.~Wit (2019a).
\newblock {\em BDgraph: Bayesian Structure Learning in Graphical Models using
  Birth-Death MCMC}.
\newblock R package version 2.62.

\bibitem[\protect\citeauthoryear{Mohammadi and Wit}{Mohammadi and
  Wit}{2019b}]{mohammadi2019bdgraph}
Mohammadi, R. and E.~C. Wit (2019b).
\newblock {BDgraph}: An {R} package for {B}ayesian structure learning in
  graphical models.
\newblock {\em Journal of Statistical Software\/}~{\em 89\/}(3), 1--30.

\bibitem[\protect\citeauthoryear{Murray, Ghahramani, and MacKay}{Murray
  et~al.}{2006}]{murray2006mcmc}
Murray, I., Z.~Ghahramani, and D.~MacKay (2006).
\newblock {MCMC} for doubly-intractable distributions.
\newblock {\em Proceedings of the 22nd Annual Conference on Uncertainty in
  Artificial Intelligence\/}, 359--366.

\bibitem[\protect\citeauthoryear{Park and Haran}{Park and
  Haran}{2018}]{park2018bayesian}
Park, J. and M.~Haran (2018).
\newblock Bayesian inference in the presence of intractable normalizing
  functions.
\newblock {\em Journal of the American Statistical Association\/}~{\em
  113\/}(523), 1372--1390.

\bibitem[\protect\citeauthoryear{Roverato}{Roverato}{2002}]{roverato2002hyper}
Roverato, A. (2002).
\newblock Hyper inverse {W}ishart distribution for non-decomposable graphs and
  its application to {B}ayesian inference for {G}aussian graphical models.
\newblock {\em Scandinavian Journal of Statistics\/}~{\em 29\/}(3), 391--411.

\bibitem[\protect\citeauthoryear{Rue and Held}{Rue and
  Held}{2005}]{rue2005gaussian}
Rue, H. and L.~Held (2005).
\newblock {\em Gaussian Markov random fields: theory and applications}.
\newblock CRC press.

\bibitem[\protect\citeauthoryear{Uhler, Lenkoski, and Richards}{Uhler
  et~al.}{2018}]{uhler2014exact}
Uhler, C., A.~Lenkoski, and D.~Richards (2018).
\newblock Exact formulas for the normalizing constants of {W}ishart
  distributions for graphical models.
\newblock {\em The Annals of Statistics\/}~{\em 46}, 90--118.

\bibitem[\protect\citeauthoryear{Wang}{Wang}{2012}]{wang12}
Wang, H. (2012).
\newblock The {B}ayesian graphical {L}asso and efficient posterior computation.
\newblock {\em Bayesian Analysis\/}~{\em 7}, 771--790.

\bibitem[\protect\citeauthoryear{Wang and Li}{Wang and
  Li}{2012}]{wang2012efficient}
Wang, H. and S.~Li (2012).
\newblock Efficient {G}aussian graphical model determination under
  {G}-{W}ishart prior distributions.
\newblock {\em Electronic Journal of Statistics\/}~{\em 6}, 168--198.

\end{thebibliography}
% - - - - - - - - - - - - - - - - - - - - - - - - - - - - - - - - - - - - - - - - - - - - - - - - - - - - - - - - - - - -|
\end{document}